%% file: HypCuspDual.tex
\documentclass[10pt]{amsart}

\usepackage[T1]{fontenc}
\usepackage[dvips]{graphicx}
\usepackage{epsfig}
\usepackage{amsmath,amssymb,amscd,amsthm}

\usepackage[latin1]{inputenc}
\usepackage{latexsym}
\usepackage{graphics}
\usepackage{colortbl}
\usepackage{color}
\usepackage{ae}
\usepackage{url}

\graphicspath{{Fig/}}
\def\H{{\mathbb H}}

\usepackage{textcomp}
\begin{document}

\setcounter{secnumdepth}{3}
\setcounter{tocdepth}{2}

\newtheorem{dfn}{Definition}[section]
\newtheorem{thm}[dfn]{Theorem}
\newtheorem{prp}[dfn]{Proposition}
\newtheorem{lem}[dfn]{Lemma}
\newtheorem{cor}[dfn]{Corollary}

\newcommand{\mf}{\mathfrak}
\newcommand{\mb}{\mathbb}
\newcommand{\ol}{\overline}
\newcommand{\la}{\langle}
\newcommand{\ra}{\rangle}

\newtheorem{thmprime}{Theorem}
\renewcommand{\thethmprime}{1.\arabic{thmprime}\textquotesingle}


\newtheorem{Alphatheorem}{Theorem}
\renewcommand{\theAlphatheorem}{\Alph{Alphatheorem}}

\newtheorem{Alphatheoremprime}{Theorem}
\renewcommand{\theAlphatheoremprime}{\Alph{Alphatheoremprime}\textquotesingle}

\theoremstyle{definition}

\newtheorem{rem}[dfn]{Remark}
\newtheorem{exl}[dfn]{Example}

\theoremstyle{plain}
\newtheorem*{flp}{Flip algorithm}
\newtheorem{clm}{Claim}

\newcommand{\EM}{\ensuremath}
\newcommand{\norm}[1]{\EM{\left\| #1 \right\|}}

\newcommand{\modul}[1]{\left| #1\right|}

\def\co{\colon\thinspace}

\def\N{{\mathbb N}}
\def\R{{\mathbb R}}
\def\Z{{\mathbb Z}}
\def\Sph{{\mathbb S}}
\def\Tor{{\mathbb T}}
\def\Disk{{\mathbb D}}
\def\Fl{{\mathbb F}}

\def\H{{\mathbb H}}
\def\RP{{\mathbb R}{\mathrm{P}}}
\def\dS{{\mathrm d}{\mathbb{S}}}

\def\phi{\varphi}
\def\epsilon{\varepsilon}
\def\V{{\mathcal V}}
\def\E{{\mathcal E}}
\def\F{{\mathcal F}}
\def\C{{\mathcal C}}

\renewcommand{\span}{\operatorname{span}}
\newcommand{\ang}{\operatorname{ang}}
\newcommand{\area}{\operatorname{area}}
\newcommand{\cone}{\operatorname{cone}}
\newcommand{\pr}{\operatorname{pr}}
\newcommand{\vol}{\operatorname{vol}}
\newcommand{\covol}{\operatorname{covol}}
\newcommand{\st}{\operatorname{st}}
\newcommand{\ost}{\operatorname{st}^\circ}
\newcommand{\inn}{\operatorname{int}}
\newcommand{\grad}{\operatorname{grad}}
\def\dist{\mathrm{dist\,}}
\def\diam{\mathrm{diam\,}}

\def\M{{\mathcal M}}
\def\Mt{{\mathcal M}_{\mathrm{tr}}}
\def\T{{\mathcal T}}
\def\Vt{V_{\mathrm{tr}}}

\title[Gauss images of convex polyhedral cusps]{Gauss images of hyperbolic cusps\\ with convex polyhedral boundary}
\author{Fran\c{c}ois Fillastre}
\address{University of Cergy-Pontoise\\ UMR CNRS 8088\\ departement of mathematics\\
F-95000 Cergy-Pontoise\\ FRANCE}
\email{francois.fillastre@u-cergy.fr}
\thanks{The first author was partially supported by Schweizerischer Nationalfonds 200020-113199/1} 

 \author{Ivan Izmestiev}
\thanks{The second author was supported by the DFG Research Unit 565 ``Polyhedral Surfaces''}
\address{Institut f\"ur Mathematik, MA 8-3 \\
Technische Universit\"at Berlin \\
Stra{\ss}e des 17. Juni 136 \\
D-10623 Berlin \\
 GERMANY}
\email{izmestiev@math.tu-berlin.de}

\date{August 13, 2009}
\maketitle

\begin{abstract}
We prove that a 3--dimensional hyperbolic cusp with convex polyhedral boundary is uniquely determined by its Gauss image. Furthermore, any spherical metric on the torus with cone singularities of negative curvature and all closed contractible geodesics of length greater than $2\pi$ is the metric of the Gauss image of some convex polyhedral cusp. This result is an analog of the Rivin-Hodgson theorem characterizing compact convex hyperbolic polyhedra in terms of their Gauss images.

The proof uses a variational method. Namely, a cusp with a given Gauss image is identified with a critical point of a functional on the space of cusps with cone-type singularities along a family of half-lines. The functional is shown to be concave and to attain maximum at an interior point of its domain. As a byproduct, we prove rigidity statements with respect to the Gauss image for cusps with or without cone-type singularities.

In a special case, our theorem is equivalent to existence of a circle pattern on the torus, with prescribed combinatorics and intersection angles. This is the genus one case of a theorem by Thurston. In fact, our theorem extends Thurston's theorem in the same way as Rivin-Hodgson's theorem extends Andreev's theorem on compact convex polyhedra with non-obtuse dihedral angles.

The functional used in the proof is the sum of a volume term and curvature term. We show that, in the situation of Thurston's theorem, it is the potential for the combinatorial Ricci flow considered by Chow and Luo.

Our theorem represents the last special case of a general statement about isometric immersions of compact surfaces.
\end{abstract}


\section{Introduction}

\subsection{Rivin-Hodgson and Andreev theorems}
\label{subsec:RHAnd}
In \cite{Riv86} and \cite{RivHod93}, the following theorem was proved.
\begin{thm}[Rivin-Hodgson]
\label{thm:RivHod}
Let $g$ be a spherical cone metric on the sphere $\Sph^2$ such that the following two conditions hold:
\begin{enumerate}
\item All cone angles of $(\Sph^2, g)$ are greater than $2\pi$.
\item All closed geodesics on $(\Sph^2, g)$ have lengths greater than $2\pi$.
\end{enumerate}
Then there exists a compact convex hyperbolic polyhedron $P$ such that the Gauss image of $\partial P$ is isometric to $(\Sph^2, g)$. Furthermore, $P$ is unique up to isometry.
\end{thm}

A \emph{spherical cone metric} is locally modelled on the standard spherical metric of curvature $1$, away from a finite set of \emph{cone points}. A cone point is said to have \emph{cone angle} $\theta$ if it has a neighborhood isometric to an angle of size $\theta \ne 2\pi$ with identified sides. For example, gluing $\Sph^2$ from a collection of spherical polygons defines a spherical cone metric with cone points at those vertices where the sum of the adjacent angles is different from $2\pi$.

Now let $P$ be a compact convex hyperbolic polyhedron. To define the \emph{Gauss image} of $\partial P$, one first defines the Gauss image of a vertex $v$ of $P$ as a spherical polygon polar dual to the link of $v$ in $P$. Then, the Gauss images of all vertices of $P$ are glued side-to-side to form the Gauss image of $\partial P$. For more details, see Section \ref{subsec:GaussImage}.

In \cite{RivHod93}, it is also shown that the Gauss image of the boundary of a compact convex hyperbolic polyhedron satisfies the two conditions of Theorem \ref{thm:RivHod}. Thus, Theorem \ref{thm:RivHod} characterizes compact convex hyperbolic polyhedra in terms of the intrinsic metrics of the Gauss images of their boundaries.

As shown in \cite{Hod92}, Theorem \ref{thm:RivHod} implies Andreev's theorem \cite{And70} on compact hyperbolic polyhedra with non-obtuse dihedral angles. For brevity, we cite here the version for acute dihedral angles.

\begin{thm}[Andreev]
\label{thm:And}
Let $C$ be a cellular subdivision of $\Sph^2$ combinatorially equivalent to a compact convex polyhedron with trivalent vertices other than the tetrahedron. Call three edges $e$, $f$, $g$ of $C$ a \emph{proper cutset}, if they are pairwise disjoint and if there exists a simple closed curve on $\Sph^2$ that intersects each of them exactly once.

Let $\phi \co e \to \phi_e$ be a map from the edge set of $C$ to the interval $(0, \frac{\pi}2)$ such that the following two conditions hold.
\begin{enumerate}
\item \label{it:1A} If $e$, $f$, $g$ are three edges incident to a vertex, then $\phi_e + \phi_f + \phi_g > \pi$.
\item \label{it:2A} If $e$, $f$, $g$ form a proper cutset, then $\phi_e + \phi_f + \phi_g < \pi$.
\end{enumerate}
Then there exists a compact convex hyperbolic polyhedron combinatorially equivalent to $C$ and with dihedral angles $\phi$.
\end{thm}
To understand the relation between Theorems \ref{thm:RivHod} and \ref{thm:And}, note that the edge lengths of the Gauss image of a vertex $v$ are determined by the values of dihedral angles at the edges incident to $v$. Since $C$ in Theorem \ref{thm:And} is trivalent, its dual $C^*$ is a triangulation of $\Sph^2$, and the map $\phi$ determines the lengths of the edges of this triangulation. Theorem \ref{thm:And} can be proved by applying Theorem~\ref{thm:RivHod} to the spherical cone metric obtained in this way.

Note that in Theorem \ref{thm:RivHod} the combinatorics of the polyhedron $P$ is not given, and cannot a priori be derived from the metric $g$.

\subsection{Analogy with the Minkowski theorem}
\label{subsec:MinkThm}
Cone points of the Gauss image of $\partial P$ naturally correspond to the faces of $P$. Besides, for a point with cone angle~$\theta_i$, the area of the corresponding face is $\theta_i - 2\pi$. Thus the metric of the Gauss image determines the number of faces of $P$ and their areas in a straightforward way.

Let us see what happens in the Euclidean limit. Consider a family $P_t$ of convex hyperbolic polyhedra shrinking to a point so that their rescaled limit is a convex Euclidean polyhedron $Q$. Let $g_t$ be the metric of the Gauss image of $\partial P_t$. It is easy to see that $g_t$ converges to the metric of the Gauss image of $Q$, which is just the standard metric on the unit sphere. Thus, from the first sight, in the limit all information is lost.

However, the cone points of the metrics $g_t$ converge to a collection of points on $\Sph^2 \subset \R^3$, namely to the unit normals of the faces of $Q$. Besides, although the angle excesses of the cone points of $g_t$ tend to zero, their ratios tend to the ratios of the face areas of $Q$. Thus, in the Euclidean limit, we are left with the directions of the face normals of a Euclidean polyhedron and with the face areas up to a common factor. There is a classical theorem saying that the polyhedron can be recovered from this information.

\begin{thm}[Minkowski]
\label{thm:Mink}
Let $\nu_1, \ldots, \nu_n$ be non-coplanar unit vectors in~$\R^3$, and let $F_1, \ldots, F_n$ be positive real numbers such that
\begin{equation}
\label{eqn:Mink}
\sum_{i=1}^n F_i\nu_i = 0.
\end{equation}
Then there exists a compact convex Euclidean polyhedron $Q$ with outer face normals~$\nu_i$ and face areas $F_i$. Besides, $Q$ is unique up to isometry.
\end{thm}

Note that the equation \eqref{eqn:Mink} is always fulfilled when $\nu$ and $F$ arise as a limit of spherical cone metrics, see \cite[Lemma 4.12]{BI08}.

Thus Theorem \ref{thm:RivHod} can be viewed as a hyperbolic analog of Theorem \ref{thm:Mink} (there are also other analogs, see \cite[Section 7.6.4]{Ale05}). However, there is a substantial difference between the proofs of Theorems \ref{thm:RivHod} and \ref{thm:Mink}. The idea in \cite{RivHod93} is to consider the space $\mathcal{P}_n$ of compact convex polyhedra with $n$ faces, the space $\mathcal{M}_n$ of spherical cone metrics with $n$ cone points satisfying the conditions of the theorem, and the map $\Gamma \co \mathcal{P}_n \to \mathcal{M}_n$ associating to a polyhedron the Gauss image of its boundary. Then some topological properties of the map $\Gamma$ are proved that imply that $\Gamma$ is a homeomorphism. We call this kind of argument Alexandrov's method, as it was extensively used by A.~D.~Alexandrov. While Theorem \ref{thm:Mink} can also be proved by Alexandrov's method, see \cite[Section 7.1]{Ale05}, more often it is proved by a variational method as outlined in the next paragraph.

Let $\mathcal{Q}(\nu)$ be the space of all convex Euclidean polyhedra with outer face normals $\nu_1, \ldots, \nu_n$. A polyhedron from $\mathcal{Q}(\nu)$ is uniquely determined by its heights $h_1, \ldots, h_n$, where $h_i$ is the (signed) distance of the plane of the $i$-th face from the origin. Let
$$
\vol \co \mathcal{Q}(\nu) \to \R
$$
be the function that associates to a polyhedron $Q(h)$ its volume. Then we have
\begin{equation}
\label{eqn:MinkDiff}
\frac{\partial \vol(Q(h))}{\partial h_i} = \area(Q_i(h)),
\end{equation}
where $Q_i(h)$ is the $i$-th face of the polyhedron $Q(h)$. Thus, every critical point of the function
\begin{equation}
\label{eqn:Restr}
\vol \co \mathcal{Q}(\nu) \cap \Big\{\sum_{i=1}^n h_iF_i = 1\Big\} \to \R
\end{equation}
satisfies $\area(Q_i(h)) = \mathrm{const} \cdot F_i$ and provides a solution to the problem, up to a scaling. One can show that the domain of \eqref{eqn:Restr} is convex and that the function is concave and achieves its maximum in the interior. Thus the solution exists and is unique.

Unlike the Alexandrov method, the variational method is constructive, as it amounts to finding a critical point of a functional. Also, the functional appearing in a variational method usually has a geometric interpretation, which throws more light on the problem and can lead to further insights.

\subsection{Results of the present paper}
\label{subsec:PresPap}
In this paper, we adapt the variational approach used in the proof of Theorem \ref{thm:Mink} to prove an analog of Theorem \ref{thm:RivHod} for spherical cone metrics on the torus.

\begin{Alphatheorem}
\label{thm:Main}
Let $g$ be a spherical cone metric on the torus $\Tor^2$ such that the following two conditions hold:
\begin{enumerate}
\item All cone angles of $(\Tor^2, g)$ are greater than $2\pi$.
\item All contractible closed geodesics on $(\Tor^2, g)$ have lengths greater than $2\pi$.
\end{enumerate}
Then there exists a convex polyhedral cusp $M$ such that the Gauss image of $\partial M$ is isometric to $(\Tor^2, g)$. Furthermore, $M$ is unique up to isometry.
\end{Alphatheorem}

A convex polyhedral cusp is a complete hyperbolic manifold of finite volume homeomorphic to $\Tor^2 \times [0, +\infty)$ and with convex polyhedral boundary, see Section~\ref{subsec:GaussImage}. In the same spirit, a compact convex hyperbolic polyhedron is a hyperbolic manifold homeomorphic to a $3$--ball and with convex polyhedral boundary.

To prove Theorem \ref{thm:Main}, we introduce the notion of a convex polyhedral cusp with coparticles, which is basically a convex polyhedral cusp with singular locus a union of half-lines orthogonal to the boundary. For a convex polyhedral cusp with coparticles, the Gauss image still can be defined. The space $\M^*(g)$ of all cusps with coparticles with the Gauss image $(\Tor^2, g)$ is an analog of the space $\mathcal{Q}(\nu)$ from Section \ref{subsec:MinkThm}. Similarly to a polyhedron from $\mathcal{Q}(\nu)$, a cusp from $\M^*(g)$ is also determined by its collection of truncated heights $\{h_i\}$. We define a function $V$ on $\M^*(g)$ through the formula
$$
V(M(h)) = -2\vol(M(h)) + \sum h_i\kappa_i,
$$
where $\kappa_i$ is the curvature (angle deficit) of the $i$-th coparticle. One can show that
\begin{equation}
\label{eqn:VarV}
\frac{\partial V(M(h))}{\partial h_i} = \kappa_i,
\end{equation}
so that the critical points of $V$ are cusps without coparticles. We show that the domain $\M^*(g)$ is contractible (although not convex) and that the function $V$ is concave. From an analysis of the behavior of $V$ at the boundary of $\M^*(g)$ and at the infinity we conclude that $V$ has a unique critical point in the interior of $\M^*(g)$. Thus a cusp with a given Gauss image exists and is unique.

Because of \eqref{eqn:VarV}, the Hessian of the function $V$ is the Jacobian of the map $h \mapsto \kappa$. Therefore $V$ can be used to study infinitesimal rigidity of cusps with coparticles.

\begin{Alphatheorem}
\label{thm:LocRig}
Convex polyhedral cusps are infinitesimally rigid with respect to their Gauss images. Convex polyhedral cusps with coparticles are locally rigid.
\end{Alphatheorem}

Infinitesimal rigidity of convex polyhedral cusps with respect to their Gauss images is also proved in \cite{Fill3}, with a different method.


It is conceivable that Theorem~\ref{thm:Main} can also be proved by an adaptation of the method used by Rivin and Hodgson. The proof would be close to those in \cite{Schpoly,Fillastre2}. However, one of the steps would be to prove the local rigidity of convex polyhedral cusps with respect to their Gauss images. The only way we can do that is by showing non-degeneracy of a matrix which is in fact the Hessian of our function $V$ or of the function from \cite{FI09}. See also \cite{Fill3}.

In the opposite direction, a variational proof of Theorem \ref{thm:RivHod} is not straigthforward. In place of cusps with coparticles one would consider cone manifolds homeomorphic to the $3$--ball whose singular locus is a union of segments with a common endpoint in the interior, the other endpoints in the faces of the polyhedron. An analog of function $V$ is not hard to find. The problem is that $V$ will almost surely have signature $(1, n-1)$. In the proof of the Minkowski Theorem, this problem is solved by restricting the volume function to a hyperplane. One could try to do the same in the Rivin-Hodgson theorem. Alternatively, one could try to consider polyhedra with singular locus a union of segments ending at a boundary vertex.

We have chosen the term ``cusp with coparticles'' in analogy with ``cusp with particles'' from our paper \cite{FI09}. Both are derived from ``manifold with particles'' which is the term in the physics literature for 3-dimensional anti-de Sitter manifolds with conical singularities along time-like lines. In the recent years, manifolds with particles have found interesting applications in the geometry of low-dimensional manifolds, \cite{BS06, KS07,MS09,BBS09}.

\subsection{Related results}
\label{subsec:RelRes}
Let us return to the Rivin-Hodgson Theorem. The Gauss image of the boundary of a polyhedron $P$ has a more straightforward definition as the boundary of the polar dual $P^*$. The polar dual to a convex hyperbolic polyhedron is a convex polyhedron in the \emph{de Sitter space}. The de Sitter space is the one-sheeted hyperboloid in the Minkowski space, and the polarity between convex polyhedra in $\H^3$ and $\dS^3$ is established via the Minkowski scalar product, see Section \ref{subsec:deSitter}.

Thus Theorem \ref{thm:RivHod} can be reformulated as follows.
\begin{thmprime}\label{thmrivinprime}
Let $g$ be a spherical cone metric on the sphere $\Sph^2$ such that the following two conditions hold:
\begin{enumerate}
\item All cone angles of $(\Sph^2, g)$ are greater than $2\pi$.
\item All closed geodesics on $(\Sph^2, g)$ have lengths greater than $2\pi$.
\end{enumerate}
Then $(\Sph^2, g)$ can be isometrically embedded in the de Sitter space as a convex polyhedral surface. This embedding is unique up to isometry, in the class of surfaces that don't bound a ball in $\dS^3$.
\end{thmprime}

In this form, the Rivin-Hodgson theorem resembles the following theorem of Alexandrov, \cite{alex42, Ale05}.
\begin{thm}[Alexandrov]
\label{thm:Alex}
Let $g$ be a Euclidean cone metric on the sphere $\Sph^2$ such that all of its cone angles are less than $2\pi$. Then $(\Sph^2, g)$ can be isometrically embedded in $\R^3$ as a convex polyhedral surface. The embedding is unique up to isometry.

Equivalently, there exists a unique compact convex Euclidean polyhedron $P$ whose boundary is isometric to $(\Sph^2, g)$.
\end{thm}


Theorems \ref{thm:RivHod}, \ref{thm:Alex}, and \ref{thm:Main} are special cases of a general statement which says roughly the following.

Let $S$ be a closed surface of an arbitrary genus, and let $g$ be a cone metric on $S$ of constant curvature with either all cone angles greater or all less than $2\pi$. Then the metric $g$ can be extended to a (Riemannian or Lorentzian) metric of constant curvature inside a manifold of a certain topological type (ball or cusp in the above examples) with convex boundary $S$. This result can be viewed as a geometrization with boundary conditions. Another equivalent formulation is that the universal cover of $(S, g)$ can be embedded as a convex polyhedral surface invariant under an appropriate action of the group $\pi_1(S)$.

In a precise form, the general statement is formulated as Problem~1 in \cite{FI09}. Other special cases are the main theorems of \cite{Fillastre1,Fillastre2,FI09}. Theorem \ref{thm:Main} deals with the last remaining case.

Alexandrov proved Theorem \ref{thm:Alex} by what we call Alexandrov's method, see the second paragraph after Theorem \ref{thm:Mink}. In his book \cite{Ale05}, Alexandrov asked whether Theorem \ref{thm:Alex} can be proved through a variational approach, similarly to Minkowski's theorem. This was finally done in \cite{BI08}. The functional used in \cite{BI08} is the \emph{discrete Hilbert-Einstein functional} on the space of certain singular polyhedra with a fixed metric on the boundary, and is in a sense dual to the volume functional used in the proof of Minkowski's theorem.

In \cite{FI09}, the discrete Hilbert-Einstein functional is used to characterize convex polyhedral cusps in terms of their boundary metrics, that is to resolve the case of genus one, hyperbolic metric, and cone angles less than $2\pi$ of the general statement mentioned above. The function $V$ used in the present paper can also be interpreted as the discrete Hilbert-Einstein functional of the polar dual ``cusp with particles'' in the de Sitter space, see Section~\ref{subsec:deSitter}.

Andreev's Theorem \ref{thm:And} can be reformulated in terms of circle patterns on the sphere with non-obtuse intersection angles. In this form, it was extended by Thurston to circle patterns on the torus and on the higher genus surfaces, \cite[Theorem 13.7.1]{Thurcour1}. Thus, our Theorem \ref{thm:Main} generalizes the torus case of Thurston's theorem in the same direction as Rivin-Hodgson's Theorem \ref{thm:RivHod} generalizes Andreev's theorem. For more details, see Section \ref{subsec:Thurston}.

Theorems on convex embeddings of convex polyhedral metrics have smooth analogs. The most renowned is the Weyl problem that asks whether $\Sph^2$ with a Riemannian metric of positive Gauss curvature can be embedded in $\R^3$ as the boundary of a convex body. The Weyl problem was solved through PDE in \cite{Lewy38, Nir53}, and through polyhedral approximation in \cite{alex42, Pog73}.

The metric on the Gauss image of a smooth immersed surface is the third fundamental form of the surface. A smooth analog of the Rivin-Hodgson theorem was proved in \cite{JMS1,JMS2} and later generalized in \cite{Schconvex}. The theorem in \cite{Schconvex} proves, for a compact hyperbolizable $3$--manifold with boundary, the existence of a hyperbolic metric with a prescribed third fundamental form of the boundary. For a related result see \cite{SchLab}. The work \cite{Schconvex} deals only with compact manifolds. Its extension to cusped manifolds would imply a smooth analog of Theorem~\ref{thm:Main}.

A polyhedral analog of the theorem from \cite{Schconvex} is an open problem.

\subsection{Plan of the paper}
Convex polyhedral cusps and their Gauss images are defined and discussed in detail in Section \ref{sec:DefPre}. In Section \ref{sec:CuspDef}, we define convex polyhedral cusps with coparticles.

The space $\M^*$ of all convex polyhedral cusps with coparticles with a given Gauss image is studied in Section \ref{sec:SpaceOfCusps}. The main results here are Proposition \ref{prp:HeightsDefCusp} stating that a cusp is uniquely determined by its heights $\{h_i\}$, Proposition \ref{prp:DescrM} describing $\M^*$ through a system of inequalities on $\{h_i\}$, and Proposition \ref{prp:StructM} showing that $\M^*$ is a contractible manifold with corners.

In Section \ref{sec:Functional}, we define the function $V$ on the space $\M^*$. We apply the Schl\"afli formula to compute its partial derivatives. This reduces Theorem \ref{thm:Main} to the statement that $V$ has a unique critical point in $\M^*$. Also we prove the concavity of $V$ that, with some additional effort, implies local rigidity of convex polyhedral cusps with coparticles, Theorem \ref{thm:LocRig}.

Section \ref{sec:Proof} contains the proof of Theorem \ref{thm:Main}. The principal task here is to analyze the behavior of $V$ at the boundary of $\M^*$ and at the infinity. That done, Morse theory on manifolds with corners can be applied. The existence and uniqueness of a critical point for $V$ follows by counting the indices of all critical points.

In Subsection \ref{subsec:AndThm}, an analog of Andreev's Theorem \ref{thm:And} for convex polyhedral cusps is stated, Theorem \ref{thm:AndreevCusps}. It follows from Theorem \ref{thm:Main} as a special case, but under assumptions of Theorem \ref{thm:AndreevCusps} some parts of the proof are simplified. We also discuss connections with circle patterns on the torus. Subsection~\ref{subsec:deSitter} puts Theorem~\ref{thm:Main} and its proof in the context of the de Sitter geometry.

Finally, Appendix~\ref{sec:FormCorn} contains some formulas used in the main text.

\subsection{Acknowledgements}
We are grateful to Jean-Marc Schlenker for bringing this problem to our attention and for the constant interest to our work. We thank Boris Springborn who derived the explicit formula for the function~$V$.

Parts of the work were done during the first author's visits to the TU Berlin and the second author's visit to the University of Cergy--Pontoise. Final touches to the text were given during the second author's stay at the Kyushu University in Fukuoka. We thank all three institutions for their hospitality.

\section{Definitions and preliminaries}
\label{sec:DefPre}
\subsection{Convex polyhedral cusps and their Gauss images}
\label{subsec:GaussImage}
Let $\Tor^2$ denote the 2--dimensional torus.
\begin{dfn} \label{def:Cusp}
A \emph{hyperbolic cusp with boundary} is a complete hyperbolic manifold of finite volume homeomorphic to $\Tor^2 \times [0, +\infty)$. We say that a cusp has \emph{convex polyhedral boundary} if at every boundary point it is locally isometric to a convex polyhedral cone in $\H^3$.
\end{dfn}
For brevity, we refer to hyperbolic cusps with convex polyhedral boundary as \emph{convex polyhedral cusps}.

\begin{figure}[ht]
\begin{center}
\includegraphics[width = 0.4\textwidth]{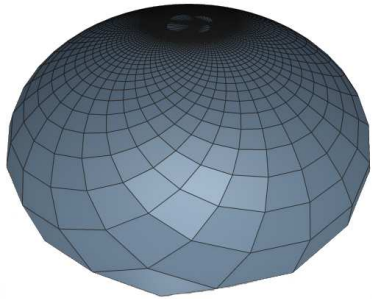}
\end{center}
\caption{A convex parabolic polyhedron in the Klein model.}
\label{fig:ParPol}
\end{figure}

In an obvious manner one can define vertices, edges, and faces of a convex polyhedral cusp $M$. A priori, an edge can be a loop, and a face can be non-simply connected. The universal cover $\widetilde{M}$ embeds in $\H^3$ as a \emph{convex parabolic polyhedron}: the convex hull of finitely many orbits of a $\Z^2$-action on $\H^3$ by parabolic isometries, see Figure \ref{fig:ParPol}. A face of $M$ lifts to a face of $\widetilde{M}$. Since every face of $\widetilde{M}$ is a convex hyperbolic polygon, every face of $M$ is also a convex hyperbolic polygon. Thus the interior of a face of a convex polyhedral cusp is simply connected, although there can be identifications on the boundary. For more details see \cite[Section 2.1]{FI09}.

By taking the $\Z^2$--quotient of $\H^3$ we obtain a complete hyperbolic manifold $N$ in which the convex polyhedral cusp $M$ is embedded isometrically. Let $v$ be a vertex of $M$. Those vectors in $T_v N$ that are tangent to curves $\gamma \co [0,\epsilon) \to M$ form the \emph{tangent cone} $C_v M$. Clearly, $C_v M$ is a convex polyhedral cone in $T_v N$. The polar dual $(C_v M)^\circ$ of $C_v M$ is defined as
$$
(C_v M)^\circ := \{x \in T_v N \,|\, \langle x, y \rangle \le 0 \mbox{ for all } y \in C_v M\}.
$$
Geometrically, $(C_v M)^\circ$ is a cone positively spanned by the outward normals at $v$ to the faces of $C_v M$.

\begin{dfn}
The convex spherical polygon
$$
\Pi_v := \{x \in (C_v M)^\circ \,|\, \|x\|=1\}
$$
is called the \emph{Gauss image} of the vertex $v$.
\end{dfn}

If $v$ and $w$ are two vertices of $M$ joined by an edge, then the parallel transport along this edge maps a side of $\Pi_v$ isometrically to a side of $\Pi_w$. Since faces of $M$ are simply connected, after performing all the gluings we obtain a surface homeomorphic to the torus $\Tor^2$. Away from the vertices of the polygons, this surface is locally isometric to the unit sphere; the vertices become cone points.

\begin{dfn}
The spherical cone-surface glued from the Gauss images of the vertices of $M$ is called the \emph{Gauss image} of $\partial M$. 
\end{dfn}
For brevity, we sometimes say `Gauss image of $M$'.

The cell structure $\bigcup_v \Pi_v$ of the Gauss image is dual to the face structure of $\partial M$: a vertex $v$ of $M$ gives rise to a face $\Pi_v$, an edge $vw$ of $M$ gives rise to a dual edge shared by the faces $\Pi_v$ and $\Pi_w$, and a face $F_i$ of $M$ gives rise to a vertex $i$ where the outward unit normals to the face $F_i$ at all of its vertices are glued together.

\begin{dfn}
The cell decomposition $\bigcup_v \Pi_v$ of the Gauss image is called the \emph{dual tesselation} associated with $M$.
\end{dfn}

For brevity, we will often say `the dual tesselation of $M$'.

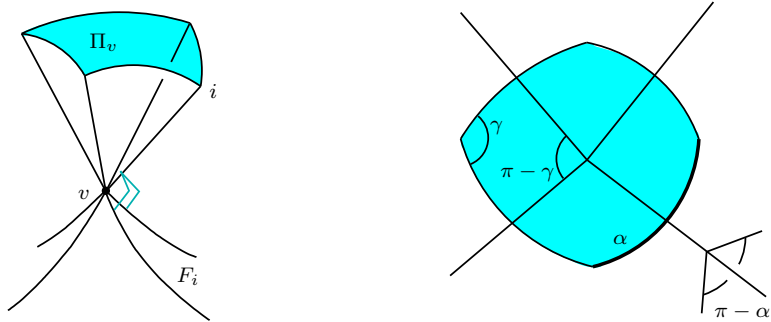
\begin{figure}[ht]
\begin{center}
\input{Fig/ComplAngles.tex}
\end{center}
\caption{Complementary angles in a convex polyhedral cone and its Gauss image.}
\label{fig:DualAngles}
\end{figure}

The length of an edge in the dual tesselation equals $\pi$ minus the dihedral angle at the corresponding edge of the cusp. An angle of a polygon $\Pi_v$ equals $\pi$ minus the corresponding plane angle at the vertex $v$. See Figure \ref{fig:DualAngles}.

In Theorem \ref{thm:Main}, we need to recover a convex polyhedral cusp $M$ from its Gauss image. Note that only the intrinsic metric of the Gauss image is given but not the dual tesselation. Thus we know the number of faces of $M$ but don't know which pairs of faces are adjacent and what are the dihedral angles between them.

However, we can compute the areas of faces. Denote by $\theta_i$ the cone angle at a vertex $i$ of the Gauss image. Since the angle of $\Pi_v$ at $i$ is equal to the exterior angle of $F_i$ at $v$, we have
\begin{equation}
\label{eqn:Omega}
\theta_i = 2\pi + \area(F_i).
\end{equation}

\begin{exl}[One-vertex cusps]
\label{exl:OneVert}
Consider an action of $\Z^2$ on $\H^3$ by parabolic isometries fixing the point $o \in \partial \overline{\H^3}$. Choose an arbitrary point in $\H^3$ and denote by $P$ the convex hull of its orbit. The orbit is a lattice in a horosphere $H$ centered at $o$. The boundary of $P$ projected from $o$ down to $H$ gives a $\Z^2$-invariant tesselation of $H$ with the lattice as the vertex set. The convexity of $P$ implies that $H$ is tessellated either by rectangles or by copies of an acute-angled triangle. Let $abcd$ be a rectangle or the union of two triangles $abc$ and $acd$ of the tesselation. See Figure \ref{fig:Cusp1Vert} left, where the tesselation is projected to $\partial \overline{\H^3}$. The convex parabolic polyhedron depicted on Figure \ref{fig:ParPol} is also the convex hull of a single orbit, and has quadrilateral faces.

\begin{figure}[ht]
\begin{center}
\input{Fig/OneVertCusp.tex}
\end{center}
\caption{A one-vertex convex polyhedral cusp in the Poincar\'e model and its Gauss image.}
\label{fig:Cusp1Vert}
\end{figure}
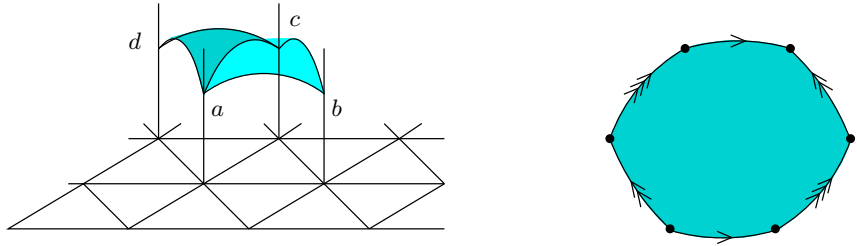

The convex hull of the points $a$, $b$, $c$, $d$, and $o$, with $o$ removed, is the fundamental domain for the $\Z^2$--action on $P$. Thus the convex hyperbolic cusp $P / \Z^2$ is isometric to the semi-ideal convex hyperbolic polyhedron $oabcd$ with $oab$ glued isometrically to $odc$ and $obc$ glued isometrically to $oad$. The points $a$, $b$, $c$, $d$ are identified to a single vertex $v$ of the cusp $M$.

The Gauss image $\Pi_v$ of the vertex $v$ is either a quadrilateral, if $abcd$ is a rectangle, or a hexagon, if $abcd$ is a parallelogram subdivided by the diagonal $ac$. In both cases, the pairs of opposite edges of $\Pi_v$ have equal length and are glued together to yield the Gauss image of $\partial M$. The dual tesselation has either one hexagonal face and two trivalent vertices or one quadrangular face and one four-valent vertex.
\end{exl}

\subsection{Closed geodesics}
\begin{dfn}
\label{dfn:Geod}
A geodesic on a spherical cone-surface is a locally length minimizing curve.
\end{dfn}
At a regular point a geodesic is locally a great circle; if a geodesic passes through a cone point, it spans angles of at least $\pi$ on both sides.

\begin{lem}
\label{lem:cont geod}
If $(\Tor^2, g)$ is the Gauss image of a convex polyhedral cusp $M$, then every contractible closed geodesic on $(\Tor^2, g)$ has length greater than $2\pi$.
\end{lem}
\begin{proof}
Consider the convex parabolic polyhedron $P = \widetilde{M}$. Clearly, the Gauss image of $\partial P$ is the universal cover of $(\Tor^2, g)$. Thus, for every contractible closed geodesic $\gamma$ on $(\Tor^2, g)$ there is a closed geodesic $\widetilde{\gamma}$ of the same length on the Gauss image of $\partial P$. By \cite[Proposition 3.6]{RivHod93}, every closed geodesic on the Gauss image of a compact convex polyhedron in $\H^3$ has length greater than $2\pi$. Although $P$ is not compact, the curve $\widetilde{\gamma}$ runs through the Gauss images of a finite number of vertices of $P$:
$$
\widetilde{\gamma} \subset \bigcup_{i=1}^n \Pi_{v_i}.
$$
Let $P'$ be the convex hull of the points $\{v_i\}_{i=1}^n$ and of all of their neighbors in $P$. Then the Gauss images of $v_i$ with respect to $P'$ are the same as with respect to $P$. Therefore, $\widetilde{\gamma}$ can be viewed as a geodesic on the Gauss image of $\partial P'$, and has length greater than~$2\pi$.
\end{proof}

Example \ref{exl:OneVert} shows that the Gauss image of a convex polyhedral cusp can contain \emph{non-contractible} closed geodesics of length less than~$2\pi$.

Lemma~\ref{lem:cont geod} and equation (\ref{eqn:Omega}) show that assumptions of Theorem~\ref{thm:Main} are necessary for a metric $g$ to be the metric of the Gauss image of a convex polyhedral cusp. In the sequel we always assume that $(\Tor^2, g)$ satisfies those assumptions. Condition on the lengths of geodesics is used in Section \ref{subsec:TessTriang} and plays a crucial role in Section~\ref{sec:VInf}.

\subsection{Geodesic tesselations and triangulations}
\label{subsec:TessTriang}
As noted in Section \ref{subsec:GaussImage}, the intrinsic metric $g$ of the Gauss image of a convex polyhedral cusp does not tell us much about the dual tesselation. Therefore, in the process of proving Theorem \ref{thm:Main}, we will consider all possible tesselations of $(\Tor^2, g)$.

Let us call a geodesic arc joining two cone points \emph{regular} if its interior does not pass through cone points.

\begin{dfn}
\label{dfn:TessTriang}
A \emph{tesselation} of $(\Tor^2, g)$ is a decomposition into vertices, edges, and faces, where
\begin{itemize}
\item vertices are exactly the cone points of $g$;
\item edges are simple regular geodesic arcs of lengths less than $\pi$ with mutually disjoint interiors;
\item faces are connected components of $\Tor^2$ minus the union of vertices and edges.
\end{itemize}

A tesselation is called a \emph{triangulation} if all faces are triangles.
\end{dfn}

A face of a tesselation can be non-simply connected. In particular, there is always a tesselation without edges, whose only face is the complement of $\Tor^2$ to the set of cone points. Edges are allowed to be loops, and there can be multiple edges. Thus a triangulation is not required to be a simplicial complex.

The dual tesselation of a convex polyhedral cusp can be refined to a triangulation by subdividing every polygon $\Pi_v$ into triangles by diagonals.


\begin{lem}
\label{lem:Troyanov}
For every compact spherical cone surface, the lengths of regular geodesics with endpoints at singularities form a discrete subset of $\R$.
\end{lem}
\begin{proof}
This lemma is an analog of \cite[Corollary 1]{ILTC} that deals with Euclidean cone surfaces. The argument in the proof of \cite[Proposition 1]{ILTC} shows that the number of regular geodesics whose length is less than $L$ and is not a multiple of $\pi$, and whose endpoints are given cone points $i$ and $j$, is finite. Since there are only finitely many cone points, it follows that the set of lengths is discrete.
\end{proof}

\begin{lem}
\label{lem:TriangFinite}
Every compact spherical cone surface has only a finite number of triangulations in the sense of Definition \ref{dfn:TessTriang}.
\end{lem}
\begin{proof}
Again, by the argument in the proof of \cite[Proposition 1]{ILTC}, the number of regular geodesic arcs of length less than $\pi$ between cone points is finite. It follows that the number of triangulations if finite.
\end{proof}

\begin{lem}
\label{lem:TriangExist}
A spherical cone-surface $(S, g)$ can be triangulated in the sense of Definition \ref{dfn:TessTriang} if and only if for every point $x \in S$ there exists a cone point $i$ at distance less than~$\pi$ from $x$.
\end{lem}
\begin{proof} (See also \cite[Proposition 3.1]{Thurart1} and \cite[Section 3]{Riv05}.)

If $(S,g)$ admits a triangulation, then every point $x \in S$ lies at distance less than $\pi$ from every vertex of a triangle containing $x$.

For the opposite direction, let us describe the \emph{Delaunay tesselation} of $(S,g)$. Define the \emph{Voronoi cell} of a cone point $i$ as
$$
V_i = \{x \in S\ |\ \dist(x,i) \le \dist(x,j)\mbox{ for all cone points }j\}.
$$
It is not hard to show that the interior of a Voronoi cell is an open disk, and the boundary is a polygonal curve. The geodesic segments of the boundary are called Voronoi edges. Let $x$ be an interior point of a Voronoi edge $e$. Then $x$ has exactly two shortest arcs to cone points (these might be arcs joining $x$ with two different cone points, or two different arcs to one point). Thus, if $\ell_x < \pi$ is the length of the two shortest arcs, then $(S,g)$ contains an immersed open disk $D_x$ of radius $\ell_x$ centered at $x$. The disk can overlap if, for example, the cone angle at one of the cone points closest to $x$ is less than $\pi$. Develop $D_x$ on the sphere and consider the geodesic arc  that joins the images of the cone points on its boundary. The image of this arc on $(S,g)$ is the \emph{Delaunay edge} dual to the Voronoi edge $e$.

Endpoints of Voronoi edges are called Voronoi vertices. Let $x$ be a Voronoi vertice. Let $xi_1, \ldots, xi_n$, in this cyclic order, be the shortest arcs from $x$ to cone points of $(S,g)$. Then $n \ge 3$. Note that some of the vertices $i_1,\ldots,i_n$ may coincide, but the arcs are different. Again, denote by $\ell_x < \pi$ the common length of these arcs and consider the immersed open disk $D_x$ of radius $\ell_x$ centered at $x$. The cone points $i_1,\ldots, i_n$ lie on the boundary of $D_x$. The Voronoi edges starting from $x$ go along the bisectors of the angles $i_sxi_{s+1}$ (note that one of these angles can be bigger than $\pi$). The dual Delaunay edges bound a spherical polygon inscribed in $D_x$. This polygon is called a \emph{Delaunay cell}.

Thus, the Delaunay edges cut the surface $(S, g)$ into Delaunay cells which are inscribed in disks and thus are convex polygons. After subdividing each Delaunay cell into triangles by diagonals, we obtain a triangulation of $(S,g)$.
\end{proof}

\begin{cor}
The spherical cone surface $(\Tor^2, g)$ can be triangulated, and this in only finitely many ways.
\end{cor}
\begin{proof}
Let us show that for every point $x \in (\Tor^2, g)$ there exists a cone point $i$ at distance less than~$\pi$ from $x$. Assume the converse, and let $x$ be a point that lies at a distance at least $\pi$ from all cone points. Let $\gamma$ be a closed curve that bounds the disk of radius $\pi$ centered at $x$. Then $\gamma$ is a contractible curve of length $2\pi$. If it does not go through cone points, then at every point it is locally a great circle. If $i$ is a cone point on $\gamma$, then one of the angles spanned by $\gamma$ at $i$ equals $\pi$, and the other is larger than $\pi$ because the cone angle at $i$ is larger than $2\pi$. Thus $\gamma$ is in any case a closed contractible geodesic of length $2\pi$. This contradicts the second assumption of Theorem \ref{thm:Main}.

Thus by Lemma \ref{lem:TriangExist} the cone surface $(\Tor^2, g)$ can be triangulated. The finiteness of the number of the triangulations follows from Lemma \ref{lem:TriangFinite}.
\end{proof}

The next example shows that the condition on the lengths of closed geodesics is not necessary for the existence of a triangulation.

\begin{exl}[A triangulable cone-surface with a short contractible  closed geodesic]
Consider a combinatorial triangulation of the torus such that every vertex has degree at least four. Remove one of the triangles from the triangulation and glue on its place an octahedron with one face removed. This yields a combinatorial triangulation of the torus with all vertices of degree at least four and with a closed contractible path $\gamma$ consisting of three edges. Realize every triangle as an equilateral spherical one with angles greater than $\frac{\pi}2$. Then every vertex becomes a cone point with angle greater than $2\pi$. Since the side lengths of the triangles are less than $\frac{2\pi}3$, the path $\gamma$ has length less than $2\pi$. Besides, $\gamma$ is a geodesic because it spans angles greater than $\frac{3\pi}2$ on both sides.
\end{exl}

Sometimes the Gauss image of a convex polyhedral cusp has only one triangulation, which implies that the dual tesselation is determined by the intrinsic metric of the Gauss image.

\begin{exl}[Cone-surface with a unique triangulation]
\label{exl:UniqueTriang}
Consider an arbitrary combinatorial triangulation $T$ of the torus. Let $\ell$ be a map from the edge set of $T$ to the interval $[\frac{\pi}2, \pi)$ such that for every triangle of $T$ there exists a spherical triangle with side lenghts given by $\ell$. Let $(\Tor^2, g)$ be the corresponding spherical cone-surface. As shown in \cite[Proposition 2.4]{Hod92}, the only geodesic arcs of length less than $\pi$ between cone points on $(\Tor^2,g)$ are the edges of $T$. It follows that $T$ is the unique triangulation of $(\Tor^2,g)$.

Conditions that the triangulation $T$ needs to fulfill in order that $(\Tor^2, g)$ satisfies the assumptions of Theorem \ref{thm:Main} are discussed in Section \ref{subsec:AndThm}. As in this special case the metric $g$ determines the tesselation uniquely, Theorem \ref{thm:Main} can be reformulated in the spirit of Andreev's theorem or in terms of circle patterns.

%
%
\end{exl}

\section{Cusps with coparticles}
\label{sec:CuspDef}
We now want to extend the notion of a convex polyhedral cusp. A cusp will be allowed to have singular lines (so called coparticles) in the interior, in such a way that the Gauss image of its boundary still can be defined. A cusp with coparticles is composed from building blocks. Let us first provide a decomposition of a usual convex polyhedral cusp into such blocks.

Let $\widetilde{M} \subset \H^3$ be the universal cover of a convex polyhedral cusp $M$, and let $o \in \partial\overline{\H^3}$ be the lift of the apex of the cusp. Assume that the polyhedron $\widetilde{M}$ is simple, that is, each vertex of $\widetilde{M}$ is shared by exactly three edges and three faces. Then the hyperbolic planes through $o$ orthogonal to the edges cut $\widetilde{M}$ into polyhedra combinatorially isomorphic to the cube with one ideal vertex. We call these pieces \emph{bricks}. The decomposition of $\widetilde{M}$ into bricks descends to $M$.

There are two points that make the above picture more complicated, and because of which we need to proceed in this section rather formally.

First, if the plane through $o$ orthogonal to an edge of $\widetilde{M}$ misses the edge and intersects its extension, then there appears a brick whose boundary has self-intersections. The decomposition of $\widetilde{M}$ into bricks is meant then in an algebraic sense, some parts of some bricks counted with the minus sign.

Second, $\widetilde{M}$ can be non-simple. In order to obtain a decomposition into bricks we proceed as follows. We resolve each vertex $v$ of degree higher than three into a trivalent tree. Every new edge has zero length and separates two faces incident to $v$. We view this edge as lying on the intersection line of the planes spanned by these faces. The decomposition of $\widetilde{M}$ must be performed with respect to all edges, including zero ones.

A cusp with coparticles is composed from a set of bricks in a similar way, with the difference that the dihedral angles of bricks around a semi-ideal edge can sum up to an angle other than $2\pi$. We require that the links of bricks at their ideal vertices form a torus and induce on it a Euclidean cone metric modulo scaling (this is the analog of the completeness of a convex polyhedral cusp).

One is tempted to define a convex polyhedral cusp with coparticles as a hyperbolic cone-manifold with properties similar to those in Definition \ref{def:Cusp} and with singular lines orthogonal to the faces. This would be too restrictive because the foot of the perpendicular to a face can lie outside the face; in this case the complex of bricks is an abstract object that does not give rise to a cone-manifold. However, the polar dual of a convex polyhedral cusp with coparticles is always a de Sitter cone-manifold, see Section \ref{subsec:deSitter}.

\subsection{Corners}
\label{subsec:Corners}
Corners are simpler than bricks and allow to avoid difficulties that arise when the boundary of a brick intersects itself.
\begin{dfn}
A \emph{co-oriented plane} is a plane $L \subset \H^3$ together with a choice of a positive half-space $L^+$ from the two half-spaces bounded by $L$. The other half-space is called negative and denoted by $L^-$.
\end{dfn}

\begin{dfn}
\label{dfn:Corner}
A \emph{corner} $(L_1, L_2, L_3; o)$ consists of three co-oriented planes $L_1$, $L_2$, $L_3$ in $\H^3$ that have exactly one common point in $\H^3$, and of a point $o \in \partial\overline{\H^3}$ that lies on the positive side of each of the planes $L_1$, $L_2$, $L_3$. 
A corner is considered up to an isometry of $\H^3$.
\end{dfn}

Basically, a corner consists of the ideal vertex of a brick and of the planes of its three compact faces.

Let $v$ be the intersection point of $L_1$, $L_2$, and $L_3$. The Gauss image of $v$ with respect to the cone $L_1^+ \cap L_2^+ \cap L_3^+$ is called the \emph{Gauss image of the corner}.

\begin{dfn}
\label{dfn:TrCorn}
A \emph{truncated corner} is a corner $(L_1, L_2, L_3; o)$ together with a horosphere $H$ centered at $o$. The numbers
$$
h_i := \dist(H, L_i),
$$
are called the \emph{heights} of the truncated corner.

Here $\dist(H, L_i)$ is the signed length of the common perpendicular between $H$ and $L_i$; the length is positive if $H \cap L_i = \emptyset$.
\end{dfn}

\begin{figure}[ht]
\begin{center}
\input{Fig/2Corners.tex}
\end{center}
\caption{Two-dimensional truncated corners. On the right hand side, the height $h_2$ is negative.}
\label{fig:TwoCorn}
\end{figure}
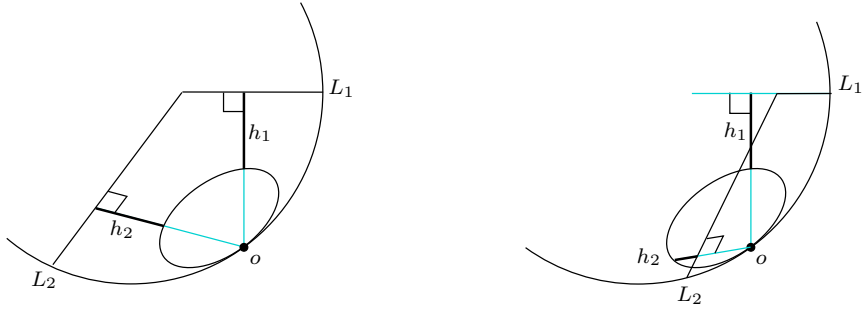

Figure \ref{fig:TwoCorn} shows examples of two-dimensional truncated corners in the Klein projective model of $\H^2$.

\begin{lem}
\label{lem:TrCorn}
Let $(123)$ be a spherical triangle with vertices $1$, $2$, $3$, and let $(h_1, h_2, h_3)$ be a triple of real numbers. Then there exists a unique truncated corner with Gauss image $(123)$ and respective heights $h_1, h_2, h_3$.
\end{lem}
\begin{proof}
The triangle $(123)$ determines the triple $(L_1, L_2, L_3)$ up to an isometry of $\H^3$. We need to show that there is a unique horosphere at the given distances from these three planes. We will prove this using the hyperbolic-de Sitter duality briefly described in \ref{subsec:deSitter}.

Hyperbolic half-spaces are identified with points of the de Sitter space by duality with respect to the Minkowski scalar product. Let $w_i$, $i=1,2,3$, be the point dual to the half-space $L_i^-$:
$$
L_i^- = \{x \in \H^3 \,|\, \langle w_i, x \rangle \le 0\}.
$$
Similarly, horospheres are identified with points in the upper half of the light cone:
\begin{equation}
\label{eqn:HoroVector}
H_u = \{p \in \H^3 \,|\, \langle p, u \rangle = -1\}
\end{equation}
is the horosphere associated with the vector $u$. It is not hard to show that
$$
\dist(H_u, L_i) = \log \langle u, w_i \rangle.
$$
(Note that $\langle u, w_i \rangle > 0$ iff the center of the horosphere lies on the positive side of $L_i$.) Thus the conditions $\dist(H_u, L_i) = h_i$ are equivalent to a non-degenerate system of three linear equations for four coordinates of the point $u$. The solution set is a time-like line not passing through the origin. Such a line has a unique intersection with the upper half of the light cone. The lemma is proved.
\end{proof}

\begin{lem}
\label{lem:ChangeTr}
Changing the truncation of a corner is equivalent to adding a constant to all of its heights:
$$
h'_i = h_i + c \quad \mbox{ for all } i \in \{1,2,3\}.
$$
\end{lem}
\begin{proof}
This follows from the fact that the set of concentric horospheres forms an equidistant family.
\end{proof}

\subsection{Gluing cusps from corners}
\label{subsec:CuspFromCorn}
Consider two spherical triangles $(123)$ and $(124)$ glued along the edge $(12)$. Any quadruple of real numbers $(h_1, h_2, h_3, h_4)$ determines two truncated corners. The corners can be fitted together so that they have a common truncating horosphere and common planes $L_1$ and $L_2$. There are two ways to do this; we choose that one where the co-orientations of the planes $L_3$ and $L_4$ induce different orientations on the line $L_1 \cap L_2$.

\begin{figure}[ht]
\begin{center}
\input{Fig/l12.tex}
\end{center}
\caption{Fitting two corners together produces an edge.}
\label{fig:l12}
\end{figure}
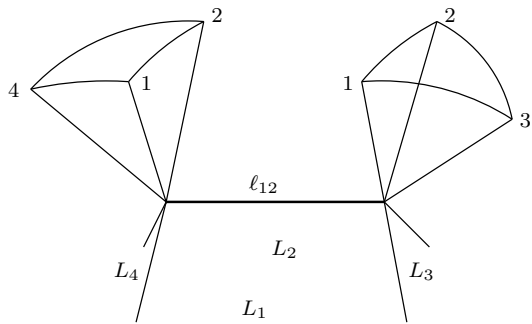

Denote by $\ell_{12}$ the signed length of the segment of the line $L_1 \cap L_2$ between the planes $L_3$ and $L_4$. The length is taken to be positive if the segment lies on the positive sides of $L_3$ and $L_4$, and negative if it lies on the negative side of both. See Figure \ref{fig:l12}, where we have $\ell_{12} > 0$.

\begin{dfn}
\label{dfn:TriangCuspPart}
Let $T$ be a geodesic triangulation of $(\Tor^2, g)$ and let
$$
h: i \mapsto h_i
$$
be an arbitrary map from the set of cone points of $g$ to $\R$.

For each triangle $ijk$ of $T$, construct the truncated corner with Gauss image $ijk$ and heights $h_i$, $h_j$, $h_k$. For each pair of adjacent triangles $ijk$ and $ijl$ of $T$, fit the corresponding corners and denote by $\ell_{ij}$ the length of the appearing edge.

If $\ell_{ij} \ge 0$ for all edges $ij$ of $T$, then the pair $(T,h)$ is called a \emph{triangulated truncated convex polyhedral cusp with coparticles} with Gauss image $(\Tor^2, g)$.
\end{dfn}

%

For brevity, we will say \emph{cusp} instead of \emph{convex polyhedral cusp with coparticles}. We are now coming to the main definitions of this section; we start by introducing an equivalence relation on the set of triangulated truncated cusps.

\begin{dfn}
A \emph{truncated cusp} is an equivalence class of triangulated truncated cusps, where two pairs $(T,h)$ and $(T',h')$ are equivalent if their heights are the same and their triangulations coincide up to edges with zero $\ell$-length:
\begin{equation}
\label{eqn:EqRel1}
(T, h) \sim (T', h') \Longleftrightarrow \left\{
\begin{array}{l}
h=h';\\
ij \in T, \ ij \notin T' \ \Rightarrow \ \ell_{ij} = 0;\\
ij \in T', \ ij \notin T \ \Rightarrow \ \ell'_{ij} = 0.
\end{array}
\right.
\end{equation}
\end{dfn}

\begin{dfn}
\label{dfn:DualTess}
Let $(T,h)$ be a triangulated truncated cusp. The tesselation of $(\Tor^2, g)$ obtained from $T$ by erasing all edges $ij$ such that $\ell_{ij} = 0$ is called the \emph{dual tesselation} of $(T,h)$.

The dual tesselation of a truncated cusp $M$ is defined as the dual tesselation of an arbitrary representative $(T,h)$ of $M$.
\end{dfn}

\begin{lem}
\label{lem:DTWD}
The dual tesselation of a truncated cusp is well-defined.
\end{lem}

Note that $(T, h) \sim (T', h')$ implies only that $T$ is a subdivision of the dual tesselation of $(T',h)$ and vice versa. The proof of Lemma \ref{lem:DTWD} will be given in Section \ref{subsec:SuppFunc}.

\begin{lem}
\label{lem:lijWD}
Let $M$ be a truncated cusp. Then, for every edge $ij$ of the dual tesselation of $M$, the length $\ell_{ij}$ does not depend on the choice of a representative $(T,h)$ of $M$.
\end{lem}
Rather discouraging, this lemma is not obvious. It is obvious, if the faces adjacent to the edge $ij$ are convex polygons: in this case the corners emerging from an arbitrary triangulation of a face fit all together to form a ``non-simple corner'' independent of the triangulation. The length $\ell_{ij}$ results from fitting two non-simple corners. But in general, a face of the dual tesselation can be non-convex and even non-simply connected. We postpone the proof of Lemma \ref{lem:lijWD} to Section \ref{subsec:SuppFunc}.

\begin{dfn}
A \emph{convex polyhedral cusp with coparticles} with Gauss image $(\Tor^2, g)$ is an equivalence class of triangulated truncated cusps with Gauss image $(\Tor^2, g)$ under an equivalence relation generated by \eqref{eqn:EqRel1} and
\begin{equation}
\label{eqn:EqRel2}
(T, h) \sim (T', h') \Longleftrightarrow \left\{
\begin{array}{l}
T = T';\\
h'_i = h_i + c  \mbox { for all } i.
\end{array}
\right.
\end{equation}
\end{dfn}

By Lemma \ref{lem:ChangeTr}, changing all heights by the same constant corresponds to a simultaneous change of all truncating horospheres. Therefore, by Lemmas \ref{lem:DTWD} and \ref{lem:lijWD}, a convex polyhedral cusp with coparticles has a well-defined dual tesselation, and well-defined edge lengths $\ell_{ij}$.

To prevent a confusion, let us stress that $\ell_{ij}$ is the length of an edge of the cusp and has nothing to do with the length of the edge $ij$ in the dual tesselation. In particular, $\ell_{ij}$ can be any positive real number, whereas the edges of the dual tesselation have lengths less than $\pi$.


\subsection{Curvatures of a cusp}
\label{subsec:CuspCurv}
Consider a truncated corner $(L_1, L_2, L_3; H)$. Let $o$ be the center of the horosphere $H$. For each $i \in \{1,2,3\}$, drop a perpendicular from the center of $H$ to $L_i$ and denote by $q_i$ its intersection point with $H$.
\begin{dfn}
The Euclidean triangle with vertices $q_1$, $q_2$, $q_3$ in $H$ is called the \emph{link} of the truncated corner $(L_1, L_2, L_3; H)$.
\end{dfn}

When two truncated corners are fitted, their links become glued along a side.
\begin{dfn}
Let $M$ be a truncated cusp represented by a pair $(T,h)$. By gluing the links of all truncated corners arising from $(T, h)$ we obtain a torus with a Euclidean cone metric that we call the \emph{link} of $M$.
\end{dfn}

\begin{lem}
\label{lem:LinkWD}
The link of a truncated cusp does not depend on the choice of a triangulation.
\end{lem}

This lemma will be proved in Section \ref{subsec:SuppFunc}.

The points $\{q_i\}$ of the link of $M$ are in a one-to-one correspondence with the cone points $\{i\}$ of $g$.

\begin{dfn}
\label{dfn:CuspCurv}
Let $M$ be a truncated cusp with Gauss image $(\Tor^2, g)$. For every cone point $i$ of $g$, denote by $\omega_i$ the total angle around the point $q_i$ in the link of $M$. The numbers
$$
\kappa_i = 2\pi - \omega_i
$$
are called the \emph{curvatures} of the truncated cusp $M$.
\end{dfn}

A change of the truncation of $M$ results in a scaling of the link of $M$. Therefore a convex polyhedral cusp with coparticles has also well-defined curvatures.

\begin{lem}
\label{lem:SumCurv}
For every convex polyhedral cusp with coparticles, the sum of the curvatures of coparticles vanishes:
\begin{equation}
\label{eqn:SumCurv}
\sum_i \kappa_i = 0.
\end{equation}
\end{lem}
\begin{proof}
Due to $\chi(\Tor^2)=0$, every triangulation of the torus has twice as many vertices as faces. In the triangulation of the link of $M$, the angle sum in every triangle is $\pi$, therefore the average angle around the vertex is $2\pi$. The lemma follows.
\end{proof}

As one can expect, making all curvatures zero is equivalent to finding a convex polyhedral cusp with a given Gauss image.
\begin{lem}
\label{lem:ZeroCurv}
Convex polyhedral cusps with Gauss image $(\Tor^2, g)$ are in a one-to-one correspondence with convex polyhedral cusps with coparticles with Gauss image $(\Tor^2, g)$ and all curvatures zero.
\end{lem}
\begin{proof}
Let $M$ be a convex polyhedral cusp. Truncate $M$ and represent it as a complex of corners. Clearly, all curvatures of thus obtained convex polyhedral cusp with coparticles are zero. Changing the truncation or decomposition produces a cusp with coparticles equivalent in the sense of \eqref{eqn:EqRel1} and \eqref{eqn:EqRel2}.

In the other direction, choose any representative $(T,h)$ of a convex polyhedral cusp with coparticles with zero curvatures. Due to $\omega_i = 2\pi$, the corners around the vertex $i$ can all be fitted together. We obtain a co-oriented plane $L_i$ and a set of co-oriented planes that cut in $L_i$ a convex polygon $F_i$. The apex $o$ lies on the positive side of all of the planes. This associates with each $i$ a semi-ideal pyramid with apex $o$ and base $F_i$. These pyramids are glued into a hyperbolic manifold which is complete because the gluing restricts to their truncations. The manifold has a convex polyhedral boundary, the topology of $\Tor^2 \times [0,+\infty)$, and finite volume. Thus it is a convex polyhedral cusp.
\end{proof}

\section{The space of cusps with coparticles}
\label{sec:SpaceOfCusps}
Everywhere in this section we mean by a cusp a convex polyhedral cusp with coparticles with Gauss image $(\Tor^2, g)$.

\subsection{Support function of a truncated cusp}
\label{subsec:SuppFunc}
Let $C = (L_1, L_2, L_3; H)$ be a truncated corner with vertex $v = L_1 \cap L_2 \cap L_3$ and Gauss image $\Delta$, see definitions in Section \ref{subsec:Corners}.

\begin{dfn}
The \emph{support function}
$$
\widetilde{h_\Delta} \co \Delta \to \R
$$
of a truncated corner $C$ is defined as follows. For a point $x \in \Delta \subset T_v\H^3$, denote by $L_x$ the plane through $v$ with normal $x$. Put
$$
\widetilde{h_\Delta}(x) = \dist(H, L_x),
$$
where $\dist(H, L_x)$ is the signed length of the common perpendicular between $H$ and $L_x$; the length is positive if $H \cap L_x = \emptyset$.
\end{dfn}

A triangulated truncated cusp $(T,h)$ is a complex of truncated corners, see Definition \ref{dfn:TriangCuspPart}.

\begin{dfn}
\label{dfn:SuppFunc}
The \emph{support function} of a triangulated truncated cusp $(T,h)$ is a function
$$
\widetilde{h_T} \co \Tor^2 \to \R
$$
that restricts to $\widetilde{h_\Delta}$ on each triangle $\Delta$ of $T$.
\end{dfn}

It is easy to see that the function $\widetilde{h_T}$ is well-defined on the edges of $T$. In particular,
$$
\widetilde{h_T}(i) = h_i
$$
for every cone point $i$ of $g$.

We are now going to characterize the support functions of triangulated truncated cusps and show that they are invariant under the equivalence relation \eqref{eqn:EqRel1}. The identification of truncated cusps with their support functions is a helpful tool that will be used several times in this section.

\begin{lem}
\label{lem:SuppFuncForm}
The support function of a truncated corner with Gauss image $\Delta \subset \Sph^2$ has the form
\begin{equation}
\label{eqn:SuppLike}
\widetilde{h_\Delta}(x) = \log \cos \dist(x, a) + b,
\end{equation}
where $a \in \Sph^2$, $b \in \R$.
\end{lem}
\begin{proof}
By the hyperbolic-de Sitter duality (see Section \ref{subsec:deSitter}), planes through the point $v \in \H^3$ correspond to the points on the dual de Sitter plane $v^*$. Namely, the point dual to the plane $L_x$ is its normal vector $x$, see the proof of Lemma \ref{lem:GaussImageDeSitter}. We have
\begin{equation}
\label{eqn:DistScal}
\widetilde{h_\Delta}(x) = \dist(H, L_x) = \log \langle u, x \rangle,
\end{equation}
where $u$ is a light-like vector in the Minkowski space $\R^{3,1}$ associated with the horosphere $H$ through \eqref{eqn:HoroVector}. The de Sitter plane $v^*$ is a 2--sphere, and $x \mapsto \langle u, x \rangle$ is the restriction of a linear function to a 2--sphere. Thus $\langle u, x \rangle$ is a multiple of $\cos \dist(x, a)$ for some point $a \in v^*$, where $\dist$ is the distance with respect to the intrinsic metric on $v^*$. Equation \eqref{eqn:SuppLike} follows.
\end{proof}


\begin{dfn}
A function on a subset of $\Sph^2$ is called \emph{support-like} if it has the form \eqref{eqn:SuppLike}. A function on $\Tor^2$ is called piecewise support-like or a \emph{PS function} if it is support-like on every triangle of some triangulation $T$ of~$(\Tor^2, g)$.
\end{dfn}


\begin{dfn}
Let $f \co \Tor^2 \to \R$ be a function which is smooth on every triangle of a triangulation $T$. Function $f$ is called \emph{Q-convex} if for every edge $e$ of $T$ and for every geodesic arc $\gamma$ that intersects $e$, the left derivative of $f|_\gamma$ at the intersection point with $e$ is less than or equal to its right derivative.
\end{dfn}

\begin{lem}
\label{lem:SuppFuncProp}
The support function of a triangulated truncated cusp $(T,h)$ is a Q-convex PS function.
\end{lem}
\begin{proof}
The support function is PS by definition. Let us show that it is Q-convex.

Let $ijk$ and $ijl$ be two adjacent triangles of $T$, and let $\gamma$ be a geodesic arc that intersects the edge $ij$. The function $\widetilde{h_{ijk}}$ can be extended in a support-like way to some neighborhood of the edge $ij$ in the triangle $ijl$. Let $x \in \gamma$ be a point in this neighborhood. It suffices to prove the inequality
\begin{equation}
\label{eqn:hijkhT}
\widetilde{h_{ijk}}(x) \le \widetilde{h_{ijl}}(x).
\end{equation}
Fit the corners corresponding to the triangles $ijk$ and $ijl$. Let $v = L_i \cap L_j \cap L_k$ and $w = L_i \cap L_j \cap L_l$ be their vertices. By definition,
$$
\widetilde{h_{ijl}}(x) = \dist(H, L_x),
$$
where $L_x$ is the plane through $w$ with normal $x$. It is easy to see that
$$
\widetilde{h_{ijk}}(x) = \dist(H, L'_x),
$$
where $L'_x$ is the parallel translate of $L_x$ along $vw$. Due to $\ell_{ij} \ge 0$, the plane $L'_x$ lies between the center of the horosphere $H$ and the plane $L_x$. This implies
$$
\dist(H, L'_x) \le \dist(H, L_x).
$$
Thus the inequality \eqref{eqn:hijkhT} holds, and the lemma is proved.
\end{proof}

We see that the sign of inequality between the left and the right derivatives of $f|_\gamma$ at the intersection point of $\gamma$ with $ij$ does not depend on the choice of geodesic $\gamma$. Therefore it makes sense to speak about PS functions \emph{convex across} $ij$ or \emph{smooth across} $ij$.

\begin{cor}
\label{cor:DualTess}
The faces of the dual tesselation of $(T,h)$ are the maximal subsets of $\Tor^2$ on which the function $\widetilde{h_T}$ is smooth.
\end{cor}
\begin{proof}
The argument in the proof of Lemma \ref{lem:SuppFuncProp} shows that $\ell_{ij} = 0$ if and only if $\widetilde{h_T}$ is smooth across $ij$.
\end{proof}

\begin{prp}
\label{prp:CuspFunc}
Truncated cusps with Gauss image $(\Tor^2, g)$ are in a one-to-one correspondence with Q-convex PS functions on $(\Tor^2, g)$. The bijection is established by associating to a truncated cusp the support function of any of its triangulations.
\end{prp}
\begin{proof}
First, let us show that every Q-convex PS function is the support function of a triangulated truncated cusp. Let $f$ be a Q-convex PS function. Choose a triangulation $T$ such that $f$ is support-like on the triangles of~$T$. Let $h$ be the restriction of $f$ to the set of cone points of the metric $g$. This yields a complex of corners $(T,h)$. The Q-convexity of $f$ implies $\ell_{ij} \ge 0$ for every edge $ij \in T$, by reversing the argument in the proof of Lemma \ref{lem:SuppFuncProp}. Thus the pair $(T,h)$ is a triangulated truncated cusp. Since a support-like function on a triangle is determined by its values at the vertices, $f$ is the support function of $(T,h)$.

Next, let us show that the support functions of triangulated truncated cusps $(T,h)$ and $(T',h')$ are equal if and only if $(T,h)$ and $(T',h')$ are equivalent in the sense of \eqref{eqn:EqRel1}.

Assume $\widetilde{h_T} = \widetilde{h'_{T'}} = f$. Then $h = h'$ since both are restrictions of $f$ to the set of cone points of $g$. If $ij \in T$ and $ij \notin T'$, then the latter implies that $f$ is smooth across the geodesic $ij$. Thus $\ell_{ij} = 0$. Similarly, $ij \in T'$ and $ij \notin T$ also implies $\ell_{ij}=0$. Therefore $(T,h) \sim (T',h')$.

In the inverse direction, let $(T,h) \sim (T',h)$. Let $R$ be the dual tesselation of $(T,h)$. Then, by the remark after Lemma \ref{lem:DTWD}, $T'$ is a subdivision of $R$. By Corollary \ref{cor:DualTess}, function $\widetilde{h_T}$ is support-like on every face of $R$, therefore it is support-like on every triangle of $T'$. Function $\widetilde{h_{T'}}$ is also support-like on every triangle of $T'$ and takes the same values at the vertices as $\widetilde{h_T}$. Thus $\widetilde{h_T} = \widetilde{h_{T'}}$.
\end{proof}

\begin{figure}[ht]
\begin{center}
\input{Fig/SingFace.tex}
\end{center}
\caption{Constructing a cusp with a non-simply connected face in the dual tesselation.}
\label{fig:Exotic}
\end{figure}
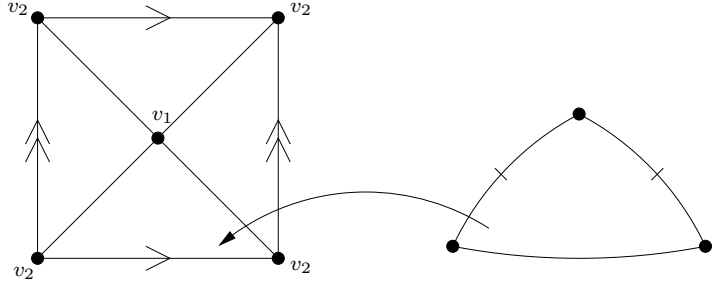

\begin{exl}[A cusp with coparticles with a non-simply connected face in the dual tesselation]
\label{exl:NSCFace}
Consider a triangulation of the torus shown on Figure \ref{fig:Exotic}, left. Equip each triangle with a spherical metric so that all of the four edges $v_1v_2$ have the same length and so that $v_1$ and $v_2$ become cone points with angles greater than $2\pi$. Define a support-like function on each triangle by taking each time the vertex $v_1$ as the point $a$ in the formula \eqref{eqn:SuppLike} and taking always the same number for $b$. It is easy to see that the PS function on the torus obtained in this way is smooth across the four $v_1v_2$ edges and convex across the two $v_1v_1$ edges. Thus the dual tesselation of the corresponding cusp has one face which is a punctured square.
\end{exl}

\begin{proof}[Proof of Lemma \ref{lem:DTWD}]
Follows from Corollary \ref{cor:DualTess} and Proposition~\ref{prp:CuspFunc}.
\end{proof}

\begin{proof}[Proofs of Lemmas \ref{lem:lijWD} and \ref{lem:LinkWD}]
Let $(T,h) \sim (T',h)$. Denote by $R$ their common dual tesselation. Let $S$ be a geodesic triangulation of $(\Tor^2, g)$ with vertices not necessarily at the cone points such that $S$ is a common subdivision of $T$ and $T'$ and every edge of $R$ is also an edge of $S$. (To obtain $S$, subdivide by diagonals the intersection $\Delta \cap \Delta'$ of every pair of triangles of $T$ and $T'$.) Let $\overline{h}$ be the restriction of the function $\widetilde{h_T} = \widetilde{h_{T'}}$ to the vertex set of $S$. The pair $(S,\overline{h})$ determines a complex of truncated corners.

For every edge $ij$ of $S$, denote by $\ell^S_{ij}$ the edge length obtained from fitting the corresponding corners of $(S,\overline{h})$. It is easy to see that
$$
\ell^S_{ij} = \ell_{ij}
$$
for all $ij \in R$. Since, in a similar way, $\ell^S_{ij} = \ell'_{ij}$, Lemma \ref{lem:lijWD} is proved.

To prove Lemma \ref{lem:LinkWD}, consider a triangle $\Delta$ of $T$ and its triangulation $S|_\Delta$ restricted from $S$. The corners that correspond to the triangles of $S|_\Delta$ can be fitted all together and all of their planes pass through one point. From this it is clear that the link of $(S|_\Delta, \overline{h})$ is isometric to the link of the corner $(\Delta, h)$. It follows that the link of $(S, \overline{h})$ is isometric to the link of $(T,h)$ and, in a similar way, to the link of $(T',h)$. Thus the links of $(T,h)$ and of $(T',h)$ coincide.
\end{proof}

%

\subsection{Heights define the cusp}
\label{subsec:HeightsToCusp}
The following construction allows to reduce the study of Q-convex PS functions to the study of convex piecewise linear functions.
\begin{dfn}
\label{dfn:Cone}
The \emph{cone} over $(\Tor^2, g)$ is a singular Riemannian manifold
$$
\cone(\Tor^2, g) = (\Tor^2 \times [0,+\infty)) / (\Tor^2 \times\{0\})
$$
with the metric $t^2 g + dt^2$ at $(x,t)$.

The \emph{cone of a function} $f \co (\Tor^2, g) \to \R$ is the function
\begin{eqnarray*}
\cone f \co  \cone(\Tor^2, g) & \to & \R,\\
(x,t) & \mapsto & tf(x).
\end{eqnarray*}
\end{dfn}

The map $x \mapsto (x, 1)$ is an isometric embedding of $(\Tor^2, g)$ in $\cone(\Tor^2, g)$ such that $(\cone f)|_{(\Tor^2, g)} = f$. With the apex $\Tor^2 \times \{0\}$ removed, the cone over $(\Tor^2, g)$ becomes an incomplete Euclidean cone-manifold. Every geodesic triangulation of $(\Tor^2, g)$ gives rise to a subdivision of $\cone(\Tor^2, g)$ into simplicial cones, which we also call a \emph{triangulation} of $\cone(\Tor^2, g)$.

\begin{lem}
\label{lem:htoe}
The map
$$
f \mapsto \cone \, \exp f
$$
establishes a bijection between Q-convex PS functions on $(\Tor^2, g)$ and positive convex PL functions on $\cone(\Tor^2, g)$.
\end{lem}
We say that a function on $\cone(\Tor^2, g)$ is positive, if it is positive everywhere except the apex. By a PL function we mean a function which is linear on every cone of some triangulation of $\cone(\Tor^2, g)$.

\begin{proof}
The exponent of a support-like function on a spherical triangle $\Delta$ has the form $r \cos \dist(x,a)$. The same form has the restriction to $\Delta$ of a linear function on $\cone\Delta$. The Q-convexity of $f$ is equivalent to the Q-convexity of $\exp f$ which is equivalent to the convexity of $\cone \exp f$.
\end{proof}

\begin{prp}
\label{prp:HeightsDefCusp}
A truncated cusp is uniquely determined by its heights.
\end{prp}
\begin{proof}
By Proposition \ref{prp:CuspFunc} and Lemma \ref{lem:htoe}, it suffices to show that every function $h$ on the set of cone points of $(\Tor^2, g)$ admits at most one extension to a convex PL function on $\cone(\Tor^2,g)$. This can be proved by a standard argument, see for example \cite[Lemma 3.8]{Izm08}.
\end{proof}

\begin{dfn}
Denote by $\M^*(g)$ the set of all convex polyhedral cusps with coparticles whose boundary has Gauss image $(\Tor^2, g)$.

Denote by $\Mt^*(g)$ the set of all truncated convex polyhedral cusps with coparticles and with Gauss image $(\Tor^2, g)$ of the boundary.
\end{dfn}
As the metric $g$ on $\Tor^2$ is fixed, we omit it from the notation and write simply $\M^*$ and $\Mt^*$.

Denote by $\Sigma \subset \Tor^2$ the set of cone points of the metric $g$. By Proposition \ref{prp:HeightsDefCusp}, we can identify $\Mt^*$ with a subset of $\R^\Sigma$ by putting
$$
\Mt^* = \{h\co \Sigma \to \R \,|\, h \mbox{ has a Q-convex PS extension to }(\Tor^2, g)\}.
$$
Since a change of truncation adds a constant to every height (Lemma \ref{lem:ChangeTr}), the space $\M^*$ can be identified either with a section of $\Mt^*$:
$$
\M^* = \Mt^* \cap \left\{\textstyle{\sum} h_i = 0\right\} \subset \R^\Sigma,
$$
or with a quotient of $\Mt^*$:
$$
\M^* = \Mt^* / \mathbb{L} \subset \R^\Sigma / \mathbb{L},
$$
where $\mathbb{L}$ is the 1-dimensional linear subspace of $\R^\Sigma$ spanned by the vector $(1,1,\ldots,1)$.

\subsection{Description of the space of cusps}
\label{subsec:DescrCusps}
Here the space $\Mt^* \subset \R^\Sigma$ is represented as the solution set of a system of inequalities. From this we extract some information on the topology and geometry of $\Mt^*$ and $\M^*$.

\begin{dfn}
\label{dfn:Quad}
A spherical quadrilateral is a (not necessarily convex) subset of $\Sph^2$ bounded by four arcs of big circles and contained in the interior of a hemisphere.

A \emph{quadrilateral} in $(\Tor^2,g)$ is a region bounded by four geodesic arcs with endpoints at cone points of $g$ and whose development on $\Sph^2$ is a spherical quadrilateral.
\end{dfn}

Some pairs of vertices or edges of a quadrilateral in $(\Tor^2,g)$ can coincide. For example, take a spherical quadrilateral with opposite pairs of sides equal in length and glue a torus $(\Tor^2, g)$ from it. The sides of the quadrilateral become two geodesic loops. These two loops, each ran twice, bound a quadrilateral in $(\Tor^2, g)$.

\begin{figure}[ht]
\begin{center}
\input{Fig/Quad.tex}
\end{center}
\caption{Every quadrilateral $ikjl$ in $(\Tor^2, g)$ with an angle at least $\pi$ at the vertex $i$ gives rise to an inequality \eqref{eqn:Quad}.}
\label{fig:ConcQuad}
\end{figure}
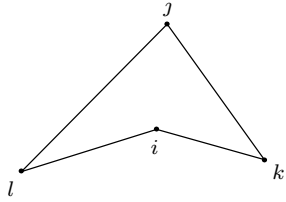

Let $ikjl$ be a quadrilateral in $(\Tor^2,g)$ with an angle at least $\pi$ at the vertex~$i$. Then it is easy to see that its angle at $j$ is less than $\pi$. Therefore, in the spherical development, the quadrilateral $ikjl$ is contained in the triangle $jkl$, see Figure \ref{fig:ConcQuad}. Given three real numbers $h_j$, $h_k$, $h_l$, let $\widetilde{h_{jkl}}$ be a support-like function on the quadrilateral $ikjl$ that takes values $h_j$, $h_k$, $h_l$ at the vertices $j$, $k$, $l$, respectively.

\begin{prp}
\label{prp:DescrM}
The space $\Mt^* \subset \R^\Sigma$ is the solution set of the system
\begin{equation}
\label{eqn:Quad}
h_i \le \widetilde{h_{jkl}}(i)
\end{equation}
that contains one inequality for every quadrilateral $ikjl$ in $(\Tor^2, g)$ with an angle at least $\pi$ at the vertex $i$.
\end{prp}


\begin{proof}
Although the proof can be carried out in terms of $h$, it is easier to work with the function
$$
e := \cone \exp h,
$$
see Definition \ref{dfn:Cone}. Thus, we have to show that a function
$$
\begin{array}{rrcl}
e\co & \Sigma & \to & (0, +\infty);\\
& i & \mapsto & e_i
\end{array}
$$
admits a convex PL extension
$$
\widetilde{e}\co \cone(\Tor^2, g) \to (0, +\infty)
$$
if and only if it satisfies the inequalities
\begin{equation}
\label{eqn:Quade}
e_i \le \widetilde{e_{jkl}}(i)
\end{equation}
for all quadrilaterals $ikjl$ with an angle at least $\pi$ at $i$. Here $\widetilde{e_{jkl}}$ is a linear function on $\cone(ikjl)$ that takes the values $e_j$, $e_k$, $e_l$ at $j$, $k$, $l$, respectively. As always, we develop $\cone(ikjl)$ in $\R^3$.

Let us prove the necessity of \eqref{eqn:Quade}. Assume that the function $\widetilde{e}$ exists. Then the functions $\widetilde{e}$ and $\widetilde{e_{jkl}}$ on $\cone(ikjl)$ agree at the points $j$, $k$, $l$. On the other hand, $\widetilde{e}$ is convex piecewise linear, and $\widetilde{e_{jkl}}$ is linear. This implies
$$
\widetilde{e}(i) \le \widetilde{e_{jkl}}(i).
$$
Since $\widetilde{e}(i) = e_i$, the inequality \eqref{eqn:Quade} holds.

Let us prove the sufficiency of \eqref{eqn:Quade}. For any triangulation $T$ of $(\Tor^2,g)$, denote by
$$
\widetilde{e_T}\co \cone(\Tor^2, g) \to (0, +\infty)
$$
the PL extension of $e$ with respect to $T$. Call an edge $ij$ of $T$ \emph{bad} if the function $\widetilde{e_T}$ is strictly concave across $\cone(ij)$. We have to show that there exists a triangulation such that all its edges are good.

A triangulation with good edges can be found using the \emph{flip algorithm}. To begin with, we take an arbitrary triangulation $T_0$ of $(\Tor^2,g)$. There is one, due to Lemma~\ref{lem:TriangExist}.

\begin{flp}
If the triangulation $T_0$ has bad edges, then let $ij$ be one of them, and let $ijk$ and $ijl$ be the triangles adjacent to $ij$. If the union of the triangles $ijk$ and $ijl$ contains a geodesic arc $kl$ in its interior, and this arc has length less than $\pi$, then replace the edge $ij$ by the edge $kl$. This operation is called a \emph{flip}.

Denote the new triangulation by $T_1$. If $T_1$ has bad edges, then choose one of them and flip it, if possible. Denote the new triangulation by $T_2$, and proceed further.
\end{flp}

In order to show that the flip algorithm outputs a triangulation with good edges, it suffices to prove the following two claims.

\begin{clm}
Any bad edge can be flipped.
\end{clm}

\begin{clm}
The flip algorithm terminates.
\end{clm}

Let us prove Claim 1. There are three reasons why an edge $ij$ can be impossible to flip; we eliminate each of them by contradiction.

First, assume that $ijl$ is the same triangle as $ijk$. Then in the triangle $ijk$ either the edge $ij$ is glued to $ik$ or $ji$ glued to $jk$. Thus either $i$ or $j$ has cone angle less than $\pi$, which contradicts the assumptions of Theorem \ref{thm:Main}.

Second, assume that the arc $kl$ inside the union of the triangles $ijk$ and $ijl$ has length at least $\pi$. Let $p$ be the point on this arc at distance $\pi$ from $k$. We can determine whether the edge $ij$ is good or bad by comparing the numbers $\widetilde{e_{ijk}}(p)$ and $\widetilde{e_{ijl}}(p)$. We have
$$
\widetilde{e_{ijk}}(p) = -e_k < 0.
$$
On the other hand,
$$
\widetilde{e_{ijl}}(p) > 0,
$$
because the point $p$ lies in the triangle $ijl$. It follows that the function $\widetilde{e_{T_0}}(p)$ is strictly convex across $\cone(ij)$. Thus the edge $ij$ is good.

Third, assume that the union of the triangles $ijk$ and $ijl$ does not contain a geodesic arc $kl$ in its interior. Then it is easy to show that $ikjl$ is a quadrilateral in $(\Tor^2,g)$ with one of the angles at $i$ and $j$ at least $\pi$. Let it be the angle at $i$. But then the inequality \eqref{eqn:Quade} implies that the function $\widetilde{e_{T_0}}$ is convex across $\cone(ij)$. Thus the edge $ij$ is good.

To prove Claim 2, note that flipping a bad edge increases the function $\widetilde{e_{T}}$ pointwise. Since, by Lemma \ref{lem:TriangFinite}, the number of triangulations of $(\Tor^2, g)$ is finite, the algorithm terminates.
\end{proof}

\begin{exl}
\label{exl:MtBig}
If $(\Tor^2, g)$ is as in Example \ref{exl:UniqueTriang}, then $(\Tor^2, g)$ contains no quadrilaterals. Thus the system \eqref{eqn:Quad} is empty, and we have $\Mt^* = \R^\Sigma$.
\end{exl}

\begin{prp}
\label{prp:StructM}
The spaces $\Mt^*$ and $\M^*$ are diffeomorphic to convex polyhedra of dimensions $|\Sigma|$ and $|\Sigma|-1$, respectively, with some faces removed.

The space $\M^* \subset \R^\Sigma / \mathbb{L}$ can be non-convex and unbounded.
\end{prp}
By definition, a convex polyhedron is an intersection of finitely many closed half-spaces. A \emph{convex polyhedron with some faces removed} is an intersection of closed and open half-spaces.

By a diffeomorphism between $\Mt^* \subset \R^\Sigma$ and ${\mathcal P} \subset \R^\Sigma$ we mean a diffeomorphism
$$
f \co \R^\Sigma \to {\mathcal U} \subset \R^\Sigma
$$
such that $f(\Mt^*) = {\mathcal P}$.

\begin{proof}
The diffeomorphism
\begin{eqnarray}
\exp \co \R^\Sigma & \to & (0,+\infty)^\Sigma,\label{eqn:Diffeo}\\
h & \mapsto & e = (e^{h_i})_{i \in \Sigma},\nonumber
\end{eqnarray}
maps $\Mt^*$ onto a subset of $\R^\Sigma$ that is the solution set of the system
\begin{eqnarray}
e_i & \le & \widetilde{e_{jkl}}(i), \label{eqn:Ineq1}\\
e_i & > & 0, \label{eqn:Ineq2}
\end{eqnarray}
where \eqref{eqn:Ineq1} is taken for every quadrilateral $ikjl$ (see Definition \ref{dfn:Quad}) with an angle at least $\pi$ at $i$. As $\widetilde{e_{jkl}}(i)$ is a linear function of $e_j$, $e_k$, and $e_l$, the system consists of linear inequalities. Thus $\exp\Mt^*$ is a convex polyhedron with some faces removed.

Let us show that the point $(1,1,\ldots,1) \in \R^\Sigma$ lies in the interior of $\exp\Mt^*$. For this we need to show that each of the inequalities \eqref{eqn:Ineq1} is strict when
$$
e_i = e_j = e_k = e_l = 1
$$
is substituted. The development of $\cone(ikjl)$ in $\R^3$ is a polyhedral cone spanned by four rays lying inside a half-space. The points $i$, $j$, $k$, $l$ lie on the spanning rays at the unit distance from the apex; the ray of $i$ is contained in the convex hull of the other three rays. The level set $\{\widetilde{e_{jkl}} = 1\}$ is a plane through the points $j$, $k$, and $l$. This plane intersects the $i$-ray at a point whose distance to apex is less than $1$. Therefore $\widetilde{e_{jkl}}(i) < 1$. Thus $(1,1,\ldots,1)$ is an interior point of the space $\exp\Mt^*$, and we have
$$
\dim \exp\Mt^* = |\Sigma|.
$$

Thus we have shown that the space $\Mt^*$ is diffeomorphic to a convex polyhedron of dimension $|\Sigma|$ with some faces removed.

As for the space $\M^*$, we use its representation in the form
$$
\M^* = \Mt^* \cap \left\{{\textstyle \sum} h_i = 0\right\}.
$$
We have
$$
\M^* \approx \exp\Mt^* \cap \left\{{\textstyle \prod} e_i = 1\right\} \approx \exp\Mt^* \cap \left\{{\textstyle \sum} e_i = 1\right\},
$$
where the last diffeomorphism takes place because $\exp\Mt^*$ is a subset of $(0,+\infty)^{\Sigma}$ and is invariant under positive scaling.
It follows that $\M^*$ is diffeomorphic to a convex polyhedron with some faces removed. The polyhedron has dimension $|\Sigma| - 1$ because the hyperplane $\{\sum e_i = 1\}$ contains an interior point $\frac{1}{|\Sigma|}(1,1,\ldots,1)$ of $\exp\Mt^*$.

Now let us proceed to the second half of Proposition \ref{prp:StructM}. In Example \ref{exl:MtBig}, we have $\M^* = \R^\Sigma / \mathbb{L}$, thus $\M^*$ in this case is unbounded provided $|\Sigma| \ge 2$. The convexity of $\M^*$ is equivalent to the convexity of $\Mt^*$. Let us show that $\Mt^*$ is not convex always when it is a proper subset of $\R^\Sigma$, that is always when the set of inequalities \eqref{eqn:Quad} is non-empty. From
$$
\widetilde{h_{jkl}} = \log\widetilde{e_{jkl}}
$$
it follows easily that
\begin{equation}
\label{eqn:QuadEx}
\widetilde{h_{jkl}}(i) = \log(ae^{h_j} + be^{h_k} + ce^{h_l})
\end{equation}
with $a \ge 0,\, b > 0, c > 0$ (we have $a = 0$ if and only if the angle at $i$ in $ikjl$ equals $\pi$). A simple computation shows that the matrix of the second partial derivatives of \eqref{eqn:QuadEx} with respect to $h_j$, $h_k$, $h_l$ is positive semidefinite and non-degenerate. Thus $\widetilde{h_{jkl}}(i)$ is convex and non-linear. Let $(h^0) \in \partial\Mt^*$ be a point where only one of inequalities \eqref{eqn:Quad} holds as an equality. Then the space $\Mt^*$ is non-convex near~$h^0$.

To construct an example of $\M^*$ that is non-convex and unbounded at the same time, one can add a singularity of a small negative curvature inside one of the triangles of Example \ref{exl:UniqueTriang}.
\end{proof}

\begin{rem}
The point $(0,0,\ldots,0) \in \Mt^*$ corresponds to a truncated cusp with equal heights, that is a cusp circumscribed around a horoball (with coparticles). It is not hard to see that the dual tesselation of such a cusp is the Delaunay tesselation constructed in the proof of Lemma \ref{lem:TriangExist}.
\end{rem}

\begin{rem}
If the metric $g$ has a single cone point $i$, then the spaces $\Mt^*$ and $\M^*$  become a real line $\R^1$ and a point $\R^0$, respectively. By \eqref{eqn:SumCurv} we have $\kappa_i = 0$ for all cusps with coparticles with Gauss image $(\Tor^2, g)$. Thus the point into which degenerates the space $\M^*$ is a unique convex polyhedral cusp with the given Gauss image. Its dual tesselation is the Delaunay tesselation of $(\Tor^2, g)$. The cusp itself has a unique face which is either a quadrilateral or a hexagon, and one, respectively, two vertices.
\end{rem}

\subsection{Properties of faces of the dual tesselation}
The next lemma will be used in Section \ref{subsec:BehaveBoundary}.
\begin{lem}
\label{lem:BoundaryM}
A cusp $M$ is an interior point of $\M^*$ if and only if all faces of its dual tesselation are strictly convex spherical polygons.

In particular, if $\kappa_i = 0$ for all $i$, then $M$ is an interior point of $\M^*$.
\end{lem}
\begin{proof}
Assume that $M$ lies on the boundary of $\M^*$. Then any truncation $M_{\mathrm{tr}}$ of $M$ lies on the boundary of $\Mt^*$. Let $h \in \R^\Sigma$ be the heights of $M_{\mathrm{tr}}$. By Proposition \ref{prp:DescrM}, one of the inequalities \eqref{eqn:Quad} holds as equality for $h$. It follows that the support function $\widetilde{h}$ of $M_{\mathrm{tr}}$ is smooth on some quadrilateral $ikjl$ with an angle at least $\pi$ at~$i$. Since the faces of the dual tesselation are regions of smoothness of the support function, there is a face that contains the quadrilateral $ijkl$. Then this face is not a strictly convex polygon.

In the opposite direction, let $F$ be a face with an angle at least $\pi$ at a vertex~$i$. (The vertex $i$ might as well be an isolated vertex surrounded by the face $F$, in this case the angle at $i$ is even bigger than $2\pi$.) Since the domain of a support-like function is a hemisphere, the development of $F$ is contained in a hemisphere. One can show that there are vertices $j$, $k$, $l$ of $F$ such that the quadrilateral $ikjl$ is contained in $F$ and has an angle at least $\pi$ at $i$. It follows that we have an equality in \eqref{eqn:Quad}, thus $M$ lies on the boundary of $\M^*$.

By Lemma \ref{lem:ZeroCurv}, if $\kappa_i = 0$ for all $i$, then $M$ is a convex polyhedral cusp without coparticles. The dual tesselation of a convex polyhedral cusp consists of the Gauss images of the vertices which are convex spherical polygons. Thus $M$ lies in the interior of $\M^*$.
\end{proof}

Now we will investigate non-simply connected faces of the dual tesselation. The next two lemmas will be needed only in the proof of Theorem \ref{thm:LocRig}.

Let $F$ be a face of the dual tesselation. Recall that the support function $\widetilde{h} \co F \to \R$ of a truncated cusp with coparticles has locally the form
\begin{equation}
\label{eqn:FormOfH}
x \mapsto \log \cos \dist(x,a) + b.
\end{equation}
It follows that for every small open set $U \subset F$ the composition of the function $\exp(\widetilde{h})$ with an embedding $U \to \Sph^2 \to \R^3$ is the restriction of a linear function.

Let $\Pi \co \widetilde{F} \to F$ be the universal covering map, and let $D: \widetilde{F} \to \Sph^2$ be the developing map.
By the previous paragraph, we have
\begin{equation}
\label{eqn:Devel}
\exp(\widetilde{h} \circ \Pi) = f \circ D,
\end{equation}
where $f \co \Sph^2 \to \R$ is the restriction of a linear function $\R^3 \to \R$. Let us call the point $\frac{\grad f}{\|\grad f\|}$ the \emph{pole} of $f$. The pole corresponds to the point $a$ in the formula \eqref{eqn:FormOfH}. Note that due to the positivity of the left hand side of \eqref{eqn:Devel} the set $D(\widetilde{F})$ is contained in a hemisphere centered at the pole of $f$. Also, the position of the pole determines the function $f$ up to a constant factor, and the function $\widetilde{h}$ up to a constant summand.

For every covering transformation $\phi_g \co \widetilde{F} \to \widetilde{F}$, there is a unique transformation $\psi_g \in \mathrm{SO}(3)$ such that
\begin{equation}
\label{eqn:Holon}
D \circ \phi_g = \psi_g \circ D.
\end{equation}
The map
$$
\begin{array}{rcl}
\pi_1(F) & \to & \mathrm{SO}(2),\\
g & \mapsto & \psi_g
\end{array}
$$
is called the \emph{holonomy} of $F$.

\begin{lem}
\label{lem:NonTrivHol}
If a face $F$ of the dual tesselation has a non-trivial holonomy, then the support function on $F$ is determined uniquely up to a constant summand.
\end{lem}
\begin{proof}
Let $\psi_g$ be a non-trivial orthogonal transformation. It follows from \eqref{eqn:Devel} and \eqref{eqn:Holon} that the function $f$ is invariant on the orbit of $x \in \Sph^2$ under the action of $\psi_g$, for all $x \in D(\widetilde{F})$. Since $D(\widetilde{F})$ has a non-empty interior, one can easily show that the pole of $f$ is one of the fixed points of the map $\psi_g$. As the hemisphere centered at the pole contains $D(\widetilde{F})$, the pole is uniquely determined. And a support function $\widetilde{h}$ is determined by the pole uniquely up to a constant summand.
\end{proof}

By the \emph{boundary} of a face $F$ we mean the complement $\widehat{F} \setminus F$, where $\widehat{F}$ is the completion of $F$ (rather than taking $\overline{F} \setminus F$, where $\overline{F} \subset \Tor^2$ is the closure of $F$). Therefore the boundary of $F$ consists of cone points and closed polygonal curves.

\begin{lem}
\label{lem:TrivHol}
If a face $F$ of the dual tesselation has a trivial holonomy, then at least one of its boundary components contains three non-collinear vertices. In particular, the support function on $F$ is uniquely determined by its values at this boundary component.
\end{lem}
\begin{proof}
If the holonomy of $F$ is trivial, then the developing map $D \co \widetilde{F} \to \Sph^2$ descends to a map $\overline{D} \co F \to \Sph^2$. If the boundary of $F$ consists of cone points only, then $F = \Tor^2 \setminus \Sigma$, and $\overline{D}$ can be extended to a branched covering $\Tor^2 \to \Sph^2$. But this contradicts to the fact that $\overline{D}(F)$ is contained in a hemisphere.

The set $\overline{D}(F) \subset \Sph^2$ is bounded by a finite collection of polygonal curves (which need not be the images of the boundary components of $F$, since those images can intersect each other and themselves). Since $\overline{D}(F)$ is contained inside a hemisphere, it has a vertex with an angle less than $\pi$. Then this vertex is the image of a vertex $j$ of a boundary component $K$ of $F$. It follows that $j$ and its neighbors $i$ and $k$ on $K$ are all different and don't lie on a big circle when mapped to $\Sph^2$. The lemma is proved.
\end{proof}

\section{A concave function $V$}
\label{sec:Functional}

\subsection{Smooth functions on $\M^*$ and Whitney's extension theorem}
Recall that we identify the spaces $\Mt^*$ and $\M^*$ with subsets of Euclidean spaces $\R^\Sigma$ and $\R^{|\Sigma|-1}$, see the end of Section \ref{subsec:HeightsToCusp}.

\begin{dfn}
Let $X$ be a closed subset of $\R^n$. A function $f \co X \to \R$ is said to belong to $C^m(X)$ if there exist an open set $Y \supset X$ and a function $g \in C^m(Y)$ such that $f = g|_X$.
\end{dfn}

\begin{lem}
\label{lem:CmPr}
A function $f \co \M^* \to \R$ is of class $C^m$ if and only if $f \circ \pr$ is of class $C^m$ on $\Mt^*$. Here $\pr \co \R^\Sigma \to \R^\Sigma / \mathbb{L}$ is the projection along the one-dimensional space spanned by $(1,1,\ldots,1)$.
\end{lem}
\begin{proof}
If $g$ is a smooth extension of $f$, then $g \circ \pr$ is a smooth extension of $f \circ \pr$. Conversely, if $\overline{g}$ extends $f \circ \pr$ smoothly, then identify $\M^*$ with a section of $\Mt^*$ by a hyperplane and consider the restriction of $\overline g$.
\end{proof}

We will need the following smoothness criterion which is a special case of the Whitney extension theorem, \cite{Whi34, KP99}.

\begin{prp}
\label{prp:Whi}
Let $X$ be a closed subset of $\R^\Sigma$. A function $f \co X \to \R$ is of class $C^2$ if and only if there exist functions $f_i$ and $f_{ij}$ on $X$ such that the following three conditions hold:
\begin{equation}
\label{eqn:VDiff1}
f(h+x) = f(h) + \sum_{i \in \Sigma} f_i(h) x_i + \sum_{i,j \in \Sigma} \frac{f_{ij}(h)}{2} x_i x_j + R_3(h,x),
\end{equation}
where
$$
\lim_{\|x\| \to 0} \frac{R_3(h,x)}{\|x\|^2} = 0
$$
holds uniformly on compact subsets of $X$;

\bigskip

\noindent for all $i \in \Sigma$,
\begin{equation}
\label{eqn:VDiff2}
f_i(h+x) = f_i(h) + \sum_{j \in \Sigma} f_{ij}(h) x_j + R_2(h,x), 
\end{equation}
where
$$
\lim_{\|x\| \to 0} \frac{R_2(h,x)}{\|x\|} = 0
$$
holds uniformly on compact subsets of $X$;
\begin{equation}
\label{eqn:VDiff3}
\mbox{for all }i,j \in \Sigma,\mbox{ the function }f_{ij}\mbox{ is continuous.}
\end{equation}
\end{prp}

\subsection{Schl\"afli formula for bricks}
\label{subsec:Bricks}
\begin{figure}[ht]
\begin{center}
\input{Fig/Brick.tex}
\end{center}
\caption{Angles and lengths in a brick.}
\label{fig:Brick}
\end{figure}
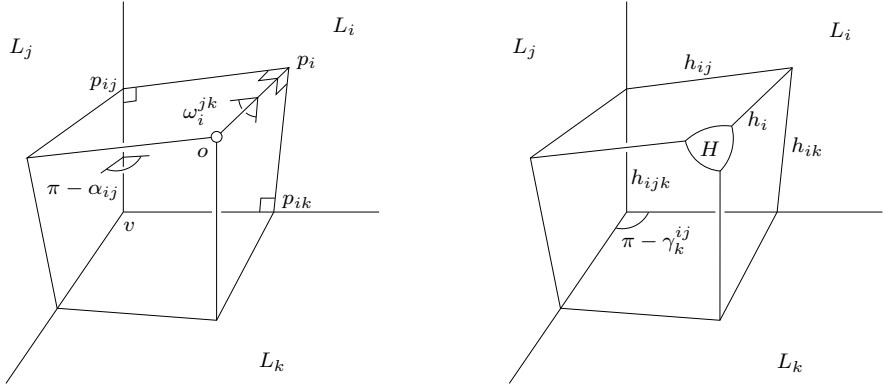

Let $M$ be a convex polyhedral cusp with coparticles. Choose a triangulation $T$ of its dual tesselation and represent $M$ as a complex of corners. For a corner $(L_i, L_j, L_k; o)$, $ijk \in T$, let $p_i$ be the foot of the perpendicular dropped from $o$ to $L_i$, and let $p_{ij}$ be the foot of the perpendicular dropped from $p_i$ to the line $L_i \cap L_j$. Define similarly $p_\cdot$ and $p_{\cdot\cdot}$ for all other indices and pairs of indices. Note that $p_{ij} = p_{ji}$. Call the polyhedron with vertices $o$, $p_i$, $p_j$, $p_k$, $p_{ij}$, $p_{jk}$, $p_{ki}$, $v = L_i \cap L_j \cap L_k$ a \emph{brick} and denote it by $B_{ijk}$. See Figure \ref{fig:Brick}, left. A brick is combinatorially equivalent to a cube, but its boundary can have self-intersections (see the discussion at the beginning of Section \ref{sec:CuspDef}).

Recall that $h_i, h_j, h_k$ are lengths of truncated edges of the brick $B_{ijk}$. Denote the lengths of the other edges by
\begin{eqnarray*}
h_{ij} & := & \dist(p_i, p_{ij}),\\
h_{ijk} & := & \dist(p_{ij}, v),
\end{eqnarray*}
see Figure \ref{fig:Brick}, right. The distances are signed, according to whether the point $p_i$, respectively $p_{ij}$, lies on the positive or on the negative side of the plane $L_j$, respectively $L_k$.  Finally, let $\alpha_{ij}$, $\gamma_i^{jk}$ be the length of the side $ij$ and the value of the angle at the vertex $i$ in the Gauss image of the corner. These are related to planar and dihedral angles of the brick at the vertex $v$, see Figure \ref{fig:Brick}.

A convenient way to look at a brick is to represent it as an algebraic sum of six orthoschemes with an ideal vertex: $o p_i p_{ij} v$ and so on.

\begin{dfn}
The volume of a brick is the algebraic sum of the volumes of its constituting orthoschemes.
\end{dfn}

Denote by $\omega_i^{jk}$ the (appropriately chosen) dihedral angle at the edge $op_i$ of the brick $B_{ijk}$. That is, $\omega_i^{jk}$ is the angle at the vertex $q_i$ in the link of the brick at $o$, see Section \ref{subsec:CuspCurv}. We have
\begin{equation}
\label{eqn:Om1}
\omega_i^{jk} + \omega_j^{ik} + \omega_k^{ij} = \pi,
\end{equation}
\begin{equation}
\label{eqn:Om2}
\omega_i = \sum_{ijk \in T} \omega_i^{jk}.
\end{equation}

\begin{lem}
For a brick $B_{ijk}$ with a fixed Gauss image, we have
\begin{equation}
\label{eqn:VolBrick}
d\vol(B_{ijk}) = -\frac{1}{2}(h_i d\omega_i^{jk} + h_j d\omega_j^{ik} + h_k d\omega_k^{ij}),
\end{equation}
where $h_i$, $h_j$, $h_k$ are heights with respect to an arbitrary truncation.
\end{lem}
\begin{proof}
The Schl\"afli formula for a compact 3-dimensional hyperbolic polyhedron $P$ says
$$
d\vol(P) = -\frac{1}{2}\sum_e \ell_e d\theta_e.
$$
The sum is over all edges of the polyhedron, $\ell_e$ is the length of the edge $e$, $\theta_e$ is the dihedral angle at $e$. If $P$ has one or several ideal vertices, then by \cite[page 294]{Mil94}, \cite[Theorem 14.5]{Riv94} the same formula holds with $\ell_e$ denoting the truncated length of $e$ with respect to arbitrary truncations.

Dihedral angles at the edges $vp_{ij}$, $vp_{jk}$, $vp_{ki}$ of the brick $B_{ijk}$ are determined by its Gauss image, therefore constant. Dihedral angles at the 6 edges $p_ip_{ij}, \ldots, p_kp_{kj}$ are equal to $\frac{\pi}{2}$, also constant. Hence the Schl\"afli formula for $B_{ijk}$ yields \eqref{eqn:VolBrick}.

If the brick has self-intersecting boundary, then \eqref{eqn:VolBrick} follows from summing up the Schl\"afli formulas for its constituting orthoschemes.
\end{proof}

\subsection{Definition of $V$}
\label{subsec:DefV}

\begin{dfn}
Let $M$ be a convex polyhedral cusp with coparticles. Put
\begin{equation}
\label{eqn:DefV}
V(M) = -2\vol(M) +\sum_{i \in \Sigma} h_i\kappa_i.
\end{equation}
Here $\vol(M)$ is defined as the sum of volumes of bricks, $h_i$ are heights with respect to an arbitrary truncation of $M$, and $\kappa_i = 2\pi - \omega_i$ are curvatures of the coparticles.
\end{dfn}

\begin{lem}
$V(M)$ is well-defined.
\end{lem}
\begin{proof}
Let us show that both summands on the right hand side of \eqref{eqn:DefV} are well-defined.

The decomposition of $M$ into bricks depends on the choice of a triangulation $T$ that subdivides the dual tesselation. To show that $\vol(M)$ is well-defined, cut the bricks into orthoschemes. Then for every edge $\ell_{ij}$ of zero length there are two equal orthoschemes whose volumes are counted with opposite signs. After eliminating all of such pairs, their remains a set of orthoschemes which is independent of the choice of triangulation $T$.

A change of truncation results in adding a constant to all of $h_i$. Due to \eqref{eqn:SumCurv}, this does not change the sum $\sum_i h_i\kappa_i$.
\end{proof}

\begin{prp}
\label{prp:FuncV}
The function $V$ defined by \eqref{eqn:DefV} belongs to $C^2(\M^*)$. Moreover,
\begin{equation}
\label{eqn:VDer}
\frac{\partial V}{\partial h_i} = \kappa_i
\end{equation}
holds everywhere in $\M^*$.
\end{prp}

Let us explain the meaning of $\frac{\partial V}{\partial h_i}$ at the boundary points of $\M^*$. By Proposition \ref{prp:StructM}, the set $\M^*$ is diffeomorphic to a convex full-dimensional polyhedron with some faces removed. Thus, every point $x \in \M^*$ possesses a full-dimensional tangent cone $C_x\M^*$. From an analog of Proposition \ref{prp:Whi} for $C^1$--functions, it follows that for every $f \in C^1(\M^*)$ and $\xi \in C_x\M^*$, the directional derivative $\frac{\partial f}{\partial \xi}(x)$ is well-defined and depends linearly on $\xi$. By linearity, this allows to define $\frac{\partial f}{\partial \xi}(x)$ for all $\xi \in \R^{|\Sigma|-1}$. The expression $\frac{\partial V}{\partial h_i}$ denotes the directional derivative with respect to $e_i + \mathbb{L}$, where $e_i$ is the $i$th basis vector of $\R^{\Sigma}$.

\begin{proof}[Proof of Proposition \ref{prp:FuncV}]
We will prove that
$$
\widetilde{V} = V \circ \pr \in C^2(\Mt^*).
$$
By Lemma \ref{lem:CmPr}, this implies $V \in C^2(\M^*)$. For \eqref{eqn:VDer}, it suffices to show that in an expansion of $\widetilde{V}$ according to Proposition \ref{prp:Whi}, $\widetilde{V}_i$ can be put equal to~$\kappa_i$.

Consider the decomposition
\begin{equation}
\label{eqn:MtDec}
\Mt^* = \bigcup_T \Mt^{*,T},
\end{equation}
where $\Mt^{*,T}$ denotes the set of truncated cusps that have a representative of the form $(T,h)$. In other words, $M$ belongs to $\Mt^{*,T}$ if the dual tesselation of its boundary can be refined to the triangulation $T$.

By Lemma \ref{lem:VDiffT}, the restriction of $\widetilde{V}$ to $\Mt^{*,T}$ is of class $C^\infty$ for all $T$. We will show that these $C^\infty$-patches fit together in a $C^2$-way.

Denote the second partial derivatives \eqref{eqn:DerK1}, \eqref{eqn:DerK2} of $\widetilde{V}$ on $\Mt^{*,T}$ by $\widetilde{V}_{ij}^T$. We claim that
\begin{equation}
\label{eqn:DerIndep}
\widetilde{V}_{ij}^T(h) = \widetilde{V}_{ij}^{T'}(h)
\end{equation}
holds for all $i,j \in \Sigma$ if $h \in \Mt^{*,T} \cap \Mt^{*,T'}$. If $ij \in T$ and $ij \in T'$, then similarly to the proof of Lemma \ref{lem:lijWD} in Section \ref{subsec:SuppFunc} one can show that $h_{ijk}$ and $h_{ijl}$ are well-defined (although the vertices $k$ and $l$ opposite to the edge $ij$ can differ for $T$ and $T'$). Thus \eqref{eqn:DerIndep} holds in this case. If $ij \in T$ but $ij \notin T'$, then we have
$$
h_{ijk} + h_{ijl} = \ell_{ij} = 0.
$$
Therefore $\tanh h_{ijk} + \tanh h_{ijl} = 0$ which implies
$$
\widetilde{V}_{ij}^T(h) = 0 = \widetilde{V}_{ij}^{T'}(h).
$$
Finally, the equality $\widetilde{V}_{ii}^T(h) = \widetilde{V}_{ii}^{T'}(h)$ follows from the equalities \eqref{eqn:DerIndep} for $j \ne i$ due to \eqref{eqn:DerK2}.

It follows that the formulas \eqref{eqn:DerK1} and \eqref{eqn:DerK2} define continuous functions $\widetilde{V}_{ij}$ on $\Mt^{*,T}$. The functions $\widetilde{V}_i = \kappa_i$ are well-defined and continuous due to Lemma \ref{lem:LinkWD}. Therefore at any point $h \in \Mt^*$ that belongs to several of $\Mt^{*,T}$, the Taylor expansions of $\widetilde{V}$ and $\kappa_i$ on $\Mt^{*,T}$ up to quadratic and linear terms respectively fit together to expansions that satisfy conditions of Proposition \ref{prp:Whi}. The uniform convergence of remainders follows from the finiteness of the decomposition \eqref{eqn:MtDec}, see Lemma \ref{lem:TriangFinite}. Thus $\widetilde{V}$ is of class $C^2$ on $\Mt^*$ with partial derivatives $\kappa_i$.
\end{proof}

\begin{rem}
In general, $V$ is not of class $C^3$.
\end{rem}

\begin{lem}
\label{lem:VDiffT}
For every triangulation $T$ of $(\Tor^2,g)$, the function $\widetilde{V}$ is of class $C^\infty$ on $\Mt^{*,T}$. Furthermore, an extension of $\widetilde{V}$ to a neighborhood of $\Mt^{*,T}$ can be chosen so that its partial derivatives are
\begin{equation}
\label{eqn:DerV}
\frac{\partial \widetilde{V}}{\partial h_i} = \kappa_i;
\end{equation}

\begin{equation}
\label{eqn:DerK1}
\frac{\partial^2 \widetilde{V}}{\partial h_i \partial h_j} = \begin{cases}
\frac{\tanh h_{ijk} + \tanh h_{ijl}} {\sin\alpha_{ij} \cosh h_{ij} \cosh h_{ji}} & \mbox{ for }i \ne j \mbox{ and } ij \in T;\\
0 & \mbox{ for }i \ne j \mbox{ and } ij \notin T;
\end{cases}
\end{equation}

\begin{equation}
\label{eqn:DerK2}
\frac{\partial^2 \widetilde{V}}{\partial h_i^2} = - \sum_{j \ne i} \frac{\partial^2 \widetilde{V}}{\partial h_i \partial h_j}.
\end{equation}
\end{lem}

\begin{rem}
In the case when the vertices $i$ and $j$ are joined by several edges, on the right hand side of the formula \eqref{eqn:DerK1} one should take the sum over \emph{all} edges between $i$ and $j$. It can be shown that loop edges don't have any effect on the computation of derivatives.
\end{rem}

\begin{proof}[Proof of Lemma \ref{lem:VDiffT}]
The volume of a cusp $M$ is the algebraic sum of the volumes of constituting orthoschemes. The volume of such an orthoscheme is a $C^\infty$-function of its dihedral angles, see \cite[Sections 3.3 and 3.4 of Chapter 7]{AVS93}. The dihedral angles of the orthoscheme $op_ip_{ij}v$ on Figure \ref{fig:Brick} can be expressed through edge lengths $h_i$, $h_{ij}$, $h_{ijk}$ using formulas of hyperbolic and spherical trigonometry. As a result, the volume of each orthoscheme is a $C^\infty$-function of $h$. Similarly, $\kappa_i$ for each $i$ is a $C^\infty$-function of $h$.

Note that the formulas for the volumes of orthoschemes and for $\kappa_i$ make sense for any values of $h$, not only for those contained in $\Mt^{*,T}$. This automatically gives a $C^\infty$-extension of $\widetilde{V}$ to a neighborhood of $\Mt^{*,T}$.

To prove the formula \eqref{eqn:DerV}, sum up the equations \eqref{eqn:VolBrick} over all bricks coming from triangulation $T$. We obtain
$$
d\vol(M) = \frac{1}{2}\sum_{i \in \Sigma} h_i d\kappa_i.
$$
It follows that for the constructed extension of $\widetilde{V}$ on a neighborhood of $\Mt^{*,T}$ we have
$$
d\widetilde{V}(M) = -2 \cdot d\vol(M) + \sum_{i \in \Sigma} h_i d\kappa_i + \sum_{i \in \Sigma} \kappa_i dh_i = \sum_{i \in \Sigma} \kappa_i dh_i,
$$
which implies \eqref{eqn:DerV}.

Formulas \eqref{eqn:DerK1} follow from \eqref{eqn:DerOmegaH}. Formula \eqref{eqn:DerK2} follows from the invariance of $\kappa_i$ under the change of truncation: $\kappa_i(h+c\mathbf{1}) = \kappa_i(h)$.
\end{proof}

\begin{rem}
The space $\Mt^{*,T}$ can have dimension less than $|\Sigma|$. Therefore partial derivatives of $\widetilde{V}$ on $\Mt^{*,T}$ are not well-defined, but depend on an extension of $\widetilde{V}$ to a neighborhood of $\Mt^{*,T}$. The formulas of Lemma \ref{lem:VDiffT} hold for the ``most natural'' extension.
\end{rem}

\subsection{The Hessian of $V$ and rigidity of cusps}
\label{subsec:LocRig}

The following is a folklore lemma.
\begin{lem}
\label{lem:PosDef}
Let $A = (a_{ij})_{i,j=1}^n$ be a symmetric $n \times n$ matrix such that
\begin{eqnarray*}
a_{ij} & \le & 0 \mbox{ for }i \ne j;\\
a_{ii} & = & \sum_{j \ne i} a_{ij}.
\end{eqnarray*}
Then the quadratic form $\sum_{i,j} a_{ij}x_i x_j$ is positive semidefinite.

Define the \emph{underlying graph} of $A$ as a graph on the vertex set $\{1,\ldots,n\}$ where $i$ and $j$ are joined by an edge if $a_{ij} \ne 0$. Then the kernel of $A$ consists of vectors $x$ such that $x_i$ is constant over every connected component of the underlying graph of~$A$.
\end{lem}
\begin{proof}
We have
$$
\sum_{i,j} a_{ij} x_i x_j = -\sum_{i < j} a_{ij}(x_i - x_j)^2 \ge 0 \mbox{ for all }x.
$$
Thus the quadratic form is positively semidefinite.

A vector $x$ lies in the kernel iff $x_i=x_j$ for all $i$ and $j$ such that $a_{ij} \ne 0$ that is for all $i$ and $j$ joined by an edge of the underlying graph. This implies the second assertion of the lemma.
\end{proof}

\begin{prp}
\label{prp:VConc}
The function $V$ is concave. Its Hessian at a cusp with coparticles $M$ is negative definite if and only if the graph of the dual tesselation of $M$ is connected.
\end{prp}
\begin{proof}
It suffices to show that $\widetilde{V}$ is concave and has Hessian of corank 1 if and only if the graph of the dual tesselation is connected.

First and second partial derivatives of $\widetilde{V}$ are given by formulas \eqref{eqn:VDer}, \eqref{eqn:DerK1}, and \eqref{eqn:DerK2}. From
$$
h_{ijk} + h_{ijl} = \ell_{ij} \ge 0,
$$
it follows that $\tanh h_{ijk} + \tanh h_{ijl} \ge 0$. Thus we have
\begin{equation}
\label{eqn:DerPos}
\frac{\partial^2 \widetilde{V}}{\partial h_i \partial h_j} \ge 0 \mbox{ for } i \ne j.
\end{equation}
Due to this and \eqref{eqn:DerK2}, the Hessian matrix of $-\widetilde{V}$ satisfies the assumptions of Lemma \ref{lem:PosDef}. Thus $\widetilde{V}$ is a concave function.

The inequality in \eqref{eqn:DerPos} is strict if and only if $\ell_{ij} > 0$. Therefore the underlying graph of $\big(\frac{\partial^2 \widetilde{V}}{\partial h_i \partial h_j}\big)$ is the graph of the dual tesselation of $M$. Thus the kernel of the Hessian consists of the multiples of $(1,1,\ldots,1)$ if and only if the graph of the dual tesselation is connected.
\end{proof}

\begin{dfn}
\label{dfn:InfRig}
A convex polyhedral cusp with coparticles $M$ is called \emph{infinitesimally rigid} if the matrix $\big(\frac{\partial \kappa_i}{\partial h_j}\big)$ has corank 1.

In other words, a cusp is infinitesimally rigid if any non-trivial first-order change of its heights, the Gauss image of the cusp being fixed, leads to a non-zero first-order change of its curvatures.
\end{dfn}

Due to \eqref{eqn:DerV}, Proposition \ref{prp:VConc} implies that $M$ is infinitesimally rigid if and only if the graph of its dual tesselation is connected. Example \ref{exl:NSCFace} shows that there exist infinitesimally flexible cusps with coparticles.

For cusps without coparticles, there is a different notion of infinitesimal rigidity. A convex polyhedral cusp $M$ is called infinitesimally rigid with respect to its Gauss image, if any first-order change of the hyperbolic metric inside $M$ leads to a non-zero first order change of the metric of the Gauss image. It can be shown that a convex polyhedral cusp $M$ is infinitesimally rigid in this sense if and only if it is infinitesimally rigid in the sense of Definition \ref{dfn:InfRig}, that is as a cusp with coparticles with fixed Gauss image.

\begin{dfn}
A convex polyhedral cusp with coparticles $M$ is called \emph{locally rigid} if for any smooth family of cusps $M(t)$, $t \in [0,1]$ such that
\begin{itemize}
\item $M(0)=M$;
\item all $M(t)$ have the same Gauss image;
\item for all $i \in \Sigma$, the curvatures $\kappa_i(t)$ of $M(t)$ are constant,
\end{itemize}
the cusp $M(t)$ is isometric to $M(0)$ for all $t$.
\end{dfn}

Recall that Theorem \ref{thm:LocRig} in Section \ref{subsec:PresPap} states that cusps without coparticles are infinitesimally rigid, and cusps with coparticles are locally rigid.

\begin{proof}[Proof of Theorem \ref{thm:LocRig}]
As noted in Section \ref{subsec:GaussImage}, faces of the dual tesselation of a cusp without coparticles are convex polygons. Therefore the 1--skeleton of the dual tesselation is connected, and the corank of $(\frac{\partial^2 \widetilde{V}}{\partial h_i \partial h_j})$ at a critical point equals 1 by Proposition \ref{prp:VConc}. Thus, cusps without coparticles are infinitesimally rigid.

Now let $M(t)$ be a smooth family of cusps with coparticles with Gauss image $(\Tor^2,g)$.

Let $G(t)$ be the 1--skeleton of the dual tesselation of $\partial M(t)$. Since $\ell_{ij} > 0$ is an open condition, the map $t \mapsto G(t)$ is lower semi-continuous:
\begin{equation}
\label{eqn:LSC}
\forall t_0 \ \exists \epsilon > 0 \mbox{ such that } G(t) \supset G(t_0) \ \forall t \in (t_0-\epsilon, t_0+\epsilon).
\end{equation}
Let $(a,b) \subset [0,1]$. Choose $t_0 \in (a,b)$ so that the number of edges of $G(t_0)$ is maximal over all $t \in (a,b)$. Then \eqref{eqn:LSC} implies $G(t) = G(t_0)$ for all $t \in (t_0-\epsilon, t_0+\epsilon)$. By Lemma \ref{lem:SameDual}, this implies $M(t) = M(t_0)$ for all $t \in (t_0-\epsilon, t_0+\epsilon)$.

Since the interval $(a,b)$ can be chosen arbitrarily, $M(t)$ is constant on every connected component of a dense open subset of $[0,1]$. It follows that $M(t)$ is constant over all of $[0,1]$.
\end{proof}

\begin{lem}
\label{lem:SameDual}
Let $M(t)$, $t \in (t_0 - \epsilon, t_0 + \epsilon)$ be a smooth family of cusps with the same Gauss images, same curvatures, and same dual tesselations. Then all of $M(t)$ are isometric.
\end{lem}
\begin{proof}
Choose truncations of all of $M(t)$ so that the heights $h_i(t)$ depend smoothly on $t$. We have to show that
\begin{equation}
\label{eqn:hihjConst}
h_i(t) - h_j(t) = \mathrm{const}
\end{equation}
holds for all $i$ and $j$.

Let $G$ be the graph of the dual tesselation. By Lemma \ref{lem:PosDef}, applied to the Hessian of $\widetilde{V}$, we have $\dot{h_i}(t) = \dot{h_j}(t)$ and consequently \eqref{eqn:hihjConst} whenever $i$ and $j$ belong to the same connected component of $G$. If we prove that \eqref{eqn:hihjConst} also holds when $i$ and $j$ lie in different components of $G$ but on the boundary of the same face $F$ of the dual tesselation, then \eqref{eqn:hihjConst} holds for all $i$ and $j$.

Thus let $F$ be a non-simply connected face of $M(t)$. Now we apply Lemmas \ref{lem:NonTrivHol} and \ref{lem:TrivHol}. Recall that the support function $\widetilde{h(t)}$ of the truncation of $M(t)$ is smooth on the face $F$, and that we have $\widetilde{h(t)}(i) = h_i(t)$ for all cone points $i$. If the holonomy of $F$ is non-trivial, then by Lemma \ref{lem:NonTrivHol} the restriction of $\widetilde{h(t)}$ to $F$ is determined uniquely up to a constant summand. Thus in this case \eqref{eqn:hihjConst} holds for $i$ and $j$ on the boundary of $F$. If the holonomy of $F$ is trivial, then by Lemma \ref{lem:TrivHol} there exists a component $K$ of the boundary of $F$ such that the values of $\widetilde{h(t)}$ at the vertices of $K$ determine $\widetilde{h(t)}|_F$ uniquely. Besides, adding a constant to $\widetilde{h(t)}|_K$ adds a constant to $\widetilde{h(t)}|_F$. Thus \eqref{eqn:hihjConst} holds for all $i$ and $j$ on the boundary of $F$.
\end{proof}

\section{Proof of the main theorem}
\label{sec:Proof}

\subsection{Morse theory on manifolds with corners}
A \emph{convex polyhedral cone} is an intersection of finitely many half-spaces whose boundary hyperplanes pass through the origin.
\begin{dfn}
A manifold with corners is a topological manifold $X$ equipped with an atlas $\{(U_\alpha,\phi_\alpha)\}$, where $U_\alpha$ is an open subset of $X$, and $\phi_\alpha$ is a homeomorphism from $U_\alpha$ to an open subset of a convex polyhedral cone. The transition maps between pairs of charts are assumed to be $C^\infty$.
\end{dfn}

\begin{rem}
Quite often one considers a more restricted class of manifolds with corners, namely those locally modelled on a \emph{simple} convex polyhedral cone $[0,+\infty)^n$.
\end{rem}

Any convex polyhedron, and any convex polyhedron with some faces removed, is a manifold with corners. Thus, by Proposition \ref{prp:StructM}, the spaces $\M^*$ and $\Mt^*$ are manifolds with corners.

It is possible to define the tangent space $T_xX$ at every point $x$ of a manifold with corners $X$. If $X$ is smoothly embedded in a smooth manifold $M$, then $T_xX$ can be identified with $T_xM$. This suffices for our needs, since $\M^*$ and $\Mt^*$ are realized as subsets of Euclidean spaces. Any $C^1$--function on $X$ has a well-defined differential $df_x \in (T_xX)^*$ at every point $x \in X$.

The \emph{tangent cone} $C_xX \subset T_xX$ consists of the vectors tangent to smooth curves in $X$. The cone $C_xX$ is a convex polyhedral cone in the vector space $T_xX$.

The $k$--dimensional \emph{stratum} of $X$ consists of all points that are mapped by some $\phi_\alpha$ into the interior of a $k$-dimensional face of a polyhedral cone. The $k$--dimensional stratum is a smooth $k$--manifold, and a manifold with corners is a disjoint union of its strata. The tangent space $T_xS$ to a stratum $S$ is the maximal linear subspace of the cone $C_xX$.

\begin{dfn}
\label{dfn:CritPoint}
Let $f \co X \to \R$ be a $C^1$--function on a manifold with corners~$X$. A point $x \in X$ is called a \emph{critical point} of $f$, if
$$
\frac{\partial f}{\partial \xi}(x) \ge 0 \mbox{ for every } \xi \in C_xX.
$$
\end{dfn}

If $x$ is an interior point of $X$, then it is critical if and only if the differential of $f$ vanishes at $x$. On the contrary, if $x$ is a boundary point, then the vanishing of $df_x$ is sufficient but not necessary for $x$ to be critical. Also note that if $x$ is a critical point of function $f$, then $x$ is not necessarily a critical point of $-f$.

\begin{dfn}
A critical point $x$ of a $C^2$--function $f \co X \to \R$ is called \emph{non-degenerate}, if the following two conditions hold:
\begin{enumerate}
\item the Hessian of $f|_S$ is non-degenerate at $x$, where $S$ is the stratum of $X$ containing $x$;
\item $\frac{\partial f}{\partial \xi} > 0$ holds for every vector $\xi \in C_xX$ not tangent to the stratum~$S$.
\end{enumerate}
The \emph{index} of a non-degenerate critical point $x$ is the index of $x$ viewed as a critical point of $f|_S$. That is, the index of $x$ is the number of negative eigenvalues of the Hessian of $f|_S$.
\end{dfn}
If $x$ is a boundary point of $X$ such that $df_x=0$, then $x$ is a degenerate critical point.

A function $f \in C^2(X)$ is called a \emph{Morse function}, if all of its critical points are non-degenerate.

Similarly to the classical Morse theory, the following theorems hold.

\begin{thm}[{\cite[Appendix B, Theorem B.4.]{Far04}}]
\label{thm:Morse1}
Let $X$ be a manifold with corners, and let $f \in C^1(X)$. Suppose that the set $f^{-1}([a,b]) \subset X$ is compact and contains no critical points of $f$. Then $M^a = f^{-1}((-\infty,a])$ is a deformation retract of $M^b = f^{-1}((-\infty,b])$.
\end{thm}

In \cite{Far04}, only \emph{simple} manifolds with corners (those modelled on $[0,+\infty)^n$) are considered. Nevertheless, the proof of Theorem \ref{thm:Morse1} given there works also for general manifolds with corners.

\begin{thm}[{\cite[Appendix B, Theorem B.5.]{Far04}}]
\label{thm:Morse2}
Let $X$ be a manifold with corners, and let $f \in C^2(X)$. Suppose that the set $f^{-1}([a,b]) \subset X$ is compact and contains a single critical point $p$ which lies in its interior, is non-degenerate, and has index $\lambda$. Then the manifold $M^b = f^{-1}((-\infty,b])$ is homotopy equivalent to $M^a \cup e^\lambda$, the result of gluing a cell of dimension $\lambda$ to $M^a = f^{-1}((-\infty,a])$.
\end{thm}

\begin{cor}[{\cite[Appendix B, Corollary B.6.]{Far04}}]
\label{cor:Morse3}
Let $f \in C^2(X)$ be a Morse function on a compact manifold with corners. Then $X$ is homotopy equivalent to a CW-complex with cells of each dimension $\lambda$ in one-to-one correspondence with the critical points of $f$ of index $\lambda$.
\end{cor}

Theorem \ref{thm:Morse2} is proved in \cite{Far04} for Morse-Bott functions on simple manifolds with corners (i.~e. function with non-degenerate critical submanifolds). Being restricted to Morse functions, the argument from \cite{Far04} works also for general manifolds with corners.

Note also that we will use only a very light version of Theorem \ref{thm:Morse2}, namely when the index of a critical point is equal to $0$.

\subsection{Proof of Theorem \ref{thm:Main}}
By Proposition \ref{prp:FuncV} and Lemma \ref{lem:ZeroCurv}, convex polyhedral cusps with Gauss image $(\Tor^2,g)$ are in one-to-one correspondence with those points in $\M^*$ where the differential of function $V$ vanishes. Therefore to establish Theorem \ref{thm:Main} it suffices to prove the following proposition.

\begin{prp}
The function $-V$ has a unique critical point (in the sense of Definition \ref{dfn:CritPoint}), and this point lies in the interior of $\M^*$.
\end{prp}
\begin{proof}
For any $d > 0$, put
$$
Q_d = \{h \in \R^\Sigma\ |\ |h_i - h_j| \le d \mbox{ for all } i,j \in \Sigma\},
$$
and denote by $[Q_d]$ the projection of $Q_d$ on $\R^\Sigma / \mathbb{L}$.

By Proposition \ref{prp:VInf} applied to $d(h) = D$, the restriction $-V|_{\M^* \cap [Q_D]}$ has no critical points on the boundary. Thus if we show that $-V|_{\M^* \cap [Q_D]}$ has a unique critical point $h^0$, then $h^0$ is also critical for the function $-V$ on $\M^*$. Besides, $h^0$ lies in the interior of $\M^*$. Proposition \ref{prp:VInf} also implies that the function $-V$ has no critical points outside $[Q_D]$. Thus $h^0$ is the unique critical point of $-V$.

So, let us consider the restriction of $-V$ to $\M^* \cap [Q_D]$. By Lemma \ref{lem:MQ}, the space $\M^* \cap [Q_D] \subset \R^\Sigma / \mathbb{L}$ is a compact full-dimensional manifold with corners homeomorphic to a ball. Every interior critical point of $-V$ corresponds to a cusp without coparticles, by Lemma \ref{lem:ZeroCurv}. Since the graph of the dual tesselation of a cusp without coparticles is connected, Proposition \ref{prp:VConc} implies that the Hessian of $-V$ at every interior critical point is positive definite. Thus $-V|_{\M^* \cap [Q_D]}$ is a Morse function on a compact contractible manifold with corners, and all critical points of $-V|_{\M^* \cap [Q_D]}$ have index~$0$. By Corollary \ref{cor:Morse3}, the critical point is unique. The proposition is proved.
\end{proof}

\begin{lem}
\label{lem:MQ}
For all $d > 0$, the space $\M^* \cap [Q_d] \subset \R^\Sigma / \mathbb{L}$ is diffeomorphic to a compact convex polyhedron of dimension $|\Sigma|-1$.
\end{lem}
\begin{proof}
Consider the space $\Mt^* \cap Q_d$. As in the proof of Proposition \ref{prp:StructM}, consider its image $\exp(\Mt^* \cap Q_d)$ under the diffeomorphism \eqref{eqn:Diffeo}. The set $\exp(\Mt^* \cap Q_d) \subset \R^\Sigma$ is the solution set of a system of linear inequalities \eqref{eqn:Ineq1}, \eqref{eqn:Ineq2}, and
\begin{equation}
\label{eqn:Ineq3}
e_i \le C e_j \quad \mbox{for all }i,j \in \Sigma,
\end{equation}
where $C = e^d > 1$. It is easy to see that in the presence of \eqref{eqn:Ineq3} the inequalities \eqref{eqn:Ineq2} can be replaced by $e \ne (0,0,\ldots,0)$. Thus $\exp(\Mt^* \cap Q_d)$ is a closed convex polyhedral cone contained in $(0,+\infty)^\Sigma$, with the apex removed. The space $M^* \cap [Q_d]$ is diffeomorphic to the section of $\exp(\Mt^* \cap Q_d)$ by the hyperplane $\{\sum h_i = 1\}$. Thus $M^* \cap [Q_d]$ is diffeomorphic to a compact convex polyhedron. Since the point $(1,1,\ldots,1)$ lies in the interior of $\Mt^* \cap Q_d$, the dimension of $\dim (M^* \cap [Q_d])$ equals $|\Sigma|-1$.
\end{proof}

\subsection{Behavior of $V$ on the boundary}
\label{subsec:BehaveBoundary}
The next proposition will be used in the proof of Proposition \ref{prp:VInf}.
\begin{prp}
\label{prp:CritPtBdry}
Function $-V$ has no critical points on the boundary of $\M^*$.
\end{prp}
\begin{proof}
Let $M = [h]$ be a convex polyhedral cusp with coparticles such that $M \in \partial\M^*$. Here $h \in \R^\Sigma$ is the vector of heights of a truncation of $M$, and $[h] = h + \mathbb{L} \in \R^\Sigma/\mathbb{L}$. By Lemma \ref{lem:BoundaryM}, the dual tesselation of $M$ has at least one face which is not a strictly convex polygon. Let us say that a face $F$ of the dual tesselation is \emph{concave at} $i$, if either $F$ has an angle at least $\pi$ at $i$ or $i$ is an isolated vertex lying in the closure of $F$. Consider the set
$$
I = \{i \in \Sigma\ |\ \mbox{there is a face concave at }i\}
$$
and the vector $\xi$ with components
\begin{equation}
\label{eqn:Xi}
\xi_i = \left\{
\begin{array}{rl}
-1, & \mbox{for } i \in I;\\
0, & \mbox{for } i \notin I.
\end{array}
\right.
\end{equation}
By Lemma \ref{lem:ConcFace}, we have
$$
\kappa_i < 0 \quad \mbox{for all }i \in I.
$$
Due to this and to $\sum_i \kappa_i = 0$, the set $I$ is a proper subset of $\Sigma$. Since $I$ is also non-empty, the vector $\xi$ has a non-zero projection $[\xi]$ on $\R^\Sigma/\mathbb{L}$.

We have
$$
\frac{\partial V}{\partial [\xi]} = - \sum_{i \in I} \kappa_i > 0.
$$
If we show that $[\xi] \in C_{[h]}\M^*$, then this will imply that $[h]$ is not a critical point of $-V$ in the sense of Definition \ref{dfn:CritPoint}.

In order to prove that $[\xi] \in C_{[h]}\M^*$, it suffices to show that the point $[h] + t[\xi]$ belongs to $\M^*$ for all sufficiently small positive $t$. Recall that $\M^*$ is the quotient by $\mathbb{L}$ of the solution set of the system
\begin{equation}
\label{eqn:QuadPrime}
h_i \le \widetilde{h_{jkl}}(i),
\end{equation}
one equation for every quadrilateral $ikjl$ with an angle at least $\pi$ at $i$, see Proposition \ref{prp:DescrM}. Let us see what happens to the inequalities \eqref{eqn:QuadPrime} when we add $t\xi_i$ to each of $h_i$. If an inequality \eqref{eqn:QuadPrime} is strict, then it remains strict for all sufficiently small $t$. If \eqref{eqn:QuadPrime} holds as equality, then by the argument in the first paragraph of the proof of Lemma \ref{lem:BoundaryM} there is a face of the dual tesselation with an angle at least $\pi$ at $i$. Thus \eqref{eqn:Xi} implies $\xi_i = -1$ and
\begin{equation}
\label{eqn:xixi}
\xi_j \ge \xi_i,\ \xi_k \ge \xi_i,\ \xi_l \ge \xi_i.
\end{equation}
Recall the following two properties of the function $\widetilde{h_{jkl}}(i)$ (see (\ref{eqn:QuadEx})). First, it is a monotone increasing function of $h_j$, $h_k$, $h_l$. Second, when a constant is added to each of $h_j$, $h_k$, $h_l$, then $\widetilde{h_{jkl}}(i)$ increases by the same constant. Thus \eqref{eqn:xixi} implies that the inequality \eqref{eqn:QuadPrime} remains valid when $t\xi$ is added to $h$. Thus $[\xi] \in C_{[h]}\M^*$ and the proposition is proved.
\end{proof}

\begin{lem}
\label{lem:ConcFace}
Let $M$ be a cusp with coparticles, and let $i$ be such that $\kappa_i \ge 0$. Then in the dual tesselation, all faces adjacent to $i$ have angles less than $\pi$ at~$i$.
\end{lem}
\begin{proof}
Let us study the geometry of the $i$-th face $F_i$ of the cusp $M$. Choose a triangulated truncation $(T,h)$ of $M$ (for definitions, see Section \ref{subsec:CuspFromCorn}). The truncated corners corresponding to the triangles from the star of $i$ are glued cyclically around the $i$-th coparticle. Their faces orthogonal to the coparticle are glued cyclically to form the face $F_i$.

Let $1,2,\ldots,n$, in this cyclic order, be the vertices adjacent to $i$. Simplify the notations from Section \ref{subsec:Bricks} by putting
$$
\begin{array}{rcl}
\gamma^{j,j+1} & := & \gamma_i^{j,j+1},\\
g_j & := & h_{ij},\\
g_{j,j \pm 1} & := & h_{i,j,j \pm 1},
\end{array}
$$
see Figure \ref{fig:FaceI}. The edges of the face $F_i$ have lengths
$$
\ell_{ij} = g_{j,j-1} + g_{j,j+1},
$$
and for every $j$ the inequality $\ell_{ij} \ge 0$ holds. Recall that the dual tesselation is obtained from the triangulation $T$ by erasing all edges $jk$ such that $\ell_{jk} = 0$. Thus we need to prove that if in the middle of Figure \ref{fig:FaceI} we have $\omega_i \le 2\pi$, then at the left of Figure \ref{fig:FaceI}, after erasing all edges with $\ell_{ij} = 0$, all angles at $i$ are still less than $\pi$.

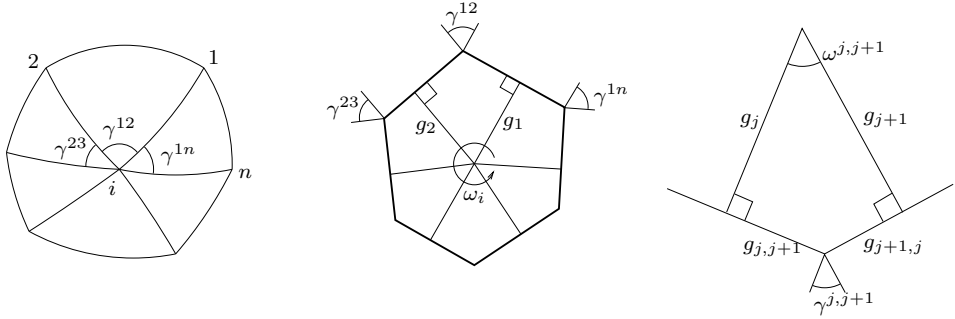
\begin{figure}[ht]
\begin{center}
\input{Fig/FaceI.tex}
\end{center}
\caption{The face $F_i$ (middle) is made up of several 2--corners (right) whose geometry is determined by the star of $i$ in the dual tesselation (left) and by support parameters $g_j$.}
\label{fig:FaceI}
\end{figure}

In the middle of Figure \ref{fig:FaceI} we depicted a very nice situation. In general, the quadrilaterals the face $F_i$ is glued from can be self-intersecting, and the cone point can lie outside the face. The following definitions provide a formal basis.

\begin{dfn}
A \emph{2--corner} in $\H^2$ consists of a point (the vertex of the corner) and two intersecting co-oriented lines (boundary lines of the corner). The angle between the boundary lines measured respecting their co-orientation and the signed distances from the vertex to the boundary lines are called the \emph{parameters} of the 2--corner.
\end{dfn}
The face $F_i$ is made up of corners with parameters $(\gamma^{j,j+1}; g_j, g_{j+1})$, see Figure~\ref{fig:FaceI}, right.

\begin{dfn}
Draw through the vertex of a 2--corner perpendiculars to its boundary lines. Orient each perpendicular in the sense opposite to that induced by the co-orientation of the corresponding boundary line. The angular region bounded by the positive halves of the perpendiculars is called the \emph{link} of the corner.

The angle measure of the link is called the \emph{central angle} of the 2--corner.
\end{dfn}
See Figure \ref{fig:CentralAngles} for examples. The central angles $\omega^{j,j+1}$ of the corners constituting the face $F_i$ sum up to the cone angle $\omega_i$.

\begin{figure}[ht]
\begin{center}
\includegraphics{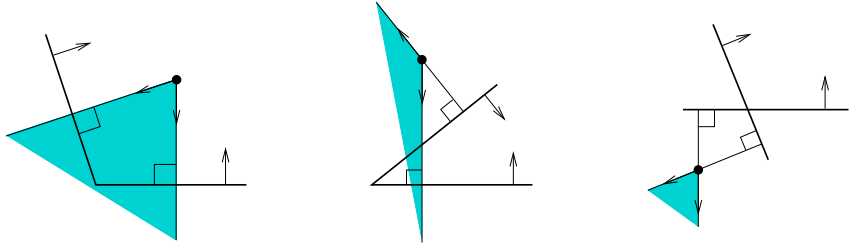}
\end{center}
\caption{Examples of 2--corners and their links.}
\label{fig:CentralAngles}
\end{figure}

Now let us proceed with the proof of Lemma \ref{lem:ConcFace}. Choose the biggest of the numbers $|g_1|$, $|g_2|, \ldots, |g_n|$. Without loss of generality, let this be $|g_1|$. If $g_j = 0$ for all $j$, then we have $\omega^{j,j+1} = \gamma^{j,j+1}$ for all $j$ and thus
\begin{equation}
\label{eqn:OmegaGamma}
\omega_i = \sum_{j=1}^n \gamma^{j,j+1}.
\end{equation}
On the other hand, we have
\begin{equation}
\label{eqn:SumGamma}
\sum_{j=1}^n \gamma^{j,j+1} = \theta_i > 2\pi
\end{equation}
because $\theta_i$ is the cone angle at the singular point $i$ of $(\Tor^2, g)$. Thus \eqref{eqn:OmegaGamma} implies $\omega_i > 2\pi$ that contradicts the assumption $\kappa_i \ge 0$ of the lemma.

If not all of $g_j$ vanish, then we have $|g_1| > 0$. By formula \eqref{eqn:h123} we have
\begin{eqnarray*}
\sinh g_{1n} & = & \frac{-\cos\gamma^{1n}\sinh g_n + \sinh g_1}{\sin\gamma^{1n} \cosh g_n},\\
\sinh g_{12} & = & \frac{-\cos\gamma^{12}\sinh g_2 + \sinh g_1}{\sin\gamma^{12} \cosh g_2}.
\end{eqnarray*}
Thus the numbers $g_{1n}$ and $g_{12}$ have the same sign as $g_1$. Since $g_{1n} + g_{12} = \ell_{i1} \ge 0$, we have
\begin{equation}
\label{eqn:PosG}
g_1 > 0,\ g_{1n} > 0,\ g_{12} > 0.
\end{equation}

Cut the face $F_i$ along the perpendicular to the first side and develop it on the plane. Since $\omega_i \le 2\pi$, the union of the links of the corners is a non-overlapping angular region (or the whole plane, if $\omega_i = 2\pi$). Due to this and to \eqref{eqn:PosG}, the boundary of $F_i$ represents a convex polygonal line, see Figure \ref{fig:CutFace}. Note that the cone point might lie inside or outside the face.

\begin{figure}[ht]
\begin{center}
\input{Fig/CutFace.tex}
\end{center}
\caption{Two examples of a face with positive curvature cut along its biggest height.}
\label{fig:CutFace}
\end{figure}
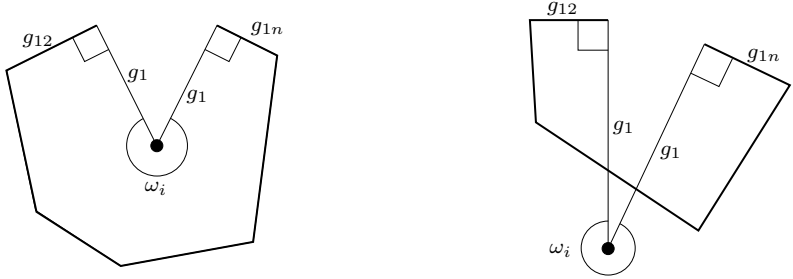

Assume that there is a face of the dual tesselation with an angle greater or equal $\pi$ at $i$. This means that several consecutive sides of the face $F_i$ have zero length, so that the exterior angle between two adjacent non-zero sides is at least $\pi$. If the angle is bigger than $\pi$, then it contradicts the convexity of the developed boundary shown in the previous paragraph. If there is an exterior angle equal to $\pi$, then the convexity of the developed boundary implies that $F_i$ degenerates into a segment so that $\omega_i = 2\pi$. Then we have
$$
\sum_{j=1}^n \gamma^{j,j+1} = \pi + \pi = 2\pi
$$
that contradicts \eqref{eqn:SumGamma}. Thus no face of the dual tesselation can have an angle greater or equal $\pi$ at $i$, and Lemma \ref{lem:ConcFace} is proved.
\end{proof}

\begin{rem}
Figure \ref{fig:ConcFace} shows an example with $\kappa_i < 0$ and an angle greater than $\pi$ at $i$ in the dual tesselation.
\begin{figure}[ht]
\begin{center}
\includegraphics{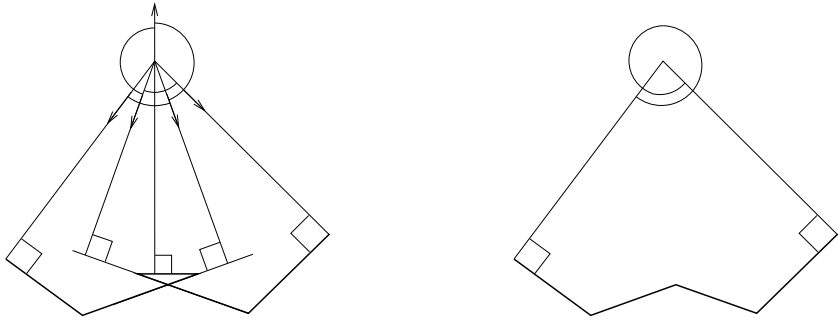}
\end{center}
\caption{The development at the right is the limit of developments like one at the left, as the shortest side shrinks to a point.}
\label{fig:ConcFace}
\end{figure}
\end{rem}

\begin{rem}
Lemma \ref{lem:ConcFace} can be viewed as an analog of Volkov's lemma from \cite{Vol60}, see also Lemma 5.3 in \cite{Izm08}.
\end{rem}

\subsection{Behavior of $V$ at infinity}
\label{sec:VInf}
\begin{prp}
\label{prp:VInf}
For every $h \in \R^\Sigma$ put
$$
d(h) = \max_{i,j \in \Sigma} |h_i - h_j|.
$$
Then there exists $D > 0$ such that for every $h \in \Mt^*$ with $d(h) \ge D$ there is a tangent vector $\xi \in C_h(\Mt^* \cap Q_{d(h)})$ such that
$$
\frac{\partial V}{\partial \xi} > 0.
$$
\end{prp}
\begin{proof}
Denote
\begin{eqnarray*}
h_+ & = & \max_{i \in \Sigma} h_i;\\
h_- & = & \min_{i \in \Sigma} h_i.
\end{eqnarray*}
In the case when $h \in \partial \Mt^*$, take $\xi$ as in \eqref{eqn:Xi}. We have already shown in the proof of Proposition \ref{prp:CritPtBdry} that $\frac{\partial V}{\partial \xi} > 0$ and that $\xi \in C_h \Mt^*$. In order to show that $\xi \in C_h Q_{d(h)}$ it suffices to prove that $h_i > h_-$ holds for all $i$ such that there is a face concave at $i$. This is indeed true because the support function of a cusp has the form \eqref{eqn:SuppLike} on every face, thus is a concave function and cannot attain its minimum at $i$.

From now on let us assume that $h$ lies in the interior of $\Mt^*$ and that for all $\xi \in C_h Q_{d(h)}$ we have $\frac{\partial V}{\partial \xi} \le 0$. The latter means that
\begin{equation}
\label{eqn:KappaCrit}
\kappa_i \left\{
\begin{array}{rl}
\le 0 & \mbox{ if }h_i = h_-;\\
= 0 & \mbox{ if }h_- < h_i < h_+;\\
\ge 0 & \mbox{ if }h_i = h_+
\end{array}
\right.
\end{equation}
holds for all $i \in \Sigma$. Our goal is to show that this is impossible if
$$
h_+ - h_- \ge D,
$$
for an appropriately chosen $D$ that depends on the geometry of $(\Tor^2,g)$.

Order the heights:
$$
h_- = h_{i_1} \le h_{i_2} \le \cdots \le h_{i_{|\Sigma|}} = h_+.
$$
Choose a big gap between two neighbors:
$$
h_{i_{s+1}} - h_{i_s} \ge \frac{D}{|\Sigma|-1},
$$
and split $\Sigma$ into two subsets:
$$
\Sigma_- = \{h_{i_1},\ldots, h_{i_s}\}, \quad \Sigma_+ = \{h_{i_{s+1}},\ldots,h_{i_{|\Sigma|}}\}.
$$
Then we have
\begin{equation}
\label{eqn:BigGap}
h_i - h_j \ge \frac{D}{|\Sigma|-1}, \quad \mbox{for all }i \in \Sigma_+, j \in \Sigma_-.
\end{equation}

Let $T$ be a triangulation of $(\Tor^2,g)$ associated with the height vector $h$, and let $G$ be the $1$--skeleton of $T$. Let $G_+$ and $G_-$ be the induced subgraphs of $G$ on the vertex sets $\Sigma_+$ and $\Sigma_-$, respectively. Consider the connected components of $G_+$ and $G_-$:
$$
G_+ = G_+^1 \sqcup \cdots \sqcup G_+^p, \quad G_- = G_-^1 \sqcup \cdots \sqcup G_-^q
$$
and their vertex sets
$$
\Sigma_+ = \Sigma_+^1 \sqcup \cdots \sqcup \Sigma_+^p, \quad \Sigma_- = \Sigma_-^1 \sqcup \cdots \sqcup \Sigma_-^q.
$$
The cell decomposition of the torus $\Tor^2$ dual to the triangulation $T$ associates with each point $i \in \Sigma$ a $2$--cell $N(i)$. For a subset $X$ of $\Sigma$, denote
$$
N(X) = \bigcup_{i \in X} N(i).
$$
Then $N(\Sigma_+^1), \ldots, N(\Sigma_+^p)$ and $N(\Sigma_-^1), \ldots, N(\Sigma_-^q)$ are the connected components of $N(\Sigma_+)$ and $N(\Sigma_-)$, respectively. Note that each of $N(\Sigma_\pm^s)$ is a compact surface with boundary. An Euler characteristic argument shows that either one of  $N(\Sigma_\pm^s)$ is homeomorphic to the disk, or all of them are homeomorphic to the annulus. The subsequent argument deals with these three cases separately.

\begin{rem}
\label{rem:Uniform}
The lengths of closed contractible geodesics on $(\Tor^2,g)$ form a discrete subset of $\R$. For closed contractible geodesics that contain at least one cone point, this follows from Lemma \ref{lem:Troyanov}. To those that contain no cone points and whose length is not a multiple of $2\pi$, the argument in the proof of \cite[Proposition~1]{ILTC} can be applied.

Since, by assumption of Theorem \ref{thm:Main}, all closed contractible geodesics have lengths bigger than $2\pi$, there exists a number $\delta > 0$ such that all closed contractible geodesics have lengths at least $2\pi + \delta$.
\end{rem}

\noindent\textbf{Case 1.} \emph{There is a component $N(\Sigma^s_+)$ homeomorphic to the disk.}
\par\noindent In this case we will show that, if $D$ is chosen sufficiently large, there is a contractible closed geodesic in $(\Tor^2,g)$ of length less than $2\pi + \delta$. This contradicts Remark \ref{rem:Uniform}.

Let $T'$ be a subcomplex of $T$ consisting of all triangles that have a non-empty intersection with $N(\Sigma^s_+)$.
There is a map
\begin{equation}
\label{eqn:DiskToTor}
\Disk^2 \to \Tor^2
\end{equation}
whose image is the union of all triangles of $T'$ and which is injective when restricted to the interior of the disk $\Disk^2$. By pulling back through \eqref{eqn:DiskToTor}, we view $T'$ as a triangulation of the disk $\Disk^2$ and draw it on the plane, see Figure~\ref{fig:Case1}.

\begin{figure}[ht]
\begin{center}
\includegraphics{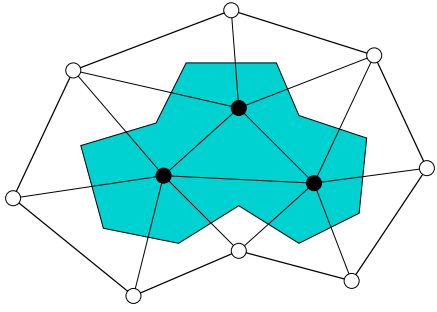}
\end{center}
\caption{The triangulation $T'$ of disk $\Disk^2$. Vertices from $\Sigma^s_+$ are colored black, vertices from $\Sigma_-$ are colored white, the set $N(\Sigma^s_+)$ is shaded.}
\label{fig:Case1}
\end{figure}

By the boundary of $T'$ we mean the polygonal curve that bounds $\Disk^2$. When mapped to $\Tor^2$ via \eqref{eqn:DiskToTor}, this curve can run twice through some edges and several times through some vertices of the triangulation $T$.

\setcounter{clm}{0}
\begin{clm}
\label{clm:BdryT}
If $D$ is sufficiently large, then the boundary of $T'$ has length less than $2\pi + \delta$.
\end{clm}

Let $V$, $E$, and $F$ denote the numbers of vertices, edges, and triangles of the triangulation $T'$, respectively. We have
$$
\begin{array}{rcl}
V & = & V_\iota + V_\partial;\\
E & = & E_\iota + E_\partial;\\
F & = & F_1 + F_2 + F_3.
\end{array}
$$
Here $V_\iota$ and $E_\iota$ denote the numbers of vertices and edges in the interior of~$T'$; $V_\partial$ and $E_\partial$ denote the numbers of vertices and edges on the boundary of~$T'$; $F_1$, $F_2$, and $F_3$ denote the numbers of triangles with $1$, $2$, and $3$ vertices in the set $\Sigma^s_+$, respectively. We have
$$
\begin{array}{c}
V - E + F = 1;\\
V_\partial = E_\partial = F_1;\\
3F = 2E_\iota + E_\partial,
\end{array}
$$
which easily implies
\begin{equation}
\label{eqn:VF}
V_\iota = \frac{1}{2}(F_2 + F_3) + 1.
\end{equation}

By the assumption \eqref{eqn:KappaCrit},
\begin{equation}
\label{eqn:SumOmega1}
\sum_{i \in \Sigma^s_+} \omega_i \le 2\pi V_\iota.
\end{equation}
On the other hand,
\begin{equation}
\label{eqn:SplitSum}
\sum_{i \in \Sigma^s_+} \omega_i = \frac12 \sum_{i,j,k \in \Sigma^s_+} \omega_i^{jk} + \sum_{\substack{i,j \in \Sigma^s_+\\ k \in \Sigma_-}} \omega_i^{jk} + \frac12\sum_{\substack{i \in \Sigma^s_+\\ j,k \in \Sigma_-}} \omega_i^{jk}.
\end{equation}
(The factors $\frac12$ are needed because of $\omega_i^{jk} = \omega_i^{kj}$, so that in the first and the third sum on the right hand side all angles appear twice.)

The first sum on the right hand side of \eqref{eqn:SplitSum} splits into $F_3$ groups of the form $\omega_i^{jk} + \omega_j^{ik} + \omega_k^{ij} = \pi$, the second and the third sum are estimated by Lemma \ref{lem:BigDiff}. As a result, we have 
\begin{eqnarray*}
\sum_{i \in \Sigma^s_+} \omega_i & > & \pi F_3 + (\pi - \epsilon)F_2 + \sum_{jk \in \partial T'} (\alpha_{jk} - \epsilon)\\
& > & \pi (F_2 + F_3) - 6|\Sigma|\epsilon + \sum_{jk \in \partial T'} \alpha_{jk},
\end{eqnarray*}
where $\epsilon$ can be made arbitrarily small by choosing large $D$. (Note that the triangulation $T'$ can depend on $D$, and the speeds of the convergencies in Lemma \ref{lem:BigDiff} depend on the shapes of the triangles. But, since the number of triangles on $(\Tor^2, g)$ is finite, $\epsilon$ can be viewed as depending only on $D$.)

By combining the last inequality with \eqref{eqn:VF} and \eqref{eqn:SumOmega1}, we obtain
$$
\sum_{jk \in \partial T'} \alpha_{jk} < 2\pi + 6|\Sigma|\epsilon,
$$
which proves Claim \ref{clm:BdryT}.

Through the map \eqref{eqn:DiskToTor}, the metric $g$ on $\Tor^2$ induces on the disk $\Disk^2$ a spherical cone metric which, by abuse of notation, is also denoted by $g$. As next, we will show that $(\Disk^2, g)$ has a locally concave boundary.

\begin{clm}
\label{clm:BigAngle}
If $D$ is sufficiently large, then all angles spanned in $(\Disk^2, g)$ by two consecutive boundary edges of the triangulation $T'$ are at least $\pi$.
\end{clm}

Assume the converse: there is a vertex $i$ on the boundary of $T'$ such that the angle in $(\Disk^2, g)$ formed by the boundary edges $ij$, $ik$ of $T'$ is less than~$\pi$. Develop on the sphere the union of all triangles of $T'$ adjacent to~$i$. We obtain a spherical polygon $P$, see Figure \ref{fig:BdryAngle}, left. The polygon $P$ is triangulated, so that all triangles have a common vertex $i$. Aside from $i$, $j$, and $k$, polygon~$P$ has at least one other vertex $l$, otherwise $i$ has no neighbor in $\Sigma^s_+$. Draw the shortest geodesic arc $jk$ on the sphere. Note that $jk$ is not necessarily contained in $P$; thus there may be no corresponding arc in $(\Tor^2, g)$. Let $l'$ be the intersection point of the arc $jk$ with the arc $il$ or with its extension.

\begin{figure}[ht]
\begin{center}
\input{Fig/BdryAngle.tex}
\end{center}
\caption{}
\label{fig:BdryAngle}
\end{figure}

The exponent of the support function $\widetilde{h_T}$ is the restriction of a convex function on $\cone(P) \subset \cone(\Sph^2) = \R^3$, linear on the cone over each triangle of the triangulation of $P$, see Section \ref{sec:SpaceOfCusps}. The values of $\exp \widetilde{h_T}$ at the vertices $i$, $j$, $k$, $l$, are $e^{h_i}$, $e^{h_j}$, $e^{h_k}$, $e^{h_l}$, respectively. If $jk$ is not contained in $P$, we extend the function $\widetilde{h_T}$ in a support-like way along the rays starting from $i$. Then the restriction of $\exp \widetilde{h_T}$ to $\cone(jk)$ is a convex PL function, and the restriction of $\exp \widetilde{h_T}$ to $\cone(il)$ is a linear function. Thus we have
\begin{eqnarray}
\widetilde{h_T}(l') & \le & \log(ae^{h_j} + be^{h_k}), \label{eqn:jkl'}\\
\widetilde{h_T}(l') & = & \log(ce^{h_i} + de^{h_l}). \label{eqn:ill'}
\end{eqnarray}
Here the numbers $a$, $b$, $c$, $d$ depend only on the relative position of the points $i$, $j$, $k$, $l$ on $\Sph^2$. Besides, $a$, $b$, and $d$ are positive, whereas $c$ is positive or negative depending on whether $jk$ intersects $il$ or its extension beyond $l$. The inequalities \eqref{eqn:jkl'} and \eqref{eqn:ill'} imply
\begin{equation}
\label{eqn:jkil}
ae^{h_j} + be^{h_k} \ge ce^{h_i} + de^{h_l}.
\end{equation}

The inequality \eqref{eqn:jkil} can be rewritten as
\begin{equation}
\label{eqn:dLess}
d \le ae^{h_j-h_l} + be^{h_k-h_l} - ce^{h_i-h_l}.
\end{equation}
Since $i,j,k \in \Sigma_-$ and $l \in \Sigma_+$, it follows from \eqref{eqn:BigGap} that each of the differences $h_j-h_l$, $h_k-h_l$, $h_i-h_l$ is smaller or equal $-\frac{D}{|\Sigma|-1}$. By choosing $D$ large, the right hand side of \eqref{eqn:dLess} can be made arbitrarily close to $0$, whereas the left hand side is a positive constant by the previous paragraph. This contradiction shows that the angle at $i$ in $(\Disk^2, g)$ cannot be less than $\pi$.

Again, while $D$ is changing, the triangulation $T'$ can vary, as well as a suspected vertex $i$ with an angle less than $\pi$. But, since the number of configurations is finite, the constants $a$, $b$, $c$, $d$ take only finitely many values. Thus, a sufficiently large value of $D$ in the above argument can be chosen uniformly. Claim 2 is proved.

\begin{rem}
Let us explain the geometric meaning of Claim 2. Consider the associated cusp with coparticles. Fit together the corners that correspond to the triangles of $T'$ adjacent to $i$. We obtain a configuration of hyperbolic planes with combinatorics outlined on Figure \ref{fig:BdryAngle}, right. The unmarked faces lie further from the truncating horosphere than the faces marked with $i$, $j$, and $k$, by the distance at least $\frac{D}{|\Sigma|-1}$. If, by making $D$ arbitrarily large and preserving the Gauss images of all of the vertices, we can still achieve that all of the edges have non-negative lengths, then it is intuitively clear that the edges shared by the face $i$ with the faces $j$ and $k$ don't intersect, when extended downwards. It follows that the sum of the exterior angles of the face $i$ at the depicted vertices is at least~$\pi$. 
\end{rem}

To finish dealing with Case 1, consider the space
$$
\Fl =  \Tor^2 \setminus \inn \Disk^2,
$$
with the metric $g$, where $\inn \Disk^2$ is the image of the interior of $\Disk^2$ under the map \eqref{eqn:DiskToTor}. The image of the boundary of $\Disk^2$ is a closed curve $\Gamma$ in $\Fl$, non-contractible in $\Fl$. By Claim \ref{clm:BdryT}, $\Gamma$ has length less than $2\pi + \delta$. By Lemma \ref{lem:CurveShort}, $\Fl$ contains a closed geodesic $\Gamma'$ homotopic to $\Gamma$, thus contractible in $\Tor^2$, and no longer than $\Gamma$. Claim~\ref{clm:BigAngle} implies that $\Gamma'$ is also a geodesic in $(\Tor^2,g)$. This leads to a contradiction with Remark \ref{rem:Uniform}.

\begin{rem}
Case 1 is similar to the situation of ``long thin tube'' appearing in closeness results in \cite{RivHod93}.
\end{rem}
\bigskip

\noindent\textbf{Case 2.} \emph{There is a component $N(\Sigma^s_-)$ homeomorphic to the disk.}
\par\noindent Define a complex $T'$ similarly to the Case 1. Now on Figure \ref{fig:Case1} the black points form the set $\Sigma^s_-$, the white points are their neighbors from $\Sigma_+$. Instead of \eqref{eqn:SumOmega1}, we have
\begin{equation}
\label{eqn:SumOmega2}
\sum_{i \in \Sigma^s_-} \omega_i \ge 2\pi V_\iota.
\end{equation}
On the other hand, by estimating the sums on the right hand side of \eqref{eqn:SplitSum}, where plus and minus signs should be interchanged, with the help of Lemma \ref{lem:BigDiff} we obtain
\begin{eqnarray}
\sum_{i \in \Sigma^s_-} \omega_i & < & \pi F_3 + \frac12 \sum_{\substack{i,j \in \Sigma^s_-\\ k \in \Sigma_+}} (\pi - \alpha_{ij} + \epsilon) + F_1 \epsilon \label{eqn:LongIneq}\\
& < & \pi(F_2 + F_3) - \frac12 \sum_{\substack{i,j \in \Sigma^s_-\\ k \in \Sigma_+}} \alpha_{ij} + 6|\Sigma|\epsilon, \nonumber
\end{eqnarray}
with $\epsilon \to 0$ as $D \to \infty$. With \eqref{eqn:VF} taken into account, this contradicts \eqref{eqn:SumOmega2}.

\bigskip

\noindent\textbf{Case 3.} \emph{All components of $N(\Sigma_+)$ and $N(\Sigma_-)$ are homeomorphic to the annulus.}
\par\noindent Let $N(\Sigma^s_-)$ be a component of $N(\Sigma_-)$. The complex $T'$, with identifications on the boundary resolved, is homeomorphic to the annulus. Therefore we have
$$
V_\iota = \frac{1}{2}(F_2 + F_3)
$$
instead of \eqref{eqn:VF}. This still suffices to arrive at a contradiction between \eqref{eqn:SumOmega2} and \eqref{eqn:LongIneq}.

\bigskip

Thus, all three cases are resolved, which means that the inequalities \eqref{eqn:KappaCrit} cannot hold if $h_+ - h_-$ is large enough, which in turn means that a tangent vector $\xi$ as required in the proposition exists. Proposition \ref{prp:VInf} is proved.
\end{proof}

\begin{lem}[Degeneration of corners]
\label{lem:BigDiff}
Let $\Delta$ be a spherical triangle with edge lengths $\alpha_{12}$, $\alpha_{23}$, $\alpha_{31}$. Consider a truncated corner with Gauss image $\Delta$ and heights $h_1$, $h_2$, $h_3$. Denote by $\omega_1^{23}$, $\omega_2^{31}$, $\omega_3^{12}$ the dihedral angles at the heights of the corresponding brick, see Figure \ref{fig:Brick}. Then we have
\begin{enumerate}
\item if $h_1 - h_2 \to +\infty$ and $h_1 - h_3 \to +\infty$, then $\omega_1^{23} \to \alpha_{23}$;
\item if $h_1 - h_3 \to +\infty$ and $h_2 - h_3 \to +\infty$, then $\omega_3^{12} \to 0$.
\end{enumerate}
\end{lem}
\begin{proof}
By \eqref{eqn:h12}, if $h_1 - h_2$ and $h_1 - h_3$ tend to $+\infty$, then $\sinh h_{12}$ and $\sinh h_{13}$ tend to $-\cot \alpha_{12}$ and $-\cot \alpha_{13}$, respectively. By substituting this into \eqref{eqn:OmegaFromH} and using the spherical Cosine Law, we obtain the first statement of the lemma. The second statement is proved in a similar way.
\end{proof}

Let $T$ be a geodesic triangulation of a spherical cone-surface $\Fl$. Let $T'$ be a subcomplex of $T$, and let $\Fl'$ be the union of all simplices from $T'$. Note that the complex $T'$ is not required to be pure, that is $T'$ can contain edges that are not adjacent to any triangle of $T'$. This means that $\Fl'$ is not necessarily a surface with boundary.

\begin{dfn}
A curve $\Gamma$ in $\Fl'$ is called a \emph{geodesic}, if $\Gamma$ is locally length minimizing.
\end{dfn}

Outside the vertices of $T'$, a geodesic looks locally like a great circle. Thus a geodesic either does not pass through any vertex or consists of a sequence of spherical geodesic arcs joining vertices. At a vertex of $T'$, a geodesic must span angles in $\Fl'$ of at least $\pi$ on both sides. Note that the \emph{angle in $\Fl'$ at a vertex} $i \in \Gamma$ on a particular side of $\Gamma$ is defined only if at that side of $\Gamma$ the vertex $i$ has a neighborhood homeomorphic to the half-plane. By excluding angles less than $\pi$ we also exclude zero angles, that is a geodesic is not recurring.

By a \emph{polygonal curve} on $\Fl'$ we mean a curve that consists of a sequence of spherical geodesic arcs joining cone points. Note that in this case we allow a spherical geodesic arc to have length $\pi$ or bigger.

\begin{lem}[Curve shortening]
\label{lem:CurveShort}
Let $\Gamma$ be an arbitrary homotopically non-trivial closed polygonal curve in $\Fl'$. Then $\Gamma$ is homotopic to a closed geodesic~$\Gamma'$ in $\Fl'$ such that $\Gamma'$ is no longer than $\Gamma$.
\end{lem}
\begin{proof}
We will show that if $\Gamma$ is not a geodesic, then it can be homotoped to a polygonal curve of smaller length. Lemma \ref{lem:Troyanov} implies that the lengths of closed polygonal curves form a discrete subset of $\R$. Because of this and since $\Gamma$ is non-contractible, by successive shortening we arrive sooner or later to a geodesic homotopic to $\Gamma$.

So assume that $\Gamma$ spans in $\Fl'$ an angle less than $\pi$ at a vertex $i$. If $\Gamma$ is recurring at $i$, then it contains two copies of the same arc run in opposite directions, and we can shorten $\Gamma$ by removing them. If the angle at $i$ is bigger than $0$, then consider a point $z$ that moves along the bisector starting from $i$, see Figure \ref{fig:CurveShort1}.

Let $ji$, $ki$ be the arcs of $\Gamma$ ending at $i$. Close to them there are geodesic arcs $jz$, $kz$ in $\Fl$. Let us show that for all $z$ sufficiently close to $i$, the arcs $jz$ and $kz$ lie in $\Fl'$. Assume that $ji$ is an edge of $T$. Then $\Fl'$ contains a neighborhood $U$ of $i$ on the side of $z$ from $\Gamma$ and the triangle $\Delta$ of $T$ adjacent to $ji$ on that side. Since the length of $ji$ is less than $\pi$, the arc $jz$ is contained in $U \cup \Delta$ for $z$ sufficiently close to $i$. If $ji$ is not an edge of $T$, then a whole two-sided neighborhood $W$ of the interior of $ji$ is in $\Fl'$, and the arc $jz$ is contained in $U \cup W$. See Figure \ref{fig:CurveShort1}.

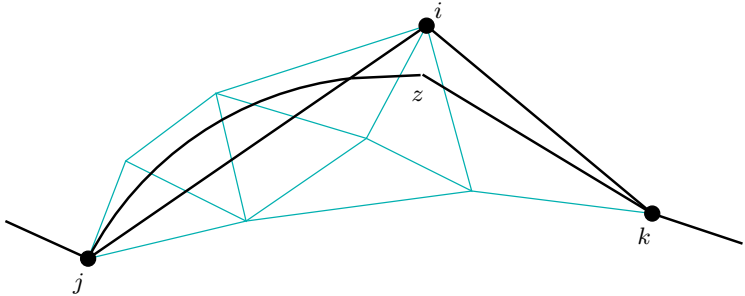
\begin{figure}[ht]
\begin{center}
\input{Fig/CurveShort1.tex}
\end{center}
\caption{Shortening a polygonal curve; here the length of $ji$ is bigger than $\pi$.}
\label{fig:CurveShort1}
\end{figure}

It is not hard to show that during the deformation the sum of the lengths of $jz$ and $kz$ decreases. We stop moving $z$ if one of the following three events occur. First, if the angle at $z$ becomes $\pi$; then we have homotoped $\Gamma$ into a shorter polygonal curve. Second, if $z$ arrives at a vertex of $T'$; then we also have homotoped~$\Gamma$ into a shorter polygonal curve. Third, if one of the geodesic arcs $jz$ or $kz$ (say, this is $jz$) meets a vertex $j_1$. In this last case we continue pushing the point $z$ inside the angle which is smaller than $\pi$, now deforming the arcs $j_1z$ and $kz$, and wait for the next event to occur. Note that if $z$ meets the boundary of $\Fl'$, then either $z$ arrives at a vertex or the angle at $z$ becomes $\pi$.

If $i$ was the only vertex on $\Gamma$, then at the beginning we are deforming a geodesic loop based at $z$, see Figure \ref{fig:CurveShort2}, left. In this case, if the angle at $z$ becomes $\pi$, we obtain a geodesic avoiding the vertices. The other two kinds of events from the previous paragraph yield the same results as described there. On Figure \ref{fig:CurveShort2}, two encounters with new vertices are depicted.

We cannot stumble upon new vertices permanently, because the length of our curve keeps decreasing, and if a closed curve on~$\Fl$ runs through $n$ vertices, then its length is at least $n$ times the length of a shortest arc between two vertices. Thus we will end up with a polygonal curve shorter than $\Gamma$.

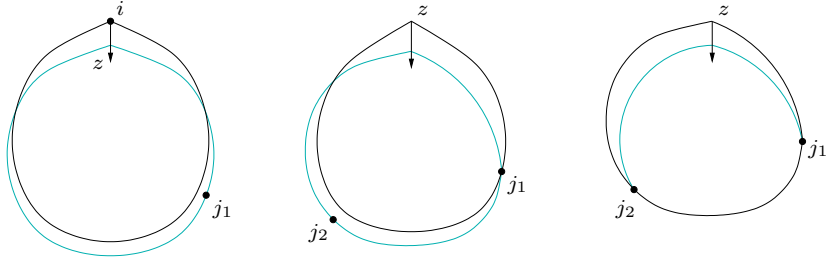
\begin{figure}[ht]
\begin{center}
\input{Fig/CurveShort2.tex}
\end{center}
\caption{Coming across new vertices during the deformation of $\Gamma$.}
\label{fig:CurveShort2}
\end{figure}

The lemma is proved.
\end{proof}

\section{Concluding remarks}

\subsection{Andreev's theorem for cusps}
\label{subsec:AndThm}
As mentioned in the Example \ref{exl:UniqueTriang}, if all edge lengths of a triangulation $T$ of $(\Tor^2,g)$ are at least $\frac{\pi}2$, then $T$ is a unique geodesic triangulation of $(\Tor^2,g)$. This allows to characterize convex polyhedral cusps with non-obtuse dihedral angles in terms of their combinatorics and values of dihedral angles. The following theorem deals only with the case of acute dihedral angles.

\begin{Alphatheorem}
\label{thm:AndreevCusps}
Let $C$ be a cellular subdivision of the torus $\Tor^2$ with trivalent vertices. Call three pairwise disjoint edges $e$, $f$, $g$ of $C$ a \emph{proper cutset}, if there exists a simple closed curve that intersects each of $e$, $f$, $g$ exactly once and bounds a disk on $\Tor^2$.

Let $\phi \co e \to \phi_e$ be a map from the edge set of $C$ to the interval $(0, \frac{\pi}2)$ such that the following conditions hold.
\begin{enumerate}
\item \label{it:1} If $e$, $f$, $g$ are three edges incident to a vertex, then $\phi_e + \phi_f + \phi_g > \pi$.
\item \label{it:2} If $e$, $f$, $g$ form a proper cutset, then $\phi_e + \phi_f + \phi_g < \pi$.
\end{enumerate}
Then there exists a unique convex polyhedral cusp with face structure $C$ and dihedral angles $\phi$.
\end{Alphatheorem}
\begin{proof}
Let $T$ be a triangulation of $\Tor^2$ dual to $C$. Denote $\alpha_e = \pi - \phi_e$. Then the condition \eqref{it:1} implies that $\alpha_e$, $\alpha_f$, $\alpha_g$ satisfy the triangle inequality. Consider the spherical cone metric $g$ on $\Tor^2$ given by gluing spherical triangles with edge lengths $\alpha$. Then, as shown in \cite{Hod92}, $\alpha_e > \frac{\pi}2$ implies that all cone angles in $g$ are bigger than $2\pi$, and condition \eqref{it:2} implies that in $(\Tor^2,g)$ there is no closed contractible geodesic of length less or equal $2\pi$. Thus by Theorem \ref{thm:Main} there exists a convex polyhedral cusp $M$ with Gauss image $(\Tor^2,g)$. If we show that $T$ is the dual tesselation of $M$, then Theorem \ref{thm:AndreevCusps} follows.

As shown in \cite[Proposition 2.4]{Hod92}, every geodesic arc of length at most $\pi$ between two cone points in  $(\Tor^2,g)$ is an edge of $T$. Thus $T$ is a refinement of the dual tesselation of $M$. Since the faces of the dual tesselation are convex spherical polygons, and the union of several triangles of $T$ is not a convex polygon, $T$ is exactly the dual tesselation of $M$. The theorem is proved.
\end{proof}

In the situation of Theorem \ref{thm:AndreevCusps}, the proof of Theorem \ref{thm:Main} is substantially simplified. Since $T$ is the unique triangulation of $(\Tor^2,g)$, all cusps with coparticles have the form $(T,h)$. Further, it is easy to show that for all $h \in \R^\Sigma$, the pair $(T,h)$ yields a convex truncated cusp. Thus, $\Mt = \R^\Sigma$, and the whole Section \ref{sec:SpaceOfCusps} is redundant. In Section \ref{sec:Functional}, we don't need to take care of different triangulations, so the differentiability of $V$ is straightforward without Whitney's extension theorem. Finally, in Section \ref{sec:Proof} we don't need Morse theory on manifold with corners and don't need to study the behavior of $V$ on the boundary. We only need the argument from Section \ref{sec:VInf} showing that $V$ attains its maximum. The assumptions on the edge lengths of the triangulation $T$ again simplify that argument at several points.

\subsection{Circle patterns on the torus}
\label{subsec:Thurston}
Theorem \ref{thm:AndreevCusps} can be reformulated in terms of circle patterns on the torus. In this form, it becomes a special case of Thurston's theorem, \cite[Theorem 13.7.1]{Thurcour1}.
\begin{thm}[Thurston]
\label{thm:thurston}
Let $C$ and $\phi$ be as in the statement of Theorem \ref{thm:AndreevCusps}. Then there is a flat metric on $\Tor^2$, unique up to scalar multiple, a unique geodesic triangulation isotopic to the dual $T$ of $C$ and a unique family of circles $K_i$, centered at the vertices of $T$, such that $K_i$ and $K_j$ intersect at an angle $\phi_{ij}$ (and are disjoint if there is no edge $ij$ in $T$, and intersect at several pairs of points if there are several edges that join $i$ with $j$).
\end{thm}

Consider a convex polyhedral cusp $M$ and its universal cover $\widetilde{M} \subset \H^3$. The plane spanned by a face of $\widetilde{M}$ intersects $\partial\overline{\H^3}$ in a circle. This gives a doubly-periodic circle pattern on a Euclidean plane $\R^2 = \partial \overline{\H^3} \setminus \{o\}$, where $o$ is the fixed point of the group of covering transformations. By taking a quotient, we obtain a circle pattern on the torus. If all dihedral angles of $M$ are acute, then the planes of non-adjacent faces don't intersect, so that the nerve of the obtained circle pattern is dual to the 1--skeleton of the cusp $M$, and the intersection angles of circles are equal to the dihedral angles of $M$. Thus Theorem \ref{thm:thurston} is equivalent to Theorem \ref{thm:AndreevCusps}.

In the original version of Thurston's theorem, \cite[Theorem 13.7.1]{Thurcour1}, condition (1) of Theorem \ref{thm:AndreevCusps} is not imposed. This means that the cusp that corresponds to the circle pattern can have ideal and hyperideal vertices. Thurston also considers circle patterns on surfaces of higher genus.

Here is an outline of Thurston's proof of Theorem \ref{thm:thurston}. Let $r_i$ be the (variable) radius of the circle $K_i$. For every triangle $ijk$ of the triangulation $T$, there is a unique arrangement of circles of radii $r_i$, $r_j$, $r_k$ so that their intersection angles are $\phi_{ij}$, $\phi_{jk}$, $\phi_{ik}$, see \cite[Lemma 13.7.2]{Thurcour1}. One can fail to simultaneously arrange all circles that intersect $K_i$. By forcing the cycle of circles around $K_i$ to close, one obtains a cone singularity at the center of $K_i$. Thus, every set of radii $\{r_i\}$ gives rise to a circle pattern on a torus with cone singularities of curvatures $\{\kappa_i\}$, and the problem is to show that there exist such radii that all of $\kappa_i$ vanish. For this, Thurston investigates the behavior of curvatures when some of the radii tend to zero.

In many aspects our approach parallels that of Thurston. First, cusps with coparticles are generalizations of circle patterns on a torus with cone singularities. An arrangement of three circles corresponds to a corner, so that our Lemma \ref{lem:TrCorn} is equivalent to Thurston's lemma 13.7.2. Note also that the radii of circles and the heights of a cusp are related through
\begin{equation}
\label{eqn:rh}
r_i = e^{-h_i},
\end{equation}
up to a constant factor. Second, Thurston's Lemma 13.7.3 essentially finds the signs of the partial derivatives of our function $V$, see Sections \ref{subsec:DefV} and \ref{subsec:LocRig}, which is crucial to prove the concavity of $V$. Thurston does not consider function $V$, but the map $\{r_i\} \mapsto \{\kappa_i\}$ that he studies is related to the gradient of $V$. Third, our considerations in Section \ref{sec:VInf} deal with the situation when some heights become small with respect to others, which is similar to Thurston's study of radii tending to zero.

Our proof of Theorem \ref{thm:Main} is more involved, because a much more general situation is considered: in terms of circle patterns, not only the intersection angles can be obtuse, but also the combinatorics is not known (but a choice of combinatorics fixes intersection angles). Theorem \ref{thm:Main} generalizes Theorem \ref{thm:thurston} in the same direction in that Rivin-Hodgson's Theorem \ref{thm:RivHod} generalizes Andreev's Theorem \ref{thm:And}.

Chow and Luo gave in \cite{CL03} a proof of Thurston's theorem based on a ``combinatorial Ricci flow''. In the genus one case, they consider the evolution of radii governed by
\begin{equation}
\label{eqn:CLFlow}
\frac{dr_i}{dt} = -\kappa_i r_i
\end{equation}
and show that the limit for $t \to \infty$ gives a circle pattern with zero curvatures. They also noted that the flow \eqref{eqn:CLFlow} becomes a gradient flow after a variable change \eqref{eqn:rh}, but have not provided a geometric interpretation for the potential function. Our Proposition \ref{prp:FuncV} shows that this function is $V$. Besides, results of the present paper show that $V$ is concave independently of whether the intersection angles are acute or obtuse, thus partially answering Question 1 in \cite[Section 7]{CL03}.

In the situation of Theorem \ref{thm:Main}, one can also consider the gradient flow of function $V$. Note that the combinatorics changes during the evolution. As $V$ is concave, the flow will converge to the unique critical point of $V$ and thus will yield a cusp with a given Gauss image, provided that the solution exists for all times. In order to show that the solution exists for all times with any initial point, one has to show that at the boundary points of the domain $\M^*$ the gradient of $V$ points inwards. It would be interesting to see if this is indeed the case (in Section \ref{subsec:BehaveBoundary} we show only that the gradient is not pointing outwards orthogonally to the boundary). On the other hand, one can be satisfied with showing that the flow exists for some particular initial point. A natural candidate for an initial point would be a cusp with equal heights.

For circle patterns on surfaces of higher genus, a corresponding 3--dimensional object would be a convex polyhedral Fuchsian manifold. Cone singularities at the centers of circles would correspond to coparticles in a Fuchsian manifold. It is conceivable that the (appropriately modified) function $V$ is again concave and attains maximum. This would yield a variational proof of the higher genus case of Thurston's theorem, and more generally that of an analog of Theorem \ref{thm:Main}. Note that the higher genus analog of Theorem \ref{thm:Main} is proved in \cite{Schpoly,Fillastre2}, see also the discussion in Section \ref{subsec:RelRes}.

\subsection{De Sitter point of view on the main theorem}
\label{subsec:deSitter}
Consider the Minkowski space $\R^{3,1}$ with the norm $\|p\|^2 = - p_0^2 + p_1^2 + p_2^2 + p_3^2$. The hyperboloid model of the hyperbolic space is well-known:
$$
\H^3 = \{p \in \R^{3,1} \,|\, \|p\|^2 = -1, p_0 > 0\}.
$$
The \emph{de Sitter space} is the one-sheeted hyperboloid
$$
\dS^3 = \{p \in \R^{3,1} \,|\, \|p\|^2 = 1\}
$$
with the metric induced from $\R^{3,1}$. Thus the de Sitter space is a \emph{Lorentzian manifold}: the Minkowski scalar product restricts as an inner product of signature $(-, +, +)$ to $T_p \dS^3$ for every $p \in \dS^3$. Note that $\dS^3$ is homeomorphic to $\Sph^2 \times \R$.

By duality with respect to the Minkowski scalar product, a point $w \in \dS^3$ corresponds to a hyperbolic half-space:
\begin{equation}
\label{eqn:wDual}
w^* = \{x \in \H^3 \,|\, \langle w, x \rangle \le 0\}.
\end{equation}
Similarly, a point $v \in \H^3$ corresponds to a de Sitter ``half-space'' $v^*$.

For a closed convex subset $K \subset \H^3$, its polar dual is defined as
$$
K^* = \{y \in \dS^3 \,|\, \langle x,y \rangle \le 0 \mbox{ for all }x \in K\} = \bigcap_{x \in K} x^*.
$$

\begin{lem}[\cite{RivHod93}, Proposition 2.6]
\label{lem:GaussImageDeSitter}
If $P$ is a compact convex polyhedron in $\H^3$, then the Gauss image of $\partial P$ is isometric to $\partial P^* \subset \dS^3$.
\end{lem}
\begin{proof}
Let $v$ be a vertex of $P$. A parallel translation yields a linear isomorphism
$$
\tau_v \co T_v \H^3 \to v^\perp,
$$
where $v^\perp$ is the orthogonal complement of $v$ in $\R^{3,1}$. By definition, $\tau_v$ preserves the scalar products. Thus the de Sitter space is a perfect habitat for the Gauss image $\Pi_v$ of the vertex $v$:
$$
\tau_v(\Pi_v) \subset v^\perp \cap \dS^3 = \partial v^*.
$$
It is not hard to show that
$$
\bigcup_v \tau_v(\Pi_v) = \partial P^*,
$$
which proves the lemma.
\end{proof}

Through Lorentz transformations, isometries of $\H^3$ correspond to isometries of $\dS^3$, so that it makes sense to speak of parabolic de Sitter isometries. Similarly to a convex parabolic polyhedron in $\H^3$ (see Section \ref{subsec:GaussImage}), a \emph{convex parabolic de Sitter polyhedron} is defined as the convex hull of finitely many orbits of a $\Z^2$-action on $\dS^3$ by parabolic isometries. (We cheat a little when saying ``the convex hull'', since the de Sitter space is not convex, and there is no convex subset of $\dS^3$ invariant under a parabolic isometry. Let us call a subset of the de Sitter space convex if it is the intersection of $\dS^3$ with a convex cone in the Minkowski space.)

The de Sitter space can be visualized in the Klein projective model of $\H^3$ as the exterior of the absolute quadric. Then the duality with respect to the Minkowski product becomes the polarity with respect to the quadric. For more details, see \cite{Fill3}.




It is easy to show that the polar dual of a convex parabolic polyhedron in $\H^3$ is a convex parabolic de Sitter polyhedron. Since Lemma \ref{lem:GaussImageDeSitter} has a straightforward generalization to parabolic polyhedra, Theorem \ref{thm:Main} can be reformulated as follows (see also Theorem~\ref{thmrivinprime}).

\begin{Alphatheoremprime}
Let $g$ be a spherical cone metric on the torus $\Tor^2$ such that the following two conditions hold:
\begin{enumerate}
\item All cone angles of $(\Tor^2, g)$ are greater than $2\pi$.
\item There are no contractible closed geodesics on $(\Tor^2, g)$ of length less or equal $2\pi$.
\end{enumerate}
Then the universal cover of $(\Tor^2, g)$ can be isometrically embedded in the de Sitter space $\dS^3$ as a convex polyhedral surface invariant under a parabolic action of $\Z^2$. This embedding is unique up to isometry.
\end{Alphatheoremprime}

Uniqueness in the statement above must be understood as uniqueness among convex parabolic polyhedra.

Our proof of Theorem \ref{thm:Main} can also be interpreted in terms of de Sitter geometry. The objects dual to hyperbolic cusps with \emph{coparticles} are de Sitter cusps with \emph{particles}. They are constructed as follows. Let $(L_1, L_2, L_3; o)$ be a corner. Since co-oriented planes are in a one-to-one correspondence with half-spaces, the equation \eqref{eqn:wDual} associates a point $w_i \in \dS^3$ to each of $L_i$. Consider the tetrahedron with vertices $w_1$, $w_2$, $w_3$, $o$. Its intersection with the de Sitter space is called a \emph{de Sitter horoprism}. The operation dual to fitting corners together (see Section \ref{subsec:CuspFromCorn}) is side-to-side gluing of de Sitter horoprisms. A cusp with particles is a Lorentzian cone-manifold  with time-like singular lines ending at the vertices. In these terms Lemmas \ref{lem:DTWD}, \ref{lem:lijWD}, \ref{lem:LinkWD} become more obvious.

Due to Lemma \ref{lem:GaussImageDeSitter}, the space $\M^*$ can be viewed as the space of all convex de Sitter cusps with particles with fixed boundary metric. On $\M^*$, we consider the function
$$
V^*(M^*) = -2 \covol(M^*) + \sum_i h^*_i \kappa_i + \sum_{ij} \ell^*_{ij}(\pi - \theta^*_{ij}).
$$
Here $h^*_i$ are the (imaginary parts of the) heights of the de Sitter horoprisms; $\ell^*_{ij}$ and $\theta^*_{ij}$ are the lengths of and the dihedral angles at the boundary edges of $M^*$; finally $\covol(M^*)$ is the total volume of semi-ideal de Sitter simplices associated with the horoprisms. Note that $h^*_i =  -h_i$. By the Schl\"afli formula in $\dS^3$, \cite{SP00}, we have
$$
\frac{\partial V^*}{\partial h^*_i} = \kappa_i.
$$
This implies the equation $V^*(M^*) + V(M) = \mathrm{const}$, see also \cite[Proposition 2.1]{SP00}. The functional $V^*$ is a discrete analog of the Hilbert-Einstein functional, see a discussion in the middle of Section \ref{subsec:RelRes}.

\begin{appendix}

\section{Some trigonometry}
\label{sec:FormCorn}

\begin{figure}[ht]
\begin{center}
\input{Fig/HypKites.tex}
\end{center}
\caption{Metric elements in hyperbolic kites.}
\label{fig:HypKites}
\end{figure}
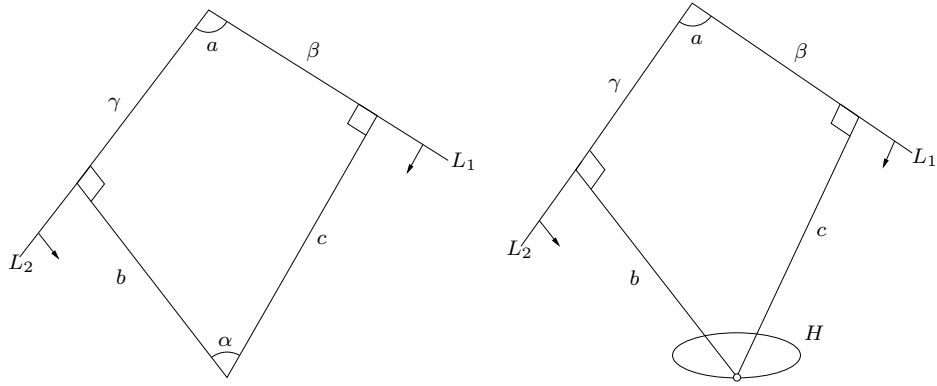

Let us call a \emph{kite} a quadrilateral with two opposite right angles. We will study two types of hyperbolic kites: compact ones, left of Figure \ref{fig:HypKites}, and semi-ideal ones, right of Figure \ref{fig:HypKites}. We denote the edge lengths and angle values in kites as shown on Figure \ref{fig:HypKites}. The lengths $b$ and $c$ in a semi-ideal kite are defined as distances to a horosphere $H$ centered at the ideal vertex. All edge lengths in a kite are allowed to be negative. Formally, a kite is a triple $(L_1, L_2; v)$, where $L_1$ and $L_2$ are intersecting co-oriented lines, and $v$ is a point from which perpendiculars to the lines are dropped. If $v$ is an ideal point, then we require that it lies on the positive sides of both $L_1$ and~$L_2$.

\begin{lem}
\label{lem:SinCosLaw}
For a compact hyperbolic kite, the following holds.
\begin{itemize}
\item
The sine law:
$$
\frac{\sin a}{\sin\alpha} = \frac{\cosh b}{\cosh\beta} = \frac{\cosh c}{\cosh\gamma}.
$$
\item
The cosine laws:
\begin{eqnarray*}
\cos a & = & \sinh b \sinh c - \cosh b \cosh c \cos\alpha,\\
\sinh c &  = & - \cos a \sinh b  + \sin a \cosh b \sinh\gamma,\\
\end{eqnarray*}
\item
When $c$ varies, we have
$$
\frac{\partial\alpha}{\partial c} = - \frac{\tanh\beta}{\cosh c}.
$$
\end{itemize}

For a semi-ideal kite, the following holds.
\begin{itemize}
\item
The sine law:
$$
\frac{e^b}{\cosh\beta} = \frac{e^c}{\cosh\gamma}.
$$
\item
The cosine laws:
\begin{eqnarray*}
e^{c-b} & = & - \cos a + \sin a \sinh\gamma\\
1 & = & \sinh \beta \sinh \gamma - \cosh \beta \cosh \gamma \cos a.
\end{eqnarray*}
\item
When $c$ varies, we have
$$
\frac{\partial\gamma}{\partial c} = \frac{e^{c-b}}{\cosh\gamma \sin a} = \frac{1}{\cosh\beta \sin a}.
$$
\end{itemize}
\end{lem}
\begin{proof}
By polarity, co-oriented lines in $\H^2$ correspond to points in the de Sitter plane $\dS^2$. Let $w_1$ and $w_2$ be the poles of the lines $L_1$ and $L_2$, respectively. Then we have a triangle $vw_1w_2$ in $\H^2 \cup \dS^2$, and formulas of Theorem \ref{lem:SinCosLaw} are just sine and cosine laws for this triangle. The formulas can be proved using an argument similar to that in \cite[Section 2.4]{Thu97}. For a compact kite, this is done in \cite[end of Section 4.3]{Cho07}.
\end{proof}

Faces of bricks, see Section \ref{subsec:Bricks}, are hyperbolic kites. By applying formulas of Lemma \ref{lem:SinCosLaw} to them, we obtain the following.

\begin{equation}
\label{eqn:h12}
\sinh h_{12} = \frac{e^{h_2-h_1} - \cos\alpha_{12}}{\sin\alpha_{12}}
\end{equation}

\begin{equation}
\label{eqn:h123}
\sinh h_{123} = \frac{-\cos\gamma_1^{23}\sinh h_{12} + \sinh h_{13}}{\sin\gamma_1^{23} \cosh h_{12}}
\end{equation}

(Note that $h_{123} = h_{213}$ by construction, although it is not obvious in the above formula.)

\begin{equation}
\label{eqn:OmegaFromH}
\cos \omega_1^{23} = \frac{\sinh h_{12} \sinh h_{13} + \cos\gamma_1^{23}}{\cosh h_{12} \cosh h_{13}}
\end{equation}

\begin{equation}
\label{eqn:DerOmegaH}
\frac{\partial\omega_1^{23}}{\partial h_2} = \frac{\partial\omega_1^{23}}{\partial h_{12}} \frac{\partial h_{12}}{\partial h_2} = - \frac{\tanh h_{123}} {\sin\alpha_{12} \cosh h_{12} \cosh h_{21}}
\end{equation}

\end{appendix}

\bibliographystyle{alpha}

\input{HypCuspDual.bbl}
\end{document}

%% file: Fig/ComplAngles.tex
\begin{picture}(0,0)%
\includegraphics{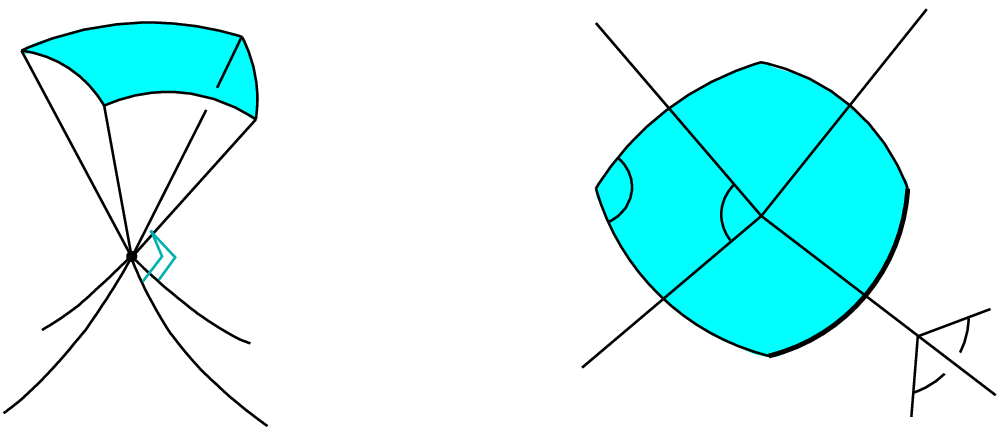}%
\end{picture}%
\setlength{\unitlength}{5801sp}%
\begingroup\makeatletter\ifx\SetFigFont\undefined%
\gdef\SetFigFont#1#2#3#4#5{%
  \reset@font\fontsize{#1}{#2pt}%
  \fontfamily{#3}\fontseries{#4}\fontshape{#5}%
  \selectfont}%
\fi\endgroup%
\begin{picture}(3263,1402)(575,-1406)
\put(3606,-1364){\makebox(0,0)[lb]{\smash{{\SetFigFont{8}{9.6}{\rmdefault}{\mddefault}{\updefault}{\color[rgb]{0,0,0}$\pi-\alpha$}%
}}}}
\put(935,-213){\makebox(0,0)[lb]{\smash{{\SetFigFont{8}{9.6}{\rmdefault}{\mddefault}{\updefault}{\color[rgb]{0,0,0}$\Pi_v$}%
}}}}
\put(2648,-565){\makebox(0,0)[lb]{\smash{{\SetFigFont{8}{9.6}{\rmdefault}{\mddefault}{\updefault}{\color[rgb]{0,0,0}$\gamma$}%
}}}}
\put(2691,-765){\makebox(0,0)[lb]{\smash{{\SetFigFont{8}{9.6}{\rmdefault}{\mddefault}{\updefault}{\color[rgb]{0,0,0}$\pi-\gamma$}%
}}}}
\put(3172,-1053){\makebox(0,0)[lb]{\smash{{\SetFigFont{8}{9.6}{\rmdefault}{\mddefault}{\updefault}{\color[rgb]{0,0,0}$\alpha$}%
}}}}
\put(887,-858){\makebox(0,0)[lb]{\smash{{\SetFigFont{8}{9.6}{\rmdefault}{\mddefault}{\updefault}{\color[rgb]{0,0,0}$v$}%
}}}}
\put(1308,-1211){\makebox(0,0)[lb]{\smash{{\SetFigFont{8}{9.6}{\rmdefault}{\mddefault}{\updefault}{\color[rgb]{0,0,0}$F_i$}%
}}}}
\put(1449,-425){\makebox(0,0)[lb]{\smash{{\SetFigFont{8}{9.6}{\rmdefault}{\mddefault}{\updefault}{\color[rgb]{0,0,0}$i$}%
}}}}
\end{picture}%

%% file: Fig/OneVertCusp.tex
\begin{picture}(0,0)%
\includegraphics{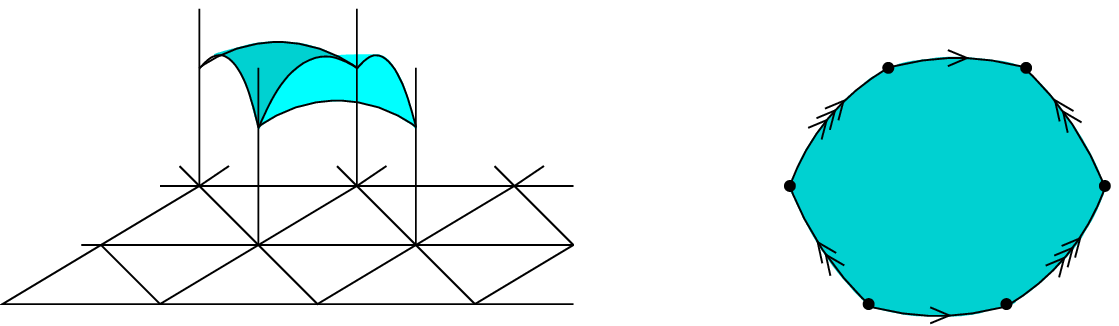}%
\end{picture}%
\setlength{\unitlength}{4144sp}%
\begingroup\makeatletter\ifx\SetFigFont\undefined%
\gdef\SetFigFont#1#2#3#4#5{%
  \reset@font\fontsize{#1}{#2pt}%
  \fontfamily{#3}\fontseries{#4}\fontshape{#5}%
  \selectfont}%
\fi\endgroup%
\begin{picture}(5083,1462)(2239,-2501)
\put(3935,-1187){\makebox(0,0)[lb]{\smash{{\SetFigFont{9}{10.8}{\rmdefault}{\mddefault}{\updefault}{\color[rgb]{0,0,0}$c$}%
}}}}
\put(4186,-1726){\makebox(0,0)[lb]{\smash{{\SetFigFont{9}{10.8}{\rmdefault}{\mddefault}{\updefault}{\color[rgb]{0,0,0}$b$}%
}}}}
\put(2971,-1321){\makebox(0,0)[lb]{\smash{{\SetFigFont{9}{10.8}{\rmdefault}{\mddefault}{\updefault}{\color[rgb]{0,0,0}$d$}%
}}}}
\put(3466,-1726){\makebox(0,0)[lb]{\smash{{\SetFigFont{9}{10.8}{\rmdefault}{\mddefault}{\updefault}{\color[rgb]{0,0,0}$a$}%
}}}}
\end{picture}%

%% file: Fig/2Corners.tex
\begin{picture}(0,0)%
\includegraphics{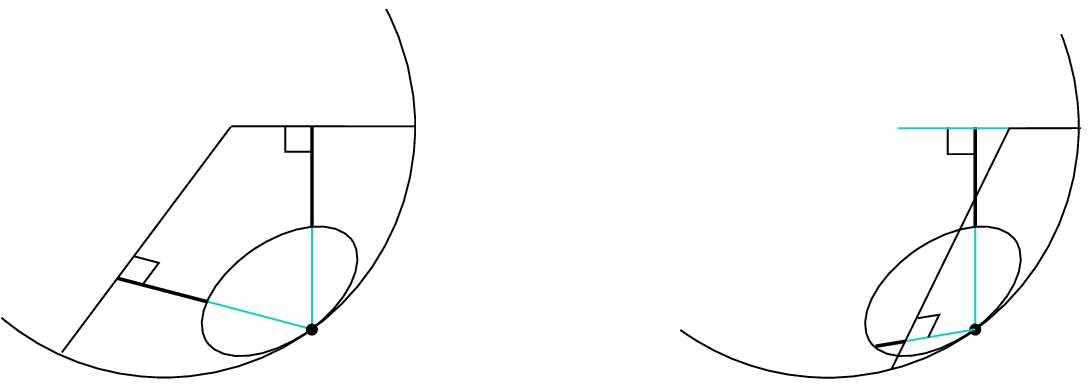}%
\end{picture}%
\setlength{\unitlength}{4144sp}%
\begingroup\makeatletter\ifx\SetFigFontNFSS\undefined%
\gdef\SetFigFontNFSS#1#2#3#4#5{%
  \reset@font\fontsize{#1}{#2pt}%
  \fontfamily{#3}\fontseries{#4}\fontshape{#5}%
  \selectfont}%
\fi\endgroup%
\begin{picture}(4996,1840)(2888,-2699)
\put(7869,-1390){\makebox(0,0)[lb]{\smash{{\SetFigFontNFSS{8}{9.6}{\rmdefault}{\mddefault}{\updefault}{\color[rgb]{0,0,0}$L_1$}%
}}}}
\put(4350,-2429){\makebox(0,0)[lb]{\smash{{\SetFigFontNFSS{8}{9.6}{\rmdefault}{\mddefault}{\updefault}{\color[rgb]{0,0,0}$o$}%
}}}}
\put(6666,-2424){\makebox(0,0)[lb]{\smash{{\SetFigFontNFSS{8}{9.6}{\rmdefault}{\mddefault}{\updefault}{\color[rgb]{0,0,0}$h_2$}%
}}}}
\put(7187,-1673){\makebox(0,0)[lb]{\smash{{\SetFigFontNFSS{8}{9.6}{\rmdefault}{\mddefault}{\updefault}{\color[rgb]{0,0,0}$h_1$}%
}}}}
\put(3514,-2260){\makebox(0,0)[lb]{\smash{{\SetFigFontNFSS{8}{9.6}{\rmdefault}{\mddefault}{\updefault}{\color[rgb]{0,0,0}$h_2$}%
}}}}
\put(4340,-1680){\makebox(0,0)[lb]{\smash{{\SetFigFontNFSS{8}{9.6}{\rmdefault}{\mddefault}{\updefault}{\color[rgb]{0,0,0}$h_1$}%
}}}}
\put(3042,-2561){\makebox(0,0)[lb]{\smash{{\SetFigFontNFSS{8}{9.6}{\rmdefault}{\mddefault}{\updefault}{\color[rgb]{0,0,0}$L_2$}%
}}}}
\put(4821,-1412){\makebox(0,0)[lb]{\smash{{\SetFigFontNFSS{8}{9.6}{\rmdefault}{\mddefault}{\updefault}{\color[rgb]{0,0,0}$L_1$}%
}}}}
\put(7377,-2429){\makebox(0,0)[lb]{\smash{{\SetFigFontNFSS{8}{9.6}{\rmdefault}{\mddefault}{\updefault}{\color[rgb]{0,0,0}$o$}%
}}}}
\put(6905,-2653){\makebox(0,0)[lb]{\smash{{\SetFigFontNFSS{8}{9.6}{\rmdefault}{\mddefault}{\updefault}{\color[rgb]{0,0,0}$L_2$}%
}}}}
\end{picture}%

%% file: Fig/l12.tex
\begin{picture}(0,0)%
\includegraphics{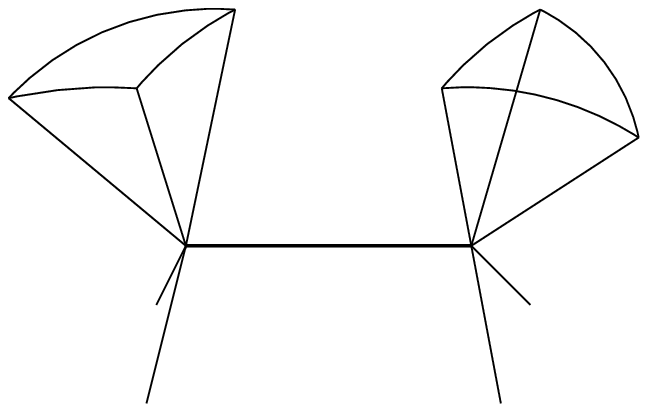}%
\end{picture}%
\setlength{\unitlength}{4144sp}%
\begingroup\makeatletter\ifx\SetFigFont\undefined%
\gdef\SetFigFont#1#2#3#4#5{%
  \reset@font\fontsize{#1}{#2pt}%
  \fontfamily{#3}\fontseries{#4}\fontshape{#5}%
  \selectfont}%
\fi\endgroup%
\begin{picture}(3090,1911)(1966,-2593)
\put(3556,-2176){\makebox(0,0)[lb]{\smash{{\SetFigFont{8}{9.6}{\rmdefault}{\mddefault}{\updefault}{\color[rgb]{0,0,0}$L_2$}%
}}}}
\put(3196,-781){\makebox(0,0)[lb]{\smash{{\SetFigFont{8}{9.6}{\rmdefault}{\mddefault}{\updefault}{\color[rgb]{0,0,0}$2$}%
}}}}
\put(4591,-781){\makebox(0,0)[lb]{\smash{{\SetFigFont{8}{9.6}{\rmdefault}{\mddefault}{\updefault}{\color[rgb]{0,0,0}$2$}%
}}}}
\put(5041,-1411){\makebox(0,0)[lb]{\smash{{\SetFigFont{8}{9.6}{\rmdefault}{\mddefault}{\updefault}{\color[rgb]{0,0,0}$3$}%
}}}}
\put(3376,-2536){\makebox(0,0)[lb]{\smash{{\SetFigFont{8}{9.6}{\rmdefault}{\mddefault}{\updefault}{\color[rgb]{0,0,0}$L_1$}%
}}}}
\put(2611,-2311){\makebox(0,0)[lb]{\smash{{\SetFigFont{8}{9.6}{\rmdefault}{\mddefault}{\updefault}{\color[rgb]{0,0,0}$L_4$}%
}}}}
\put(1981,-1231){\makebox(0,0)[lb]{\smash{{\SetFigFont{8}{9.6}{\rmdefault}{\mddefault}{\updefault}{\color[rgb]{0,0,0}$4$}%
}}}}
\put(2773,-1190){\makebox(0,0)[lb]{\smash{{\SetFigFont{8}{9.6}{\rmdefault}{\mddefault}{\updefault}{\color[rgb]{0,0,0}$1$}%
}}}}
\put(3981,-1191){\makebox(0,0)[lb]{\smash{{\SetFigFont{8}{9.6}{\rmdefault}{\mddefault}{\updefault}{\color[rgb]{0,0,0}$1$}%
}}}}
\put(4375,-2311){\makebox(0,0)[lb]{\smash{{\SetFigFont{8}{9.6}{\rmdefault}{\mddefault}{\updefault}{\color[rgb]{0,0,0}$L_3$}%
}}}}
\put(3421,-1788){\makebox(0,0)[lb]{\smash{{\SetFigFont{8}{9.6}{\rmdefault}{\mddefault}{\updefault}{\color[rgb]{0,0,0}$\ell_{12}$}%
}}}}
\end{picture}%

%% file: Fig/SingFace.tex
\begin{picture}(0,0)%
\includegraphics{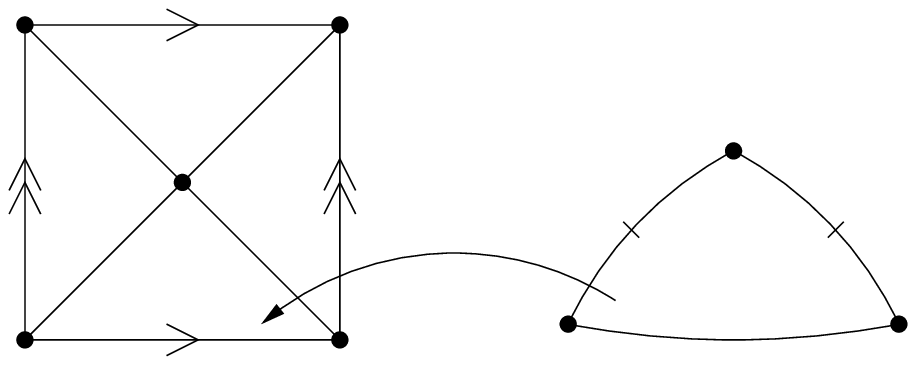}%
\end{picture}%
\setlength{\unitlength}{3315sp}%
\begingroup\makeatletter\ifx\SetFigFont\undefined%
\gdef\SetFigFont#1#2#3#4#5{%
  \reset@font\fontsize{#1}{#2pt}%
  \fontfamily{#3}\fontseries{#4}\fontshape{#5}%
  \selectfont}%
\fi\endgroup%
\begin{picture}(5288,2175)(211,-1606)
\put(1306,-376){\makebox(0,0)[lb]{\smash{{\SetFigFont{8}{9.6}{\rmdefault}{\mddefault}{\updefault}{\color[rgb]{0,0,0}$v_1$}%
}}}}
\put(226,434){\makebox(0,0)[lb]{\smash{{\SetFigFont{8}{9.6}{\rmdefault}{\mddefault}{\updefault}{\color[rgb]{0,0,0}$v_2$}%
}}}}
\put(271,-1546){\makebox(0,0)[lb]{\smash{{\SetFigFont{8}{9.6}{\rmdefault}{\mddefault}{\updefault}{\color[rgb]{0,0,0}$v_2$}%
}}}}
\put(2341,434){\makebox(0,0)[lb]{\smash{{\SetFigFont{8}{9.6}{\rmdefault}{\mddefault}{\updefault}{\color[rgb]{0,0,0}$v_2$}%
}}}}
\put(2341,-1501){\makebox(0,0)[lb]{\smash{{\SetFigFont{8}{9.6}{\rmdefault}{\mddefault}{\updefault}{\color[rgb]{0,0,0}$v_2$}%
}}}}
\end{picture}%

%% file: Fig/Quad.tex
\begin{picture}(0,0)%
\includegraphics{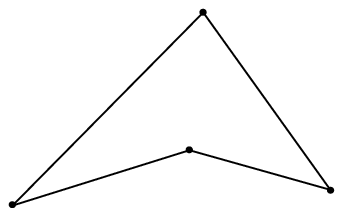}%
\end{picture}%
\setlength{\unitlength}{4144sp}%
\begingroup\makeatletter\ifx\SetFigFont\undefined%
\gdef\SetFigFont#1#2#3#4#5{%
  \reset@font\fontsize{#1}{#2pt}%
  \fontfamily{#3}\fontseries{#4}\fontshape{#5}%
  \selectfont}%
\fi\endgroup%
\begin{picture}(1614,1254)(5976,-1354)
\put(6924,-211){\makebox(0,0)[lb]{\smash{{\SetFigFont{8}{9.6}{\familydefault}{\mddefault}{\updefault}{\color[rgb]{0,0,0}$j$}%
}}}}
\put(7575,-1207){\makebox(0,0)[lb]{\smash{{\SetFigFont{8}{9.6}{\familydefault}{\mddefault}{\updefault}{\color[rgb]{0,0,0}$k$}%
}}}}
\put(6847,-1060){\makebox(0,0)[lb]{\smash{{\SetFigFont{8}{9.6}{\familydefault}{\mddefault}{\updefault}{\color[rgb]{0,0,0}$i$}%
}}}}
\put(5991,-1308){\makebox(0,0)[lb]{\smash{{\SetFigFont{8}{9.6}{\familydefault}{\mddefault}{\updefault}{\color[rgb]{0,0,0}$l$}%
}}}}
\end{picture}%

%% file: Fig/Brick.tex
\begin{picture}(0,0)%
\includegraphics{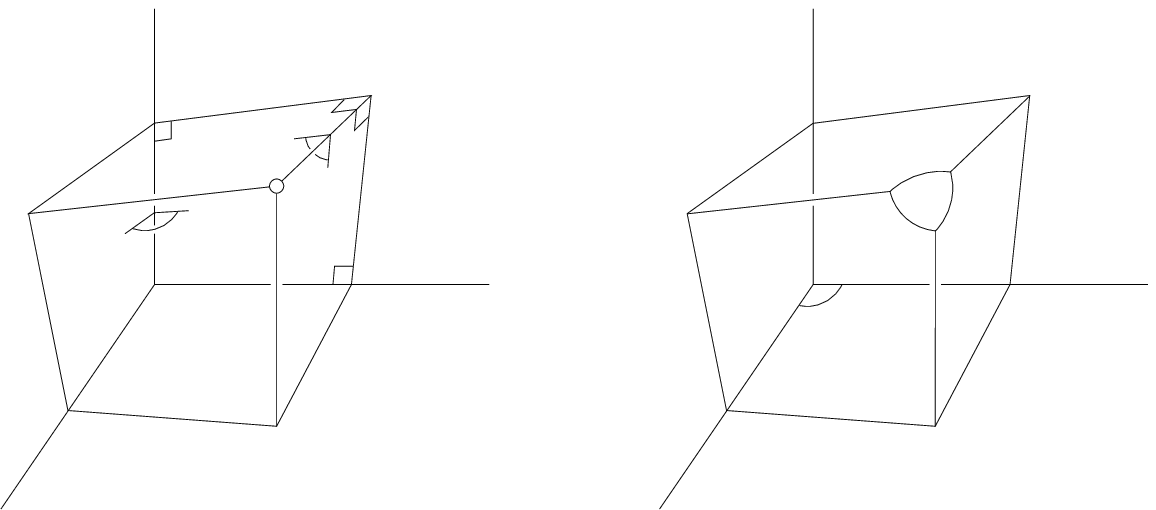}%
\end{picture}%
\setlength{\unitlength}{1657sp}%
\begingroup\makeatletter\ifx\SetFigFont\undefined%
\gdef\SetFigFont#1#2#3#4#5{%
  \reset@font\fontsize{#1}{#2pt}%
  \fontfamily{#3}\fontseries{#4}\fontshape{#5}%
  \selectfont}%
\fi\endgroup%
\begin{picture}(13131,5739)(1429,-5923)
\put(1486,-961){\makebox(0,0)[lb]{\smash{{\SetFigFont{8}{9.6}{\rmdefault}{\mddefault}{\updefault}{\color[rgb]{0,0,0}$L_j$}%
}}}}
\put(2056,-3076){\makebox(0,0)[lb]{\smash{{\SetFigFont{8}{9.6}{\rmdefault}{\mddefault}{\updefault}{\color[rgb]{0,0,0}$\pi - \alpha_{ij}$}%
}}}}
\put(4094,-1936){\makebox(0,0)[lb]{\smash{{\SetFigFont{8}{9.6}{\rmdefault}{\mddefault}{\updefault}{\color[rgb]{0,0,0}$\omega_{i}^{jk}$}%
}}}}
\put(5221,-5671){\makebox(0,0)[lb]{\smash{{\SetFigFont{8}{9.6}{\rmdefault}{\mddefault}{\updefault}{\color[rgb]{0,0,0}$L_k$}%
}}}}
\put(12973,-5686){\makebox(0,0)[lb]{\smash{{\SetFigFont{8}{9.6}{\rmdefault}{\mddefault}{\updefault}{\color[rgb]{0,0,0}$L_k$}%
}}}}
\put(9013,-961){\makebox(0,0)[lb]{\smash{{\SetFigFont{8}{9.6}{\rmdefault}{\mddefault}{\updefault}{\color[rgb]{0,0,0}$L_j$}%
}}}}
\put(6326,-646){\makebox(0,0)[lb]{\smash{{\SetFigFont{8}{9.6}{\rmdefault}{\mddefault}{\updefault}{\color[rgb]{0,0,0}$L_i$}%
}}}}
\put(13753,-726){\makebox(0,0)[lb]{\smash{{\SetFigFont{8}{9.6}{\rmdefault}{\mddefault}{\updefault}{\color[rgb]{0,0,0}$L_i$}%
}}}}
\put(5799,-1138){\makebox(0,0)[lb]{\smash{{\SetFigFont{8}{9.6}{\rmdefault}{\mddefault}{\updefault}{\color[rgb]{0,0,0}$p_i$}%
}}}}
\put(13193,-2447){\makebox(0,0)[lb]{\smash{{\SetFigFont{8}{9.6}{\rmdefault}{\mddefault}{\updefault}{\color[rgb]{0,0,0}$h_{ik}$}%
}}}}
\put(12533,-2064){\makebox(0,0)[lb]{\smash{{\SetFigFont{8}{9.6}{\rmdefault}{\mddefault}{\updefault}{\color[rgb]{0,0,0}$h_i$}%
}}}}
\put(10648,-3868){\makebox(0,0)[lb]{\smash{{\SetFigFont{8}{9.6}{\rmdefault}{\mddefault}{\updefault}{\color[rgb]{0,0,0}$\pi-\gamma_k^{ij}$}%
}}}}
\put(3200,-3620){\makebox(0,0)[lb]{\smash{{\SetFigFont{8}{9.6}{\rmdefault}{\mddefault}{\updefault}{\color[rgb]{0,0,0}$v$}%
}}}}
\put(11821,-2513){\makebox(0,0)[lb]{\smash{{\SetFigFont{8}{9.6}{\rmdefault}{\mddefault}{\updefault}{\color[rgb]{0,0,0}$H$}%
}}}}
\put(11583,-1204){\makebox(0,0)[lb]{\smash{{\SetFigFont{8}{9.6}{\rmdefault}{\mddefault}{\updefault}{\color[rgb]{0,0,0}$h_{ij}$}%
}}}}
\put(5579,-3228){\makebox(0,0)[lb]{\smash{{\SetFigFont{8}{9.6}{\rmdefault}{\mddefault}{\updefault}{\color[rgb]{0,0,0}$p_{ik}$}%
}}}}
\put(2704,-1458){\makebox(0,0)[lb]{\smash{{\SetFigFont{8}{9.6}{\rmdefault}{\mddefault}{\updefault}{\color[rgb]{0,0,0}$p_{ij}$}%
}}}}
\put(10793,-2924){\makebox(0,0)[lb]{\smash{{\SetFigFont{8}{9.6}{\rmdefault}{\mddefault}{\updefault}{\color[rgb]{0,0,0}$h_{ijk}$}%
}}}}
\put(4301,-2515){\makebox(0,0)[lb]{\smash{{\SetFigFont{8}{9.6}{\rmdefault}{\mddefault}{\updefault}{\color[rgb]{0,0,0}$o$}%
}}}}
\end{picture}%

%% file: Fig/FaceI.tex
\begin{picture}(0,0)%
\includegraphics{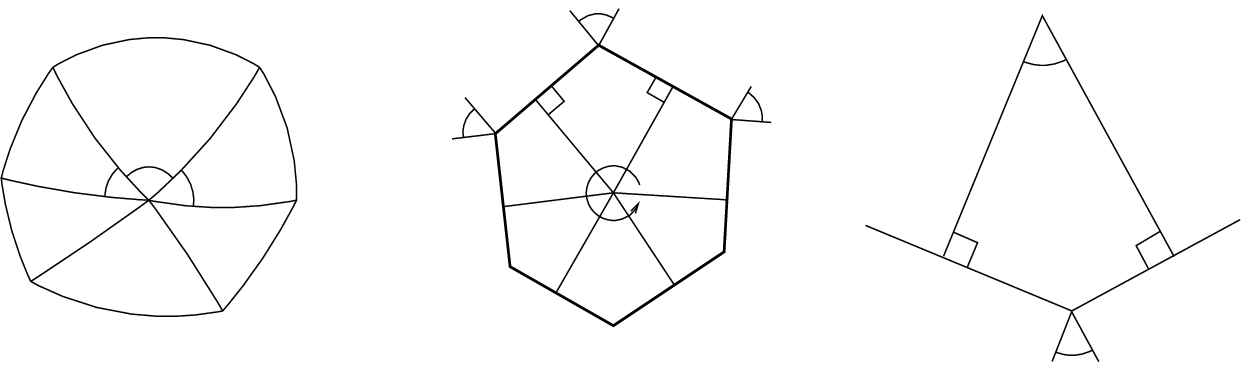}%
\end{picture}%
\setlength{\unitlength}{3108sp}%
\begingroup\makeatletter\ifx\SetFigFont\undefined%
\gdef\SetFigFont#1#2#3#4#5{%
  \reset@font\fontsize{#1}{#2pt}%
  \fontfamily{#3}\fontseries{#4}\fontshape{#5}%
  \selectfont}%
\fi\endgroup%
\begin{picture}(7572,2524)(443,-2289)
\put(5135,-557){\makebox(0,0)[lb]{\smash{{\SetFigFont{8}{9.6}{\rmdefault}{\mddefault}{\updefault}{\color[rgb]{0,0,0}$\gamma^{1n}$}%
}}}}
\put(3957, 88){\makebox(0,0)[lb]{\smash{{\SetFigFont{8}{9.6}{\rmdefault}{\mddefault}{\updefault}{\color[rgb]{0,0,0}$\gamma^{12}$}%
}}}}
\put(4098,-1339){\makebox(0,0)[lb]{\smash{{\SetFigFont{8}{9.6}{\rmdefault}{\mddefault}{\updefault}{\color[rgb]{0,0,0}$\omega_i$}%
}}}}
\put(3008,-659){\makebox(0,0)[lb]{\smash{{\SetFigFont{8}{9.6}{\rmdefault}{\mddefault}{\updefault}{\color[rgb]{0,0,0}$\gamma^{23}$}%
}}}}
\put(4413,-752){\makebox(0,0)[lb]{\smash{{\SetFigFont{8}{9.6}{\rmdefault}{\mddefault}{\updefault}{\color[rgb]{0,0,0}$g_1$}%
}}}}
\put(3713,-782){\makebox(0,0)[lb]{\smash{{\SetFigFont{8}{9.6}{\rmdefault}{\mddefault}{\updefault}{\color[rgb]{0,0,0}$g_2$}%
}}}}
\put(6890,-2225){\makebox(0,0)[lb]{\smash{{\SetFigFont{8}{9.6}{\rmdefault}{\mddefault}{\updefault}{\color[rgb]{0,0,0}$\gamma^{j,j+1}$}%
}}}}
\put(6961,-219){\makebox(0,0)[lb]{\smash{{\SetFigFont{8}{9.6}{\rmdefault}{\mddefault}{\updefault}{\color[rgb]{0,0,0}$\omega^{j,j+1}$}%
}}}}
\put(6328,-1736){\makebox(0,0)[lb]{\smash{{\SetFigFont{8}{9.6}{\rmdefault}{\mddefault}{\updefault}{\color[rgb]{0,0,0}$g_{j,j+1}$}%
}}}}
\put(7287,-1727){\makebox(0,0)[lb]{\smash{{\SetFigFont{8}{9.6}{\rmdefault}{\mddefault}{\updefault}{\color[rgb]{0,0,0}$g_{j+1,j}$}%
}}}}
\put(6305,-723){\makebox(0,0)[lb]{\smash{{\SetFigFont{8}{9.6}{\rmdefault}{\mddefault}{\updefault}{\color[rgb]{0,0,0}$g_j$}%
}}}}
\put(7287,-713){\makebox(0,0)[lb]{\smash{{\SetFigFont{8}{9.6}{\rmdefault}{\mddefault}{\updefault}{\color[rgb]{0,0,0}$g_{j+1}$}%
}}}}
\put(2061,-255){\makebox(0,0)[lb]{\smash{{\SetFigFont{8}{9.6}{\rmdefault}{\mddefault}{\updefault}{\color[rgb]{0,0,0}$1$}%
}}}}
\put(2305,-1165){\makebox(0,0)[lb]{\smash{{\SetFigFont{8}{9.6}{\rmdefault}{\mddefault}{\updefault}{\color[rgb]{0,0,0}$n$}%
}}}}
\put(1684,-1032){\makebox(0,0)[lb]{\smash{{\SetFigFont{8}{9.6}{\rmdefault}{\mddefault}{\updefault}{\color[rgb]{0,0,0}$\gamma^{1n}$}%
}}}}
\put(623,-290){\makebox(0,0)[lb]{\smash{{\SetFigFont{8}{9.6}{\rmdefault}{\mddefault}{\updefault}{\color[rgb]{0,0,0}$2$}%
}}}}
\put(1222,-823){\makebox(0,0)[lb]{\smash{{\SetFigFont{8}{9.6}{\rmdefault}{\mddefault}{\updefault}{\color[rgb]{0,0,0}$\gamma^{12}$}%
}}}}
\put(836,-983){\makebox(0,0)[lb]{\smash{{\SetFigFont{8}{9.6}{\rmdefault}{\mddefault}{\updefault}{\color[rgb]{0,0,0}$\gamma^{23}$}%
}}}}
\put(1265,-1299){\makebox(0,0)[lb]{\smash{{\SetFigFont{8}{9.6}{\rmdefault}{\mddefault}{\updefault}{\color[rgb]{0,0,0}$i$}%
}}}}
\end{picture}%

%% file: Fig/CutFace.tex
\begin{picture}(0,0)%
\includegraphics{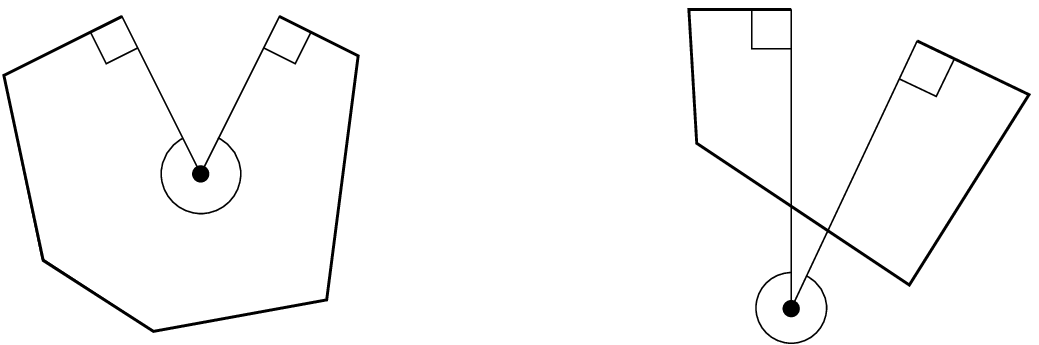}%
\end{picture}%
\setlength{\unitlength}{3315sp}%
\begingroup\makeatletter\ifx\SetFigFont\undefined%
\gdef\SetFigFont#1#2#3#4#5{%
  \reset@font\fontsize{#1}{#2pt}%
  \fontfamily{#3}\fontseries{#4}\fontshape{#5}%
  \selectfont}%
\fi\endgroup%
\begin{picture}(5903,2133)(654,-1488)
\put(1585, -2){\makebox(0,0)[lb]{\smash{{\SetFigFont{8}{9.6}{\rmdefault}{\mddefault}{\updefault}{\color[rgb]{0,0,0}$g_1$}%
}}}}
\put(2025,-162){\makebox(0,0)[lb]{\smash{{\SetFigFont{8}{9.6}{\rmdefault}{\mddefault}{\updefault}{\color[rgb]{0,0,0}$g_1$}%
}}}}
\put(2505,344){\makebox(0,0)[lb]{\smash{{\SetFigFont{8}{9.6}{\rmdefault}{\mddefault}{\updefault}{\color[rgb]{0,0,0}$g_{1n}$}%
}}}}
\put(798,258){\makebox(0,0)[lb]{\smash{{\SetFigFont{8}{9.6}{\rmdefault}{\mddefault}{\updefault}{\color[rgb]{0,0,0}$g_{12}$}%
}}}}
\put(5565,-582){\makebox(0,0)[lb]{\smash{{\SetFigFont{8}{9.6}{\rmdefault}{\mddefault}{\updefault}{\color[rgb]{0,0,0}$g_1$}%
}}}}
\put(6218,144){\makebox(0,0)[lb]{\smash{{\SetFigFont{8}{9.6}{\rmdefault}{\mddefault}{\updefault}{\color[rgb]{0,0,0}$g_{1n}$}%
}}}}
\put(4719,498){\makebox(0,0)[lb]{\smash{{\SetFigFont{8}{9.6}{\rmdefault}{\mddefault}{\updefault}{\color[rgb]{0,0,0}$g_{12}$}%
}}}}
\put(5212,-395){\makebox(0,0)[lb]{\smash{{\SetFigFont{8}{9.6}{\rmdefault}{\mddefault}{\updefault}{\color[rgb]{0,0,0}$g_1$}%
}}}}
\put(1713,-853){\makebox(0,0)[lb]{\smash{{\SetFigFont{8}{9.6}{\rmdefault}{\mddefault}{\updefault}{\color[rgb]{0,0,0}$\omega_i$}%
}}}}
\put(4740,-1317){\makebox(0,0)[lb]{\smash{{\SetFigFont{8}{9.6}{\rmdefault}{\mddefault}{\updefault}{\color[rgb]{0,0,0}$\omega_i$}%
}}}}
\end{picture}%

%% file: Fig/BdryAngle.tex
\begin{picture}(0,0)%
\includegraphics{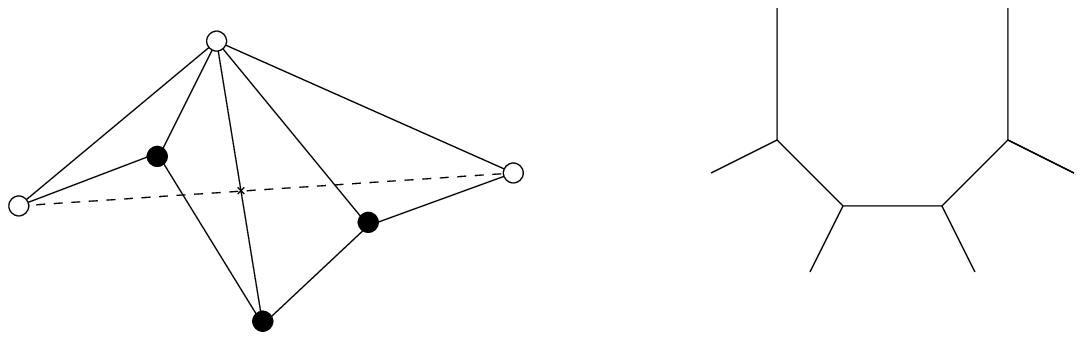}%
\end{picture}%
\setlength{\unitlength}{2776sp}%
\begingroup\makeatletter\ifx\SetFigFontNFSS\undefined%
\gdef\SetFigFontNFSS#1#2#3#4#5{%
  \reset@font\fontsize{#1}{#2pt}%
  \fontfamily{#3}\fontseries{#4}\fontshape{#5}%
  \selectfont}%
\fi\endgroup%
\begin{picture}(7362,2484)(5476,-2251)
\put(11521,-421){\makebox(0,0)[lb]{\smash{{\SetFigFontNFSS{8}{9.6}{\rmdefault}{\mddefault}{\updefault}{\color[rgb]{0,0,0}$i$}%
}}}}
\put(12601,-286){\makebox(0,0)[lb]{\smash{{\SetFigFontNFSS{8}{9.6}{\rmdefault}{\mddefault}{\updefault}{\color[rgb]{0,0,0}$k$}%
}}}}
\put(10351,-286){\makebox(0,0)[lb]{\smash{{\SetFigFontNFSS{8}{9.6}{\rmdefault}{\mddefault}{\updefault}{\color[rgb]{0,0,0}$j$}%
}}}}
\put(6931, 74){\makebox(0,0)[lb]{\smash{{\SetFigFontNFSS{8}{9.6}{\rmdefault}{\mddefault}{\updefault}{\color[rgb]{0,0,0}$i$}%
}}}}
\put(9046,-826){\makebox(0,0)[lb]{\smash{{\SetFigFontNFSS{8}{9.6}{\rmdefault}{\mddefault}{\updefault}{\color[rgb]{0,0,0}$k$}%
}}}}
\put(5491,-1096){\makebox(0,0)[lb]{\smash{{\SetFigFontNFSS{8}{9.6}{\rmdefault}{\mddefault}{\updefault}{\color[rgb]{0,0,0}$j$}%
}}}}
\put(7354,-2182){\makebox(0,0)[lb]{\smash{{\SetFigFontNFSS{8}{9.6}{\rmdefault}{\mddefault}{\updefault}{\color[rgb]{0,0,0}$l$}%
}}}}
\put(7183,-1022){\makebox(0,0)[lb]{\smash{{\SetFigFontNFSS{8}{9.6}{\rmdefault}{\mddefault}{\updefault}{\color[rgb]{0,0,0}$l'$}%
}}}}
\end{picture}%

%% file: Fig/CurveShort1.tex
\begin{picture}(0,0)%
\includegraphics{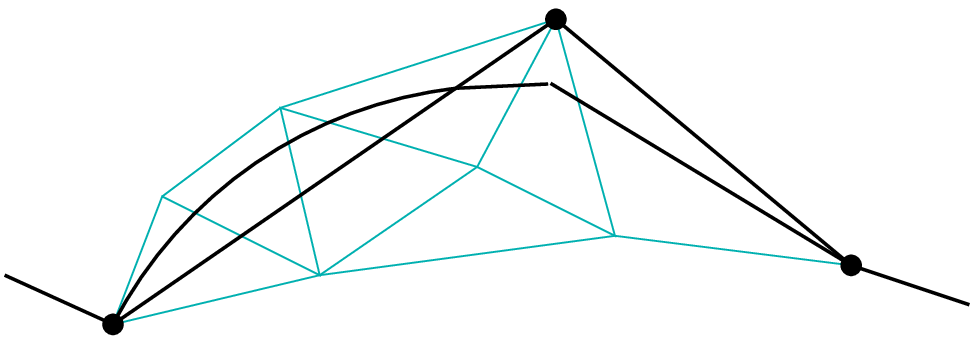}%
\end{picture}%
\setlength{\unitlength}{4144sp}%
\begingroup\makeatletter\ifx\SetFigFontNFSS\undefined%
\gdef\SetFigFontNFSS#1#2#3#4#5{%
  \reset@font\fontsize{#1}{#2pt}%
  \fontfamily{#3}\fontseries{#4}\fontshape{#5}%
  \selectfont}%
\fi\endgroup%
\begin{picture}(4454,1803)(-66,-1291)
\put(2521,389){\makebox(0,0)[lb]{\smash{{\SetFigFontNFSS{9}{10.8}{\rmdefault}{\mddefault}{\updefault}{\color[rgb]{0,0,0}$i$}%
}}}}
\put(2386,-106){\makebox(0,0)[lb]{\smash{{\SetFigFontNFSS{9}{10.8}{\rmdefault}{\mddefault}{\updefault}{\color[rgb]{0,0,0}$z$}%
}}}}
\put(3736,-961){\makebox(0,0)[lb]{\smash{{\SetFigFontNFSS{9}{10.8}{\rmdefault}{\mddefault}{\updefault}{\color[rgb]{0,0,0}$k$}%
}}}}
\put(361,-1231){\makebox(0,0)[lb]{\smash{{\SetFigFontNFSS{9}{10.8}{\rmdefault}{\mddefault}{\updefault}{\color[rgb]{0,0,0}$j$}%
}}}}
\end{picture}%

%% file: Fig/CurveShort2.tex
\begin{picture}(0,0)%
\includegraphics{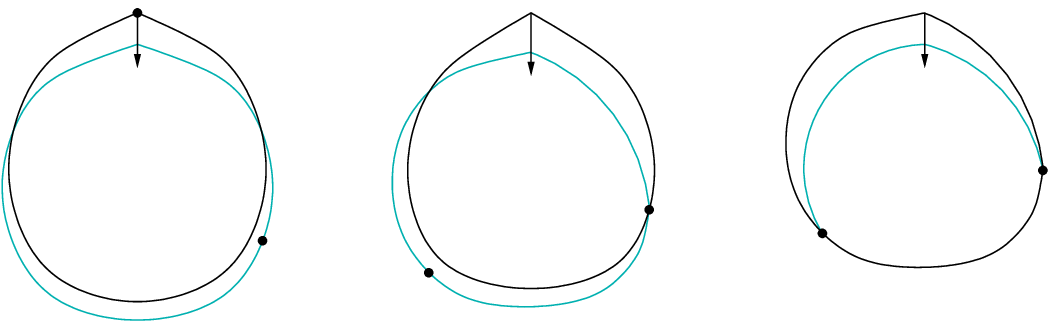}%
\end{picture}%
\setlength{\unitlength}{3315sp}%
\begingroup\makeatletter\ifx\SetFigFont\undefined%
\gdef\SetFigFont#1#2#3#4#5{%
  \reset@font\fontsize{#1}{#2pt}%
  \fontfamily{#3}\fontseries{#4}\fontshape{#5}%
  \selectfont}%
\fi\endgroup%
\begin{picture}(6020,1935)(1466,-1783)
\put(2296, 29){\makebox(0,0)[lb]{\smash{{\SetFigFont{8}{9.6}{\rmdefault}{\mddefault}{\updefault}{\color[rgb]{0,0,0}$i$}%
}}}}
\put(5221,-1276){\makebox(0,0)[lb]{\smash{{\SetFigFont{8}{9.6}{\rmdefault}{\mddefault}{\updefault}{\color[rgb]{0,0,0}$j_1$}%
}}}}
\put(4546, 29){\makebox(0,0)[lb]{\smash{{\SetFigFont{8}{9.6}{\rmdefault}{\mddefault}{\updefault}{\color[rgb]{0,0,0}$z$}%
}}}}
\put(6796, 29){\makebox(0,0)[lb]{\smash{{\SetFigFont{8}{9.6}{\rmdefault}{\mddefault}{\updefault}{\color[rgb]{0,0,0}$z$}%
}}}}
\put(3016,-1456){\makebox(0,0)[lb]{\smash{{\SetFigFont{8}{9.6}{\rmdefault}{\mddefault}{\updefault}{\color[rgb]{0,0,0}$j_1$}%
}}}}
\put(3736,-1636){\makebox(0,0)[lb]{\smash{{\SetFigFont{8}{9.6}{\rmdefault}{\mddefault}{\updefault}{\color[rgb]{0,0,0}$j_2$}%
}}}}
\put(6031,-1456){\makebox(0,0)[lb]{\smash{{\SetFigFont{8}{9.6}{\rmdefault}{\mddefault}{\updefault}{\color[rgb]{0,0,0}$j_2$}%
}}}}
\put(7471,-1006){\makebox(0,0)[lb]{\smash{{\SetFigFont{8}{9.6}{\rmdefault}{\mddefault}{\updefault}{\color[rgb]{0,0,0}$j_1$}%
}}}}
\put(2116,-376){\makebox(0,0)[lb]{\smash{{\SetFigFont{8}{9.6}{\rmdefault}{\mddefault}{\updefault}{\color[rgb]{0,0,0}$z$}%
}}}}
\end{picture}%

%% file: Fig/HypKites.tex
\begin{picture}(0,0)%
\includegraphics{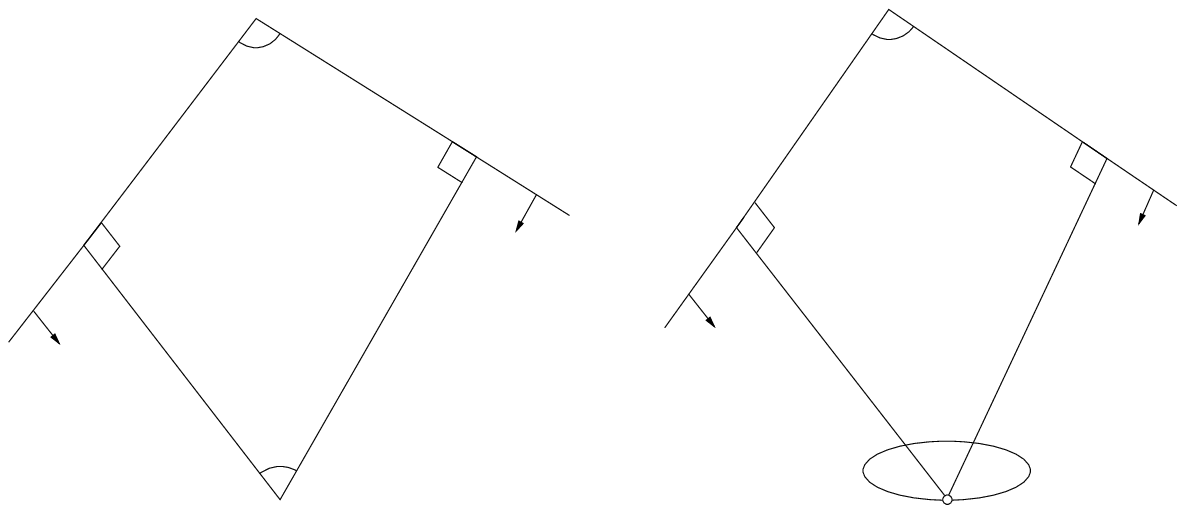}%
\end{picture}%
\setlength{\unitlength}{2072sp}%
\begingroup\makeatletter\ifx\SetFigFont\undefined%
\gdef\SetFigFont#1#2#3#4#5{%
  \reset@font\fontsize{#1}{#2pt}%
  \fontfamily{#3}\fontseries{#4}\fontshape{#5}%
  \selectfont}%
\fi\endgroup%
\begin{picture}(10842,4545)(948,-5957)
\put(8146,-2446){\makebox(0,0)[lb]{\smash{{\SetFigFont{8}{9.6}{\rmdefault}{\mddefault}{\updefault}{\color[rgb]{0,0,0}$\gamma$}%
}}}}
\put(4656,-4316){\makebox(0,0)[lb]{\smash{{\SetFigFont{8}{9.6}{\rmdefault}{\mddefault}{\updefault}{\color[rgb]{0,0,0}$c$}%
}}}}
\put(8401,-4806){\makebox(0,0)[lb]{\smash{{\SetFigFont{8}{9.6}{\rmdefault}{\mddefault}{\updefault}{\color[rgb]{0,0,0}$b$}%
}}}}
\put(2256,-4791){\makebox(0,0)[lb]{\smash{{\SetFigFont{8}{9.6}{\rmdefault}{\mddefault}{\updefault}{\color[rgb]{0,0,0}$b$}%
}}}}
\put(2161,-2676){\makebox(0,0)[lb]{\smash{{\SetFigFont{8}{9.6}{\rmdefault}{\mddefault}{\updefault}{\color[rgb]{0,0,0}$\gamma$}%
}}}}
\put(4531,-2071){\makebox(0,0)[lb]{\smash{{\SetFigFont{8}{9.6}{\rmdefault}{\mddefault}{\updefault}{\color[rgb]{0,0,0}$\beta$}%
}}}}
\put(10366,-2086){\makebox(0,0)[lb]{\smash{{\SetFigFont{8}{9.6}{\rmdefault}{\mddefault}{\updefault}{\color[rgb]{0,0,0}$\beta$}%
}}}}
\put(10636,-4216){\makebox(0,0)[lb]{\smash{{\SetFigFont{8}{9.6}{\rmdefault}{\mddefault}{\updefault}{\color[rgb]{0,0,0}$c$}%
}}}}
\put(10486,-5461){\makebox(0,0)[lb]{\smash{{\SetFigFont{8}{9.6}{\rmdefault}{\mddefault}{\updefault}{\color[rgb]{0,0,0}$H$}%
}}}}
\put(6925,-4442){\makebox(0,0)[lb]{\smash{{\SetFigFont{8}{9.6}{\rmdefault}{\mddefault}{\updefault}{\color[rgb]{0,0,0}$L_2$}%
}}}}
\put(6259,-3395){\makebox(0,0)[lb]{\smash{{\SetFigFont{8}{9.6}{\rmdefault}{\mddefault}{\updefault}{\color[rgb]{0,0,0}$L_1$}%
}}}}
\put(11775,-3345){\makebox(0,0)[lb]{\smash{{\SetFigFont{8}{9.6}{\rmdefault}{\mddefault}{\updefault}{\color[rgb]{0,0,0}$L_1$}%
}}}}
\put(963,-4615){\makebox(0,0)[lb]{\smash{{\SetFigFont{8}{9.6}{\rmdefault}{\mddefault}{\updefault}{\color[rgb]{0,0,0}$L_2$}%
}}}}
\put(3479,-5529){\makebox(0,0)[lb]{\smash{{\SetFigFont{8}{9.6}{\rmdefault}{\mddefault}{\updefault}{\color[rgb]{0,0,0}$\alpha$}%
}}}}
\put(3336,-1996){\makebox(0,0)[lb]{\smash{{\SetFigFont{8}{9.6}{\rmdefault}{\mddefault}{\updefault}{\color[rgb]{0,0,0}$a$}%
}}}}
\put(9131,-1921){\makebox(0,0)[lb]{\smash{{\SetFigFont{8}{9.6}{\rmdefault}{\mddefault}{\updefault}{\color[rgb]{0,0,0}$a$}%
}}}}
\end{picture}%